\documentclass[12pt,a4paper]{article}
\usepackage{bbm}
\usepackage{stmaryrd}
\usepackage{}
\usepackage{bbding}
\usepackage{mathrsfs}
\usepackage{amssymb}
\usepackage{amsfonts}
\usepackage{amsmath}
\usepackage{cases}
\usepackage{amsthm}
\usepackage[square, comma, sort&compress, numbers]{natbib}

\usepackage[top=1in,bottom=1in,left=1.20in,right=1.20in]{geometry}

\newtheorem{theorem}{Theorem}[section]
\newtheorem{lemma}{Lemma}[section]
\newtheorem{definition}{Definition}[section]
\newtheorem{proposition}[lemma]{Proposition}
\newtheorem{remark}{Remark}[section]

\newtheorem{corollary}{Corollary}[section]

\numberwithin{equation}{section}
\theoremstyle{plain}

\newcommand{\me}{\mathrm{e}}
\newcommand{\mi}{\mathrm{i}}

\makeatletter

\newcommand{\Rmnum}[1]{\expandafter\@slowromancap\romannumeral #1@}
\makeatother

\DeclareMathOperator{\mes}{mes}
\DeclareMathOperator{\card}{card}
\DeclareMathOperator{\sgn}{sgn}
\DeclareMathOperator{\diam}{diam}

\begin{document}
\title{Reducibility of linear quasi-periodic  Hamiltonian derivative wave equations and half-wave equations under  the Brjuno
conditions}
\author{Zhaowei Lou$^{1,2}$\footnote{ E-mail: zwlou@nuaa.edu.cn}
\\\\$^1$\small{School of Mathematics, Nanjing University of Aeronautics and Astronautics},\\
\small{Nanjing 211106, China.}
\\\\$^2$\small{MIIT Key Laboratory of Mathematical Modelling}\\
\small{and High Performance Computing of Air Vehicles}, \\
\small{Nanjing 211106,  China.}}
\date{}
\maketitle
\begin{abstract}
In this paper, we prove the reducibility for some linear quasi$-$periodic Hamiltonian derivative wave  and half$-$wave equations under the  Brjuno$-$R\"{u}ssmann non$-$resonance conditions.
 This generalizes  KAM theory by P\"{o}schel  in  $[38]$ from the finite dimensional Hamiltonian systems to Hamiltonian PDEs.
\end{abstract}

\textbf{Keywords.}  Reducibility, KAM, Brjuno-R\"{u}ssmann condition.

\textbf{2010 Mathematics Subject Classification:} 37K55, 35L05, 35Q55.

\section[]{Introduction and Main Result}
%

The reducibility theory of linear  quasi-periodic systems is the generalization of the classical Floquet theory for  linear  periodic systems.
It is important both  in the linear problems (spectral analysis of operator, growth of Sobolev norms) and  in  the non-linear case(linear stability analysis  of  quasi-periodic solutions of non-linear systems).
The first  reducibility result via  Kolmogorov-Arnold-Moser (KAM) theory was due to Bogoljubov,   Mitropoliskii and Samoilenko  \cite{Bogoljubov}, Dinaburg and Sinai  \cite{Dinaburg}  for finite
degrees of freedom systems.
Since then KAM theory has been a powerful tool to study reducibility theory.
In the late 1980s and early 1990s, KAM theory was extended to non-linear partial differential equations (PDEs) by Kuksin \cite{Kuksin87} and Wayne \cite{Wayne90}.
 See also \cite{Kuksin96, Poschel96a, Poschel96b} for further developments. As a corollary, these results  imply the reducibility of
the variational equations for quasi-periodic solutions of non-linear PDEs. In fact, ``reducibility is not only an important outcome of KAM but also
an essential ingredient in the proof" (\cite{Eliasson10}).

The first pure reducibility result for linear quasi-periodic PDEs was given by  Bambusi and Graffi \cite{Bambusi2001}.
They proved the reducibility of linear Schr\"{o}dinger  equations with unbounded perturbations.
  Eliasson and Kuksin    \cite{Eliasson-Kuksin2009} investigated the  reducibility of higher dimensional linear quasi-periodic {S}chr\"{o}dinger equations.
  Combining the pseudo-differential calculus, Baldi, Berti and Montalto \cite{Baldi141, Baldi142}
obtained the reducibility of quasi-linear  forced perturbations of Airy equation and
quasi-linear  KdV equation.
Thereafter,   these  results are developed and extended widely.
One  could refer to  \cite{Bambusi16, Bambusi17, Grebert162, Grebert11, Grebert14,  Liang17, Bambusi18,  Bambusi18Montalo, LiuJJ10,  Montalto,  BambusiLangellaMontalto19, BambusiLangellaMontalto22} and the references therein.

Consider a linear quasi-periodic PDE of the form
\begin{equation}\label{qppde}
\partial_tu=(A+P(\omega t))u,
\end{equation}
where $A$ is a positive self-adjoint operator and $P$ is a operator-valued function with the basic frequencies $\omega.$
It is well known that KAM reducibility requires a lower bound  on small divisors of the form
\begin{equation}\label{sd}
 | k\cdot \omega+ \lambda_i(\omega)-\lambda_i(\omega)|,
\end{equation}
where  $k\cdot \omega=\sum^n_{i=1} k_i \omega_i $ and $\{\lambda_i\}$ are  the  eigenvalues of the operator $A.$
In all the above-mentioned papers, the  lower bound of  Diophantine  type was used. Namely, the following non-resonance conditions holds:
 $ | k\cdot \omega+ \lambda_i(\omega)-\lambda_i(\omega)|\geq \frac{\gamma}{|k|^\tau},$
where the constants $\gamma>0,\,\, \tau>n-1.$
On the other hand, thanks to    the pioneering works  of Brjuno \cite{Brjuno}, the Diophantine  conditions  can be  weakened  to the Brjuno conditions.
To make it applicable in KAM scheme,  R\"{u}ssmann \cite{Russmann2, Russmann}    introduced the notion of an approximation function  to
characterize  the Brjuno conditions.
Under such Brjuno-R\"{u}ssmann type  conditions, P\"{o}schel \cite{Poschel89} proved the persistence of elliptic lower dimensional tori
in finite dimensional Hamiltonian systems. In \cite{Poschel90}, P\"{o}schel
also  proved the   existence of infinite dimensional invariant tori
in infinite dimensional Hamiltonian systems  of the form $H=\omega\cdot I+P(\theta,I).$
Later on, Xu and You \cite{XuYou}  and  Chavaudret and Marmi  \cite{Chavaudret} proved the reducibility of linear ODEs with almost periodic coefficients and   quasi-periodic cocycles   under such Brjuno-R\"{u}ssmann type  conditions, respectively.
See also   \cite{HuangLi, SiSi1, SiSi2}  for nonlinear forced  ODEs.
We also mention some   Brjuno type  quasi-periodic  results of  Corsi and Gentile \cite{CorsiGentile} and  Gentile \cite{Gentile} for forced non-Hamiltonian  ODEs  without using  approximation function.

To the best of our knowledge, there has been no   Brjuno-R\"{u}ssmann type  results in  KAM theory for PDEs.
In this paper, we establish a  reducibility theorem for some linear Hamiltonian  PDEs  under Brjuno-R\"{u}ssmann non-resonance  conditions. More precisely,
we consider the following  linear quasi-periodic derivative wave equations
\begin{equation}\label{eq1}
      \partial_{tt}u-\partial_{xx}u+mu+\epsilon V(\omega t,x) D_mu=0,\,\,x\in [0, \pi]
\end{equation}
and linear quasi-periodic half-wave equations
\begin{equation}\label{eq2}
     \mi  \partial_{t}u+D_0u+\epsilon V(\omega t,x)u=0,\,\,x\in [0, \pi],
\end{equation}
  under Dirichlet boundary conditions,
where  the first order pseudo-differential operators $D_m:=\sqrt{-\partial_{xx}+m},\,m\geq 0.$
The basic fequencies $\omega$ of the potential $V$ satisfy the Brjuno-R\"{u}ssmann non-resonance  conditions.
The wave equation \eqref{eq1} covers the variational equation around the quasi-periodic solutions  of nonlinear Hamiltonian derivative wave equation
$\partial_{tt}u-\partial_{xx}u+mu+ (D_mu)^3=0.$
Quasi-periodic solutions with  Diophantine frequencies of this nonlinear wave equation under  periodic boundary conditions
have been  obtained in \cite{Berti13}.
The half-wave equation \eqref{eq2} is an important class of PDEs arising  in various physical problems \cite{Cai2001, Kirkpatrick, Majda1997, Elgart, Froehlich2010}.
There are two  main difficulties when studying the  reducibility theory of the equations \eqref{eq1} and \eqref{eq2}.
The first one is the weak   dispersion relation since  the eigenvalues  $\lambda_j\sim  j,\,j\rightarrow \infty.$
The second  one is the bad  smoothness of the  perturbations.
To overcome this, we introduce a simplified version of T\"{o}plitz-Lipschitz functions and T\"{o}plitz-Lipschitz matrices, which were first proposed by Eliasson and Kuksin \cite{Eliasson10}  in KAM theory for  the  higher dimensional {S}chr\"{o}dinger equations.   Such simplified form is more suitable to the equations  \eqref{eq1} and  \eqref{eq2} and it was also used in \cite{Geng11, Geng13}.  Different from that in \cite{Geng11, Geng13},  we  characterize the T\"{o}plitz-Lipschitz functions in a  way  of semi-norm. We also mention  the quasi-T\"{o}plitz functions introduced in \cite {Berti13} for  nonlinear Hamiltonian derivative wave equations, which is also  an improved  version of  Eliasson-Kuksin's T\"{o}plitz-Lipschitz functions.  Comparing to  the quasi-T\"{o}plitz functions,  our  simplified form is more  easy to handle.

%

To state our main results, we introduce some definitions and assumptions on the potentials $V$ in the equations  \eqref{eq1} and  \eqref{eq2}.
\begin{definition}[Approximation function, \cite{Poschel89, Russmann}]
A non-decreasing function
$$\Delta: [0, \infty)\rightarrow  [1, \infty)$$
is called an approximation function, if

\begin{equation}\label{decay1}
 \frac{\Delta(t)}{t}\downarrow0,\,\,  0 \leq t\rightarrow\infty
\end{equation}
and
\begin{equation}\label{decay2}
 \int\limits^{\infty}_{1}\frac{\log\Delta(t)}{t}dt<\infty.
\end{equation}
in addition, the normalization  $\Delta(0)=1$ is imposed for definition.

\end{definition}

\begin{definition}[Brjuno-R\"{u}ssmann frequency]
Let  $\Delta$
be an approximation function. A vector $\omega\in \mathbb{R}^n$
is called Brjuno-R\"{u}ssmann frequency vector
if it
satisfies
\begin{equation}\label{brf}
 |k\cdot \omega|\geq \frac{\gamma}{\Delta(|k|)},
\end{equation}
for some constant $\gamma>0.$
\end{definition}

 $Assumption\, 1.$ Suppose the  function   $V:\mathbb{T}^n\times [0, \pi]\rightarrow \mathbb{R}$
is real analytic in $(\theta,x).$
For $\theta\in\mathbb{T}^n,$ $V(\theta, \cdot)$ is a $2\pi-$periodic, even function $V(\theta, x)=V(\theta, -x).$
Then it  can be written  as
\begin{equation}\label{V1}
  V(\theta, x)=\sum\limits_{j\geq 0} \widetilde{V}_{j}(\theta)\cos jx.
\end{equation}
Moreover, suppose  for all  $\theta,$ the function $V(\theta, \cdot)$ extends to a complex analytic function on a strip $|Im x|<2a$ for some $a>0.$
For all  $x,$ the function $V(\cdot, x)$ extends to a complex analytic function on a strip on $|Im \theta|<2r$ for some $r>0.$
Then there is a positive  constant $C_V>0$  such that  for  $p\geq 0,$
 \begin{equation}\label{V2}
  \|V\|_{D(2r), 2a, p}:=\|\widetilde{V}_{0}\|_{D(2r)}+ \sum\limits_{j\geq 1} j^{p} \me^{2aj} \|\widetilde{V}_{j}\|_{D(2r)}\leq C_V,
\end{equation}
where the norm $ \|\cdot\|_{D(2r)}$ is defined in Section 2.

Let\,$\phi_j(x)=\sqrt{\frac{2}{\pi}}\sin  jx,\,j\geq1$\, be the normalized Dirichlet eigenfunctions of the operator \,$D^2_m:=-\partial_{xx}+m$\, associated  to the eigenvalues $ \lambda_j^2=j^2+m,\,j\geq1.$\,
We consider the equations \eqref{eq1} and \eqref{eq2} in  the following  function  space
 \begin{equation}\label{fB1}
    \mathcal{H}^{a,p}_0=\left\{u=\sum\limits_{j\geq1}q_j\phi_j: \|u\|_{a,p}=\sum\limits_{j\geq1}j^p\me^{a j}|q_j|<\infty\right\}.
 \end{equation}

Our main result is stated as follows.
\begin{theorem}\label{maintheorem}
Under the Assumption 1 on the potential functions   $V,$   there is $\epsilon_0$ so that
for all $0<\epsilon<\epsilon_0$  there exists $\mathcal{O}_\epsilon \subseteq [0, 2\pi)^{n}$ of positive Lebesgue measure
such that for all $\omega\in \mathcal{O}_\epsilon $  satisfying Brjuno-R\"{u}ssmann non-resonance conditions,
the above linear quasi-periodic wave equation \eqref{eq1} and  half-wave equation \eqref{eq2}
reduce to the  linear equations  with constant coefficients with respect to the time variable.
\end{theorem}
In  Section 5, we prove  this theorem  by the reducibility  Theorem \ref{Redb}.

As a corollary of Theorem \ref{maintheorem}, we have the following conclusion  concerning the solutions of the equations \eqref{eq1} and  \eqref{eq2}:
\begin{corollary}\label{cor1}
Let the  initial data  $u_0\in \mathcal{H}^{a,p}_0,$ $v_0\in \mathcal{H}^{a,p-1 }_0.$
Under the Assumption 1, there is $\epsilon_0$ so that
for all $0<\epsilon<\epsilon_0$  and   $\omega\in \mathcal{O}_\epsilon $,

\emph{(i)}  there exists a unique solution $(u(t,x), u_t(t,x))\in \mathcal{H}^{a,p}_0\times \mathcal{H}^{a,p-1 }_0$
of  the  wave equation \eqref{eq1} with $(u(0,x), u_t(0,x))=(u_0, v_0).$
Moreover, $u(t,x)$ is almost-periodic in time and stable, i.e.,
$$\sup_{t\in \mathbb{R}}(\|u(t,\cdot)\|_{a,p}+\|u_t(t,\cdot)\|_{a,p-1})\leq C( \|u_0\|_{a,p}+\|v_0\|_{a,p-1})$$
for some constant $C=C(a,p,\omega)>0.$

\emph{(ii)}  there exists a unique solution $u(t,x)\in \mathcal{H}^{a,p}_0$
of  the  half-wave equation \eqref{eq2} with $u(0,x)=u_0.$
Moreover, $u(t,x)$ is almost-periodic in time and stable, i.e.,
$$\sup_{t\in \mathbb{R}}\|u(t,\cdot)\|_{a,p}\leq C \|u_0\|_{a,p}$$
for some constant $C=C(a,p,\omega)>0.$

\end{corollary}

\section[]{Functional setting}

 Let  $\mathcal{O} \subset \mathbb{R}^n$ be a parameter set of positive Lebesgue measure. Throughout the paper, for any real or complex valued function
depending on parameters $\xi\in \mathcal{O},$ its  derivatives with respect to $\xi$ are understood in the sense of Whitney.
We denote by  $C^1_W(\mathcal{O})$ the class of $C^1$ Whitney differentiable
 functions on $\mathcal{O}.$

Suppose $f\in C^1_W(\mathcal{O}),$ we define its norm as
\begin{equation*}
  \begin{split}
       \| f\|_{\mathcal{O}}:=&\sup\limits_{\xi\in \mathcal{O}}(|f(\xi)|+|\frac{\partial f}{\partial \xi}(\xi)|),\\
   \end{split}
\end{equation*}
where $|\cdot|$ denotes the sup-norm of complex vectors.

 Given an $n$-torus $\mathbb{T}^n=\mathbb{R}^n/{(2\pi \mathbb{Z})^n}$
and its complex neighborhood
$$ D(r)=\{\theta\in \mathbb{C}^n: |\hbox{Im} \theta|<r,\,\,r>0\}.$$

Consider a real analytic function $f(\theta;\xi)$
on  $\theta\in D(r).$  It is also
$C^1_W$ on $\xi\in \mathcal{O}.$ Its Fourier expansion reads
$f(\theta;\xi)=\sum\limits_{k\in \mathbb{Z}^n}\widehat{f}(k;\xi)\me^{\mi  k\cdot\theta},$ then
we define its norm as
 \begin{equation*}
       \| f\|_{D(r)\times\mathcal{O}}:=\sum\limits_{k\in \mathbb{Z}^n}|\widehat{f}(k;\cdot)|_{\mathcal{O}}\me^{|k|r},
\end{equation*}
where $k\cdot\theta=\sum\limits^{n}_{i=1}k_i\theta_i$ and
$|k|=\sum\limits^{n}_{i=1}|k_i|$.

Let $K>0.$ For  $f(\theta;\xi)$ above, its  $K-$order Fourier truncation $\mathcal{T}_Kf$ is defined as
$$(\mathcal{T}_Kf)(\theta):=\sum\limits_{k\in \mathbb{Z}^n,\,|k|< K}\widehat{f}(k)\me^{\mi k\cdot\theta}.$$
The remainder  $\mathcal{R}_Kf$ of $f$ is defined by
$(\mathcal{R}_Kf)(\theta):=f(\theta)-\mathcal{T}_Kf(\theta).$
Suppose $0<2\sigma<r,$ we have the following estimate for $\mathcal{R}_Kf:$
\begin{equation}\label{remaind}
    \|\mathcal{R}_Kf\|_{D(r-2\sigma)\times\mathcal{O}}\leq 32\sigma^{-2}\me^{-K\sigma}\| f\|_{D(r)\times\mathcal{O}}.
\end{equation}

The average $[f]$ of  $f$ on $\mathbb{T}^n$ is defined as
 $$[f]:=\widehat{f}(0)=(2\pi)^{-n}\int_{\mathbb{T}^n}f(\theta)d\theta.$$

Let $a, p> 0$, we introduce the Banach space $\ell^{a,p}$ of all real or complex sequences
$z=(z_j)_{j\in \mathbb{Z}}$ with
$$\|z\|_{a,p}=\sum\limits_{j\in \mathbb{Z}}\me^{aj}j^{p}|z_j|<\infty.$$

Given  $r, s>0,$ we define the  phase space $$\mathcal{P}^{a,p}:=\mathbb{T}^n \times\mathbb{R}^n\times \ell^{a,p}\times \ell^{a,p}\ni w:=(\theta, I,  z,  \bar{z})$$
and a complex neighborhood
\begin{equation*}
  D(r,s)=\{w:|\hbox{Im} \theta|<r,|I|<s^2,
   \|z\|_{a,p}<s, \|\bar{z}\|_{a,p}<s\}
\end{equation*}
of $\mathcal{T}^{n}_0:=\mathbb{T}^{n}\times\{I=0\}\times \{z=0\}\times\{\bar{z}=0\}$ in $\mathcal{P}^{a,p}_{\mathbb{C}}:=\mathbb{C}^n \times\mathbb{C}^n\times \ell^{a,p}\times \ell^{a,p}.$

Consider a real analytic function $f(\theta,I, z,\bar{z}; \xi)$
on  $D(r,s),$  which  is also
$C^1_W$ on $\xi\in \mathcal{O}.$ Its Taylor-Fourier expansion reads
\begin{equation*}
  \begin{split}
        f(\theta,I, z,\bar{z}; \xi)=\sum\limits_{l,\alpha, \beta}f_{l\alpha \beta}(\theta;\xi)I^l z^{\alpha}\bar{z}^{\beta}
            =\sum\limits_{k\in \mathbb{Z}^n,l,   \alpha, \beta}\widehat{f}_{l\alpha \beta}(k;\xi)\me^{\mi k\cdot\theta}I^l z^{\alpha}\bar{z}^{\beta} ,
   \end{split}
\end{equation*}
where  we use the multi-index notations  $l=(l_j)^n_{j=1},$ $\alpha=(\alpha_j)_{j\geq1},$ $\beta=(\beta_j)_{j\geq1}$ with  $l_j, \alpha_j,\beta_j\in \mathbb{N}.$
 $\alpha $ and $\beta$ have only finitely many nonzero components.
 $I^lz^{\alpha}\bar{z}^{\beta}=(\prod\limits^n_{i=1}I^{l_i}_i )(\prod\limits_{j\in \mathbb{Z}}z^{\alpha_j}_j\bar{z}^{\beta_j}_j).$

 We define the majorant of $f$  as
 \begin{equation*}
  \begin{split}
        \lfloor f\rceil_{D(r)\times\mathcal {O}}\equiv\lfloor f(\cdot, I, z, \bar{z};\cdot)\rceil_{D(r)\times \mathcal{O}}:=\sum\limits_{l,\alpha, \beta}\|f_{l\alpha \beta}\|_{D(r)\times \mathcal{O}}|I^{l}||z^{\alpha}||\bar{z}^{\beta}|
   \end{split}
\end{equation*}
 and the  norm of $f$ as
\begin{equation*}
  \begin{split}
        \|f\|_{D(r,s)\times\mathcal {O}}:=&\sup\limits_{|I|<s^2, \|z\|_{a,p}<s,\atop \|\bar{z}\|_{a,p}<s} \lfloor f\rceil_{D(r)\times \mathcal{O}}\\
        =&\sup\limits_{|I|<s^2,\|z\|_{a,p}<s,\atop \|\bar{z}\|_{a,p}<s} \sum\limits_{l,\alpha, \beta}\|f_{l\alpha \beta}\|_{D(r)\times \mathcal{O}}|I^{l}||z^{\alpha}||\bar{z}^{\beta}|.\\
   \end{split}
\end{equation*}

Consider an infinite dimensional  dynamical system on $D(r, s):$
\begin{equation*}
    \dot{w}=X(w),\,\,w=(\theta, I, z,\bar{z})\in D(r, s),
\end{equation*}
where the vector field
\begin{equation*}
  \begin{split}
        X(w)=&(X^{(\theta)}(w),X^{(I)}(w),X^{(z)}(w), X^{(\bar{z})}(w)),
   \end{split}
\end{equation*}

 Suppose vector field $X(w;\xi)$ is real analytic on $D(r, s)$  and  $C^1_W$ smooth on  $\mathcal {O},$
we define the weighted norm of $X$ as follows
\begin{equation*}
 \begin{split}
\|X\|_{s;D(r,s)\times\mathcal {O}}&\\
=\sup\limits_{|I|<s^2,\atop \|z\|_{a,p}<s, \|\bar{z}\|_{a,p}<s}\Big\{ & \sum\limits^n_{i=1}\lfloor X^{(\theta_i)}\rceil_{D(r)\times\mathcal {O}}+\frac{1}{s^2}\sum\limits^n_{i=1}\lfloor X^{(I_i)}\rceil_{D(r)\times\mathcal {O}}\\
&+\frac{1}{s}\sum\limits_{j\in\mathbb{ Z}}\me^{aj}j^p(\lfloor X^{(z_j)}\rceil_{D(r)\times\mathcal {O}}
+\lfloor X^{(\bar{z}_j)}\rceil_{D(r)\times\mathcal {O}})\Big\}.
\end{split}
\end{equation*}

\section[]{T\"{o}plitz-Lipschitz Functions}\label{tlf}

 \subsection{Definitions}\label{vecpoldedef}
In this section, we introduce a class  of real analytic functions with exponentially off-diagonal decay.

\begin{definition}
Let  $r,  s, \rho>0.$   Suppose  $P(\theta,z,\bar{z}; \xi)$   is  real analytic  on $(\theta,z,\bar{z})\in D(r,s)$ and $C^1_W-$smooth on parameters  $\xi\in \mathcal{O}.$
We say that  $P$   is T\"{o}plitz-Lipschitz and write $P\in \mathcal{T}^\rho_{ D(r,s)\times\mathcal{O}}$ if
\begin{equation}\label{DecaySp}
        \langle P\rangle_{\rho,D(r,s)\times\mathcal{O}}<\infty,
\end{equation}
where
the semi-norm $ \langle P\rangle_{\rho,D(r,s)}$ is defined by  the following conditions
\begin{description}
  \item[(T1) Exponentially off-diagonal decay.]
  \begin{equation}\label{tl11}
\left\|\frac{\partial^2 P}{\partial z_{i}\partial z_{j}}\right\|_{D(r,s)\times\mathcal{O}}
\leq \langle P\rangle_{\rho,D(r,s)\times \mathcal{O}}\me^{-\rho|i+j|}.
\end{equation}
\begin{equation}\label{tl12}
\left\|\frac{\partial^2 P}{\partial z_{i}\partial \bar{z}_{j}}\right\|_{D(r,s)\times\mathcal{O}}
\leq \langle P\rangle_{\rho,D(r,s)\times \mathcal{O}}\me^{-\rho|i-j|}.
\end{equation}
\begin{equation}\label{tl13}
\left\|\frac{\partial^2 P}{\partial \bar{z}_{i}\partial \bar{z}_{j}}\right\|_{D(r,s)\times\mathcal{O}}
\leq \langle P\rangle_{\rho,D(r,s)\times \mathcal{O}}\me^{-\rho|i+j|}.
\end{equation}
  \item[(T2) Asymptotically T\"{o}plitz.] The limits
  $$\lim_{t\in\mathbb{Z},\,t\rightarrow\infty}\frac{\partial^2 P}{\partial z_{i+t}\partial z_{j- t}},\,\,\lim_{t\in\mathbb{Z},\,t\rightarrow\infty}\frac{\partial^2 P}{\partial z_{i+t}\partial \bar{z}_{j+ t}}\,\, \hbox{and}\,\, \lim_{t\in\mathbb{Z},\,t\rightarrow\infty}\frac{\partial^2 P}{\partial \bar{z}_{i+t}\partial \bar{z}_{j- t}}$$
  exist and are finite for all $i,j\in\mathbb{Z}.$
  \item[(T3) Lipschitz at infinity.] For sufficiently large $|t|, t\in\mathbb{Z},$ the following   hold.
\begin{equation}\label{tl21}
\left\|\frac{\partial^2 P}{\partial z_{i+t}\partial z_{j- t}}-\lim_{t\rightarrow\infty}\frac{\partial^2 P}{\partial z_{i+t}\partial z_{j- t}}\right\|_{D(r,s)\times\mathcal{O}}
\leq |t|^{-1}\langle P\rangle_{\rho,D(r,s)\times \mathcal{O}}\me^{-\rho|i+j|}.,
\end{equation}
\begin{equation}\label{tl22}
\left\|\frac{\partial^2 P}{\partial z_{i+t}\partial \bar{z}_{j+ t}}-\lim_{t\rightarrow\infty}\frac{\partial^2 P}{\partial z_{i+t}\partial \bar{z}_{j+ t}}\right\|_{D(r,s)\times\mathcal{O}}
\leq |t|^{-1}\langle P\rangle_{\rho,D(r,s)\times \mathcal{O}}\me^{-\rho|i-j|},
\end{equation}
\begin{equation}\label{tl23}
\left\|\frac{\partial^2 P}{\partial \bar{z}_{i+t}\partial \bar{z}_{j- t}}-\lim_{t\rightarrow\infty}\frac{\partial^2 P}{\partial \bar{z}_{i+t}\partial \bar{z}_{j- t}}\right\|_{D(r,s)\times\mathcal{O}}\leq |t|^{-1}\langle P\rangle_{\rho,D(r,s)\times \mathcal{O}}\me^{-\rho|i+j|},
\end{equation}
\end{description}

\end{definition}


\begin{remark}
In  fact, $\langle P\rangle_{\rho,D(r,s)\times\mathcal{O}}$ is the smallest non-negative real number   that satisfies the above conditions
$(T1)-(T3).$ Moreover, it satisfies
\begin{itemize}
  \item  $\langle P\rangle_{\rho,D(r,s)\times \mathcal{O}}\geq 0;$
  \item $\langle \lambda P\rangle_{\rho,D(r,s)\times \mathcal{O}}=|\lambda|\langle P\rangle_{\rho,D(r,s)\times \mathcal{O}}$ for all $\lambda\in \mathbb{C};$
  \item $$\langle P+F\rangle_{\rho,D(r,s)\times \mathcal{O}}\leq\langle P\rangle_{\rho,D(r,s)\times \mathcal{O}}+ \langle F\rangle_{\rho,D(r,s)\times \mathcal{O}}.$$
\end{itemize}
Note that  $\langle P\rangle_{\rho,D(r,s)\times \mathcal{O}}= 0$  could  not imply
  $P=0.$ This means  $\langle \cdot\rangle_{\rho,D(r,s)\times \mathcal{O}}$ is only a semi-norm.
\end{remark}

\begin{remark}
From (T1) and (T3), the limits in (T3) satisfy
\begin{equation}\label{tl110}
\left\|\lim_{t\rightarrow\infty}\frac{\partial^2 P}{\partial z_{i+t}\partial z_{j- t}}\right\|_{D(r,s)\times\mathcal{O}}
\leq \langle P\rangle_{\rho,D(r,s)\times \mathcal{O}}\me^{-\rho|i+j|};
\end{equation}
\begin{equation}\label{tl120}
\left\|\lim_{t\rightarrow\infty}\frac{\partial^2 P}{\partial z_{i+t}\partial \bar{z}_{j+ t}}\right\|_{D(r,s)\times\mathcal{O}}
\leq \langle P\rangle_{\rho,D(r,s)\times \mathcal{O}}\me^{-\rho|i-j|};
\end{equation}
\begin{equation}\label{tl130}
\left\|\lim_{t\rightarrow\infty}\frac{\partial^2 P}{\partial \bar{z}_{i+t}\partial \bar{z}_{j- t}}\right\|_{D(r,s)\times\mathcal{O}}
\leq \langle P\rangle_{\rho,D(r,s)\times \mathcal{O}}\me^{-\rho|i+j|}.
\end{equation}
\end{remark}

 \begin{remark}
By the definition of the semi-norm $\langle \cdot\rangle_{\rho,D(r,s)\times \mathcal{O}}$, it is not difficulty to verify that
 the following conclusions hold:
 \begin{description}
      \item[(1)]
    $\langle P\rangle_{\rho,D(r',s')\times \mathcal{O}}\leq  \langle P\rangle_{\rho,D(r,s)\times \mathcal{O}}$ if $0<r'\leq r,\, 0<s'\leq s;$
    \item[(2)]
      $\langle P\rangle_{\rho',D(r,s)\times \mathcal{O}}\leq  \langle P\rangle_{\rho,D(r,s)\times \mathcal{O}}$ if $0< \rho'\leq \rho;$
       \item[(3)] Let $K>0,$  then  the  Fourier truncation  $\mathcal{T}_KP$ of $P$ satisfies
      $$\langle \mathcal{T}_KP\rangle_{\rho,D(r,s)\times \mathcal{O}}\leq  \langle P\rangle_{\rho,D(r,s)\times \mathcal{O}}$$
      and the remainder $\mathcal{R}_KP$ of $P$ satisfies
      $$\langle \mathcal{R}_KP\rangle_{\rho,D(r',s)\times \mathcal{O}}\leq  \me^{-K(r-r')}\langle P\rangle_{\rho,D(r,s)\times \mathcal{O}}$$
         if $0< r'\leq r.$
 \end{description}
 \end{remark}

\begin{definition}
Let $\ell^{a,p}_{0}$ be the unilateral infinite sequences space  defined by
\begin{equation}\label{uniseq}
\ell^{a,p}_{0}=\{z=(z_j)_{j\geq1}: \|z\|_{a,p}=\sum\limits_{j\geq1}|z_j||j|^p\me^{a|j|}<\infty\}.
\end{equation}
Given a  real analytic function $P(\theta, z, \bar{z})$ with   $(z, \bar{z})\in  \ell^{a,p}_{0}\times\ell^{a,p}_{0},$   we
lift it from  $\ell^{a,p}_{0}\times\ell^{a,p}_{0}$ to $\ell^{a,p}\times\ell^{a,p} $ by
$\widetilde{P}(\theta, \tilde{z}, \bar{\tilde{z}})= P(\theta, z, \bar{z}),$
where $(\tilde{z}, \bar{\tilde{z}})\in \ell^{a,p}\times\ell^{a,p}$ and $\tilde{z}=z_j, \bar{\tilde{z}}=\bar{z}_j$ for all $ j\geq1.$

We say that the function  $P$   is T\"{o}plitz-Lipschitz and write $P\in \mathcal{T}^\rho_{ D(r,s)\times\mathcal{O}}$ if
$\widetilde{P}(\theta, \tilde{z}, \bar{\tilde{z}})$ is T\"{o}plitz-Lipschitz and define
\begin{equation}\label{DecaySp2}
      \langle P\rangle_{\rho,D(r,s)\times\mathcal{O}} := \langle \widetilde{P}\rangle_{\rho,D(r,s)\times\mathcal{O}}<\infty.
\end{equation}
\end{definition}
%
%
%
%



Below we focus on  a class of quadratic functions on $(z, \bar{z})$ of the form
$$P(\theta,z, \bar{z}; \xi)=\sum\limits_{|\alpha|+|\beta|=2}P_{\alpha \beta}(\theta; \xi)z^{\alpha} \bar{z}^{\beta}.$$
We study the T\"{o}plitz-Lipschitz  property for these functions under the action of the Poisson bracket, the flow of linear Hamiltonian system and the canonical transformation.

\begin{proposition}[\textbf{Poisson bracket}]\label{Poissonbracket}
Let $0<\delta<\rho$ and  $0<\sigma<r.$  Suppose  the quadratic functions   $R,\,F\in \mathcal{T}^\rho_{ D(r,s)\times\mathcal{O}},$
then $\{R,F\}\in \mathcal{T}^{\rho-\delta}_{ D(r,s)\times\mathcal{O}}$ and
 there exists a constant $C>0$ so that
\begin{equation}\label{pobest}
\langle \{R,F\}\rangle_{\rho-\delta,D(r,s)\times \mathcal{O}}\leq  \frac{C}{\delta}\langle R\rangle_{\rho,D(r,s)\times \mathcal{O}}\langle F\rangle_{\rho,D(r,s)\times \mathcal{O}}.
\end{equation}
\end{proposition}

\begin{proof}

The Poisson bracket $\{R, F\}$ reads
\begin{equation*}
 \{R, F\}=\mi \sum\limits_{k\in \mathbb{Z}}\left( \frac{\partial R}{\partial z_k}  \frac{\partial F}{\partial \bar{z}_k} -  \frac{\partial R}{\partial \bar{z}_k}  \frac{\partial F}{\partial z_k} \right).
\end{equation*}

In what follows, it remains to analysis  the second derivative  $\frac{\partial^2 \{R, F\}}{\partial z_i \partial \bar{z}_j}$ with respect to  $z_i , \bar{z}_j,$
and  the  other second derivatives  could be  similarly done.

Since  the functions $R$ and $ F$  are both  quadratic  on $(z,\bar{z}),$ their third derivatives vanish.
Then we have
\begin{equation}\label{pb2d}
\begin{split}
 &\frac{\partial^2 \{R, F\}}{\partial z_i \partial \bar{z}_j}\\
 =&\sum\limits_{k\in \mathbb{Z}}\mi \left(\frac{\partial^2 R}{\partial z_k \partial \bar{z}_j} \frac{\partial^2 F}{\partial z_i \partial \bar{z}_k}
 + \frac{\partial^2 R}{\partial z_i \partial z_k} \frac{\partial^2 F}{\partial \bar{z}_k \partial \bar{z}_j}
 - \frac{\partial^2 R}{\partial \bar{z}_k \partial \bar{z}_j} \frac{\partial^2 F}{\partial z_i \partial z_k}
  - \frac{\partial^2 R}{\partial \bar{z}_k \partial z_i} \frac{\partial^2 F}{\partial z_k \partial \bar{z}_j}
 \right).\\
  \end{split}
\end{equation}

$\bullet$ We first verify the property (T1) for $\frac{\partial^2 \{R, F\}}{\partial z_i \partial \bar{z}_j}.$
It suffices  to  consider the sums
$\sum\limits_{k\geq1}\frac{\partial^2 R}{\partial z_k \partial \bar{z}_j} \frac{\partial^2 F}{\partial z_i \partial \bar{z}_k}$
and
  $\sum\limits_{k\geq1}\frac{\partial^2 R}{\partial z_i \partial z_k} \frac{\partial^2 F}{\partial \bar{z}_k \partial \bar{z}_j}$  in \eqref{pb2d}, and   the others can be similarly done.

Since  the functions $R$ and $ F$  satisfy  the property (T1), then we have
\begin{equation*}
\begin{split}
\left\|\sum\limits_{k}\frac{\partial^2 R}{\partial z_{k} \partial \bar{z}_{j}} \frac{\partial^2 F}{\partial z_{i} \partial \bar{z}_{k}}\right\|_{D{(r,s)}\times\mathcal{O}}
\leq&\sum\limits_{k}\left\|\frac{\partial^2 R}{\partial z_{k} \partial \bar{z}_{j}}\right\|_{D{(r,s)}\times\mathcal{O}}
\left\|\frac{\partial^2 F}{\partial z_{i} \partial \bar{z}_{k}}\right\|_{D{(r,s)}\times\mathcal{O}} \\
\leq&  \langle R\rangle_{\rho,D(r,s)\times \mathcal{O}}\langle F\rangle_{\rho,D(r,s)\times \mathcal{O}}\sum\limits_{k} \me^{-\rho(|i-k|+|k-j|)}\\
\leq&  \langle R\rangle_{\rho,D(r,s)\times \mathcal{O}}\langle F\rangle_{\rho,D(r,s)\times \mathcal{O}}\me^{-(\rho-\delta)(|i-j|)}\sum\limits_{k} \me^{-\delta(|i-k|+|k-j|)}\\
\leq&  C\delta^{-1} \langle R\rangle_{\rho,D(r,s)\times \mathcal{O}}\langle F\rangle_{\rho,D(r,s)\times \mathcal{O}}\me^{-(\rho-\delta)(|i-j|)}
 \end{split}
\end{equation*}
and
\begin{equation*}
\begin{split}
\left\|\sum\limits_{k}\frac{\partial^2 R}{\partial z_{i} \partial z_{k}} \frac{\partial^2 F}{\partial \bar{z}_{j} \partial \bar{z}_{k}}\right\|_{D{(r,s)}\times\mathcal{O}}
\leq&\sum\limits_{k}\left\|\frac{\partial^2 R}{\partial z_{i} \partial z_{k}}\right\|_{D{(r,s)}\times\mathcal{O}}
\left\|\frac{\partial^2 F}{\partial \bar{z}_{j} \partial \bar{z}_{k}}\right\|_{D{(r,s)}\times\mathcal{O}} \\
\leq&  \langle R\rangle_{\rho,D(r,s)\times \mathcal{O}}\langle F\rangle_{\rho,D(r,s)\times \mathcal{O}}\sum\limits_{k} \me^{-\rho(|i+k|+|k+j|)}\\
\leq&  C\delta^{-1} \langle R\rangle_{\rho,D(r,s)\times \mathcal{O}}\langle F\rangle_{\rho,D(r,s)\times \mathcal{O}}\me^{-(\rho-\delta)(|i-j|)}.
 \end{split}
\end{equation*}
here we use the inequality $\sum\limits_{k}\me^{-\delta(|i-k|+|k-j|)}\leq C\delta^{-1} $ ( see Lemma \ref{expineq}, Appendix).

$\bullet$ Verifying the  property (T2)  for $\frac{\partial^2 \{R, F\}}{\partial z_i \partial \bar{z}_j}.$
From the above analysis, we know that the functional series
$\sum\limits_{k\geq1}\frac{\partial^2 R}{\partial z_k \partial \bar{z}_j} \frac{\partial^2 F}{\partial z_i \partial \bar{z}_k}$
and
  $\sum\limits_{k\geq1}\frac{\partial^2 R}{\partial z_i \partial z_k} \frac{\partial^2 F}{\partial \bar{z}_k \partial \bar{z}_j}$
converge uniformly on $D(r,s)\times \mathcal{O}.$ Since the limits
$ \lim\limits_{t\rightarrow\infty} \frac{\partial^2 P}{\partial z_{i+t} \partial \bar{z}_{j+t}},$
 $\lim\limits_{t\rightarrow\infty} \frac{\partial^2 P}{\partial z_{i+t} \partial z_{j-t}}$ and
$ \lim\limits_{t\rightarrow\infty} \frac{\partial^2 P}{\partial \bar{z}_{i+t} \partial z_{j-t}}$ exist  and are finite,  then the limits
$$\lim\limits_{t\rightarrow\infty}\sum\limits_{k}\frac{\partial^2 R}{\partial z_{k+t} \partial \bar{z}_{j+t}} \frac{\partial^2 F}{\partial z_{i+t} \partial \bar{z}_{k+t}}$$
and
$$\lim\limits_{t\rightarrow\infty}\sum\limits_{k}\frac{\partial^2 R}{\partial z_{i+t} \partial z_{k-t}} \frac{\partial^2 F}{\partial \bar{z}_{j+t} \partial \bar{z}_{k-t}}$$
also exist  and are finite.
This implies the property (T2) holds  for $\frac{\partial^2 \{R, F\}}{\partial z_i \partial \bar{z}_j}.$

$\bullet$ Finally, we verify  the property  (T3) for $\frac{\partial^2 \{R, F\}}{\partial z_i \partial \bar{z}_j}.$
For the sake of convenience, we introduce  the notations
 $$ P^{11}_{ij,\infty}:=\lim\limits_{t\rightarrow\infty} \frac{\partial^2 P}{\partial z_{i+t} \partial \bar{z}_{j+t}},$$
 $$ P^{20}_{ij,\infty}:=\lim\limits_{t\rightarrow\infty} \frac{\partial^2 P}{\partial z_{i+t} \partial z_{j-t}}$$
 and
$$ P^{02}_{ij,\infty}:=\lim\limits_{t\rightarrow\infty} \frac{\partial^2 P}{\partial \bar{z}_{i+t} \partial z_{j-t}}.$$
In view of   $R,\,F\in \mathcal{T}^\rho_{ D(r,s)\times\mathcal{O}}$  and  thanks to  the difference equality
\begin{equation}\label{difference}
  AB-ab=(A-a)b+a(B-b)+(A-a)(B-b),
\end{equation}
and the inequality in Lemma \ref{expineq},
we have
\begin{equation*}
\begin{split}
&\left\|\sum\limits_{k}\frac{\partial^2 R}{\partial z_{k+t} \partial \bar{z}_{j+t}} \frac{\partial^2 F}{\partial z_{i+t} \partial \bar{z}_{k+t}}-\lim\limits_{t\rightarrow\infty}\sum\limits_{k}\frac{\partial^2 R}{\partial z_{k+t} \partial \bar{z}_{j+t}} \frac{\partial^2 F}{\partial z_{i+t} \partial \bar{z}_{k+t}}\right\|_{D{(r,s)}\times\mathcal{O}}\\
\leq&\sum\limits_{k}\left\|\frac{\partial^2 R}{\partial z_{k+t} \partial \bar{z}_{j+t}}-R^{11}_{kj,\infty}\right\|_{D{(r,s)}\times\mathcal{O}}
\left\|F^{11}_{ik,\infty}\right\|_{D{(r,s)}\times\mathcal{O}} \\
&+\sum\limits_{k}\left\|R^{11}_{kj,\infty}\right\|_{D{(r,s)}\times\mathcal{O}}
\left\|\frac{\partial^2 F}{\partial z_{i+t} \partial \bar{z}_{k+t}}-F^{11}_{ik,\infty}\right\|_{D{(r,s)}\times\mathcal{O}} \\
&+\sum\limits_{k}\left\|\frac{\partial^2 R}{\partial z_{k+t} \partial \bar{z}_{j+t}}-R^{11}_{kj,\infty}\right\|_{D{(r,s)}\times\mathcal{O}}
\left\|\frac{\partial^2 F}{\partial z_{i+t} \partial \bar{z}_{k+t}}-F^{11}_{ik,\infty}\right\|_{D{(r,s)}\times\mathcal{O}}\\
\leq& |t|^{-1} \langle R\rangle_{\rho,D(r,s)\times \mathcal{O}}\langle F\rangle_{\rho,D(r,s)\times \mathcal{O}}\sum\limits_{k} \me^{-\rho(|i-k|+|k-j|)}\\
& + |t|^{-1} \langle R\rangle_{\rho,D(r,s)\times \mathcal{O}}\langle F\rangle_{\rho,D(r,s)\times \mathcal{O}}\sum\limits_{k} \me^{-\rho(|i-k|+|k-j|)}\\
& +|t|^{-2}\langle R\rangle_{\rho,D(r,s)\times \mathcal{O}}\langle F\rangle_{\rho,D(r,s)\times \mathcal{O}}\sum\limits_{k} \me^{-\rho(|i-k|+|k-j|)}\\
\leq&|t|^{-1}C\delta^{-1} \langle R\rangle_{\rho,D(r,s)\times \mathcal{O}}\langle F\rangle_{\rho,D(r,s)\times \mathcal{O}}\me^{-(\rho-\delta)(|i-j|)}
 \end{split}
\end{equation*}
and
\begin{equation*}
\begin{split}
&\left\|\sum\limits_{k}\frac{\partial^2 R}{\partial z_{i+t} \partial z_{k-t}} \frac{\partial^2 F}{\partial \bar{z}_{j+t} \partial \bar{z}_{k-t}}-\lim\limits_{t\rightarrow\infty}\sum\limits_{k}\frac{\partial^2 R}{\partial z_{i+t} \partial z_{k-t}} \frac{\partial^2 F}{\partial \bar{z}_{j+t} \partial \bar{z}_{k-t}}\right\|_{D{(r,s)}\times\mathcal{O}}\\
\leq&\sum\limits_{k}\left\|\frac{\partial^2 R}{\partial z_{i+t} \partial z_{k-t}}-R^{20}_{ik,\infty}\right\|_{D{(r,s)}\times\mathcal{O}}
\left\|F^{02}_{jk,\infty}\right\|_{D{(r,s)}\times\mathcal{O}} \\
&+\sum\limits_{k}\left\|R^{20}_{ik,\infty}\right\|_{D{(r,s)}\times\mathcal{O}}
\left\|\frac{\partial^2 F}{\partial \bar{z}_{j+t} \partial \bar{z}_{k-t}}-F^{02}_{jk,\infty}\right\|_{D{(r,s)}\times\mathcal{O}} \\
&+\sum\limits_{k}\left\|\frac{\partial^2 R}{\partial z_{i+t} \partial z_{k-t}}-R^{20}_{ik,\infty}\right\|_{D{(r,s)}\times\mathcal{O}}
\left\|\frac{\partial^2 F}{\partial \bar{z}_{j+t} \partial \bar{z}_{k-t}}-F^{02}_{jk,\infty}\right\|_{D{(r,s)}\times\mathcal{O}} \\
\leq&|t|^{-1} C\delta^{-1}\langle R\rangle_{\rho,D(r,s)\times \mathcal{O}}\langle F\rangle_{\rho,D(r,s)\times \mathcal{O}}\me^{-(\rho-\delta)(|i-j|)}.
 \end{split}
\end{equation*}
These imply that
\begin{equation*}
\begin{split}
&\left\|\frac{\partial^2 \{R, F\}}{\partial z_{i+t} \partial \bar{z}_{j+t}}-\lim\limits_{t\rightarrow\infty}\frac{\partial^2 \{R, F\}}{\partial z_{i+t} \partial \bar{z}_{j+t}}\right\|_{D{(r,s)}\times\mathcal{O}}\\
\leq &|t|^{-1}C\delta^{-1}\langle R\rangle_{\rho,D(r,s)\times \mathcal{O}}\langle F\rangle_{\rho,D(r,s)\times \mathcal{O}}\me^{-(\rho-\delta)(|i-j|)}.
 \end{split}
\end{equation*}

\end{proof}

\subsection{T\"{o}plitz-Lipschitz Matrices}

Denote by $\mathcal{M}_{2}(\mathbb{C})$ the space of $2\times 2$ complex matrices.
Let $\|\cdot\|$ be any sub-multiplicative norm on $\mathcal{M}_{2}(\mathbb{C}).$
Consider  a bilateral infinite dimensional $\mathcal{M}_{2}(\mathbb{C})-$valued matrix
 $$A:\mathbb{Z}\times \mathbb{Z}\rightarrow \mathcal{M}_{2}(\mathbb{C}):$$
 $$(i,j)\mapsto A_{ij}=\left(
                         \begin{array}{cc}
                           A^{11}_{ij} & A^{12}_{ij} \\
                           A^{21}_{ij} & A^{22}_{ij} \\
                         \end{array}
                       \right).
 $$
The matrix multiplication   is defined by
 $(AB)_{ij}=\sum_{k\in \mathbb{Z}} A_{ik}B_{kj}.$

Now we consider the matrices depend on  $(\theta, \xi)\in D(r)\times \mathcal{O}.$
\begin{definition}[\textbf{Matrices with T\"{o}plitz-Lipschitz property}]
Let  $r, \rho>0.$   We say that a matrix  $A=A(\theta, \xi)$ on $D(r)\times \mathcal{O}$ is T\"{o}plitz-Lipschitz and  write  $A\in \mathfrak{M}^{\rho}_{r,\mathcal{O}}$ if $\langle\langle A\rangle\rangle_{\rho,r,\mathcal{O}}<\infty, $ where the norm $\langle\langle A\rangle\rangle_{\rho,r,\mathcal{O}}$
is defined  by the following conditions:
\begin{description}
  \item[(T1$^{\prime}$) Exponentially off-diagonal decay]
  \begin{equation}\label{mtl11}
\| A^{11}_{ij}\|_{D(r)\times\mathcal{O}}
\leq \langle\langle A\rangle\rangle_{\rho,r,\mathcal{O}}\me^{-\rho|i-j|},
\end{equation}
\begin{equation}\label{mtl12}
\| A^{12}_{ij}\|_{D(r)\times\mathcal{O}}
\leq \langle\langle A\rangle\rangle_{\rho,r,\mathcal{O}}\me^{-\rho|i+j|},
\end{equation}
\begin{equation}\label{mtl13}
\| A^{21}_{ij}\|_{D(r)\times\mathcal{O}}
\leq \langle\langle A\rangle\rangle_{\rho,r,\mathcal{O}}\me^{-\rho|i+j|},
\end{equation}
\begin{equation}\label{mtl14}
\| A^{22}_{ij}\|_{D(r)\times\mathcal{O}}
\leq \langle\langle A\rangle\rangle_{\rho,r,\mathcal{O}}\me^{-\rho|i-j|}.
\end{equation}
  \item[(T2$^{\prime}$) Asymptotically T\"{o}plitz] The limits
  $$\lim_{t\in\mathbb{Z},\, t\rightarrow\infty} A^{11}_{i+t,j+t},\,\,  \lim_{t\in\mathbb{Z},\,t\rightarrow\infty} A^{12}_{i+t,  j-t},\,\, \lim_{t\in\mathbb{Z},\,t\rightarrow\infty} A^{21}_{i+t,  j-t}\,\, \hbox{and}\,\,  \lim_{t\in\mathbb{Z},\,t\rightarrow\infty} A^{22}_{i+t, j+t}$$
  exist and are finite for all $i,j\in\mathbb{Z}.$
  \item[(T3$^{\prime}$) Lipschitz at infinity] For sufficiently large $|t|,\,\,t\in\mathbb{Z},$ the following  hold.
\begin{equation}\label{mtl11}
\|A^{11}_{i+t, j+t}-\lim_{t\rightarrow\infty} A^{11}_{i+t, j+t}\|_{D(r)\times\mathcal{O}}
\leq |t|^{-1}\langle\langle A\rangle\rangle_{\rho,r,\mathcal{O}}\me^{-\rho|i-j|}.
\end{equation}
\begin{equation}\label{mtl12}
\|A^{12}_{i+t, j-t}-\lim_{t\rightarrow\infty} A^{12}_{i+t, j-t}\|_{D(r)\times\mathcal{O}}
\leq |t|^{-1}\langle\langle A\rangle\rangle_{\rho,r,\mathcal{O}}\me^{-\rho|i+j|}.
\end{equation}
\begin{equation}\label{mtl13}
\|A^{21}_{i+t, j-t}-\lim_{t\rightarrow\infty} A^{21}_{i+t, j-t}\|_{D(r)\times\mathcal{O}}
\leq |t|^{-1}\langle\langle A\rangle\rangle_{\rho,r,\mathcal{O}}\me^{-\rho|i+j|}.
\end{equation}
\begin{equation}\label{mtl14}
\|A^{22}_{i+t, j+t}-\lim_{t\rightarrow\infty} A^{22}_{i+t, j+t}\|_{D(r)\times\mathcal{O}}
\leq |t|^{-1}\langle\langle A\rangle\rangle_{\rho,r,\mathcal{O}}\me^{-\rho|i-j|}.
\end{equation}
\end{description}

\end{definition}

\begin{definition}
Given a unilateral  infinite dimensional $\mathcal{M}_{2}(\mathbb{C})-$valued matrix
 $$A:\mathbb{N}\times \mathbb{N}\rightarrow \mathcal{M}_{2}(\mathbb{C}),$$
we lift it from  $\mathbb{N}\times \mathbb{N}$ to $\mathbb{Z}\times \mathbb{Z} $ by
 \begin{equation}\label{infMat2}
 \tilde{A}_{ij}=
\begin{cases}
A_{ij},  i\geq1,j\geq1,\\
0,\,\, \hbox{otherwise}.\\
\end{cases}
\end{equation}

We say that  $A$   is T\"{o}plitz-Lipschitz and write $A\in \mathfrak{M}^{\rho}_{r,\mathcal{O}}$ if
$\tilde{A}$ is T\"{o}plitz-Lipschitz and define
\begin{equation}\label{DecaySp2}
     \langle\langle A\rangle\rangle_{\rho,r,\mathcal{O}}:=\langle\langle \tilde{A}\rangle\rangle_{\rho,r,\mathcal{O}}<\infty.
\end{equation}
\end{definition}

The following conclusion indicate that  $\mathfrak{M}^{\rho}_{r,\mathcal{O}}$ is an algebra. This important   property   will be applied to  Proposition \ref{Flowest}.
\begin{proposition}\label{Malg}
Let $0<\delta< \rho.$
Suppose the  matrices  $A,B\in \mathfrak{M}^{\rho}_{r,\mathcal{O}}.$ Then their product $AB\in \mathfrak{M}^{\rho-\delta}_{r,\mathcal{O}}$  and there exists a constant $C>0$ so that
$$\langle\langle AB\rangle\rangle_{\rho-\delta,r,\mathcal{O}}\leq C\delta^{-1}\langle\langle A\rangle\rangle_{\rho,r,\mathcal{O}}\langle\langle B\rangle\rangle_{\rho,r,\mathcal{O}}.$$
\end{proposition}
The proof is given in Section \ref{pfMalg}, Appendix.

\subsection{Flow of linear Hamiltonian system}\label{flow}

In this section, we study the Hamiltonian flow generated
by a  quadratic  T\"{o}plitz-Lipschitz function $F(\theta,z,\bar{z};\xi)\in \mathcal{T}^\rho_{ D(r,s)\times\mathcal{O}}.$

In the sequel, we use the notations  $Z=(Z_j)^T_{j\in \mathbb{Z}}$ with $Z_j=(z_j, \bar{z}_j)^T.$
The Hessian $\partial^2_{Z}F$ of $F$ with respect to $Z$ reads
$$\partial^2_{Z}F=\left( \nabla_{Z_j}\nabla_{Z_i}F \right)_{i,j\in \mathbb{Z}}$$
where
$$\nabla_{Z_j}\nabla_{Z_i}F =\left(
                               \begin{array}{cc}
                                 \frac{\partial^2 F}{\partial z_i\partial z_j} &  \frac{\partial^2 F}{\partial z_i\partial \bar{z}_j} \\
                                  \frac{\partial^2 F}{\partial \bar{z}_i\partial z_j} &  \frac{\partial^2 F}{\partial \bar{z}_i\partial \bar{z}_j} \\
                               \end{array}
                             \right).
$$


Denote $A=J\partial^2_{Z}F,$
where
 $$J=diag\left\{J_j=\left(
                   \begin{array}{cc}
                     0 & 1 \\
                     -1 & 0 \\
                   \end{array}
                 \right)
\right\}_{j\in \mathbb{Z}},$$
 then
 \begin{equation}\label{A}
   A_{ij}=\left(
                         \begin{array}{cc}
                          \frac{\partial^2 F}{\partial \bar{z}_i\partial z_j} &  \frac{\partial^2 F}{\partial \bar{z}_i\partial \bar{z}_j} \\
                         -\frac{\partial^2 F}{\partial z_i\partial z_j} &  -\frac{\partial^2 F}{\partial z_i\partial \bar{z}_j} \\
                         \end{array}
                       \right).
 \end{equation}
By the definitions of   T\"{o}plitz-Lipschitz function and  T\"{o}plitz-Lipschitz matrix,  they have the following relation.
\begin{lemma}\label{AF}
Let  $\rho>0.$   Suppose  $F(\theta, z, \bar{z},; \xi)$ is a  quadratic function  on $D(r,s)\times\mathcal{O}.$
Then $F\in \mathcal{T}^\rho_{ D(r,s)\times\mathcal{O}}$ if and only if
 $A=J\partial^2_{Z}F\in \mathfrak{M}^{\rho}_{r,\mathcal{O}}.$
Moreover,
 \begin{equation}\label{equanor}
\langle\langle A\rangle\rangle_{\rho,r,\mathcal{O}}=\langle F\rangle_{\rho,D(r,s)\times\mathcal{O}}.
 \end{equation}
\end{lemma}

The Hamiltonian equation associated to  the quadratic function $F$ reads
\begin{equation}\label{Feq}
\begin{cases}
(\dot{\theta}(t), \dot{I}(t), \dot{z}(t), \dot{\bar{z}}(t))=X_F(\theta(t), I(t), z(t),  \bar{z}(t)),\\
(\theta(0),I(0),z(0), \bar{z}(0))=(\theta^0,I^0,z^0, \bar{z}^0).\\
\end{cases}
\end{equation}
Under the new notation $Z$,   the quadratic function
  $F\in \mathcal{T}^\rho_{ D(r,s)\times\mathcal{O}}$ can be rewritten as
 \begin{equation}\label{FZ}
F(\theta,Z)=\frac{1}{2}Z^TA(\theta)Z=\frac{1}{2}Z^T\partial^2_{Z}F(\theta,0)Z
 \end{equation}
and the
 equation  \eqref{Feq} reads
 \begin{equation}\label{Feq2}
\begin{cases}
\dot{\theta}(t)=0,\\
\dot{I}(t)=-\partial_\theta F(\theta(t), Z(t)),\\
\dot{Z}(t)=A(\theta(t))Z=J\partial_{Z}F(\theta(t),0)Z(t),\\
(\theta(0),I(0),Z(0))=(\theta^0,I^0,Z^0).\\
\end{cases}
\end{equation}
The Jacobian  $\partial_{Z^0} Z$ (the derivative  of $Z(t)$ with respect to $Z^0$) is
$$\partial_{Z^0} Z=\left(\partial_{Z^0_j} Z_i\right)_{i,j\in \mathbb{Z}}=\left(\left(
                          \begin{array}{cc}
                            \frac{\partial z_i}{\partial z^0_j} & \frac{\partial z_i}{\partial \bar{z}^0_j} \\
                            \frac{\partial \bar{z}_i}{\partial z^0_j} & \frac{\partial \bar{z}_i}{\partial \bar{z}^0_j} \\
                          \end{array}
                        \right)\right)_{i,j\in \mathbb{Z}}.$$

\begin{proposition}\label{Flowest}
Let  $0<\delta< \rho$ and
  $0<\sigma<r/3.$ Suppose the  quadratic  function $F\in \mathcal{T}^\rho_{ D(r,s)\times\mathcal{O}}$ and
 \begin{equation}\label{SmallF}
   \|X_F\|_{s;D(r-\sigma,s)\times\mathcal{O}}+\langle F\rangle_{\rho,D(r-\sigma,s)\times\mathcal{O}}<C \sigma.
 \end{equation}
 Then the
solution $(\theta(t), I(t),Z(t))$ of the equation \eqref{Feq2} with initial condition $(\theta^0,I^0,Z^0)\in D(r-\sigma,\frac{s}{4})$ satisfies
$(\theta(t), I(t), Z(t))$ $\in D(r,\frac{s}{2})$ for all $0\leq t\leq 1.$ Moreover,   the Jacobian  $\partial_{Z^0} Z(t)$ satisfies
 \begin{equation}\label{z11}
\langle\langle \partial_{Z^0} Z(t)-Id\rangle\rangle_{\rho-\delta,r-\sigma, \mathcal{O}}
 \leq C\langle F\rangle_{\rho,D(r-\sigma,s)\times\mathcal{O}}.
\end{equation}
where the notation  $Id$ is the identity mapping.

\end{proposition}

\begin{proof}

Since $\dot{\theta}(t)=0,$ then $\theta(t)\equiv \theta^0\in D(r-\sigma)$ remains unchanged.

Consider the equation for $Z:$
 \begin{equation}\label{Feq}
\begin{cases}
\dot{Z}=A(\theta^0)Z:=J\partial_{Z}F(\theta^0,0)Z,\\
Z(0)=Z^0.\\
\end{cases}
\end{equation}
It is a linear system with constant coefficients, thus its   solution is
 \begin{equation}\label{FeqZsol}
Z(t)=\me^{tA(\theta^0)}Z^0.
\end{equation}
By \eqref{A} and \eqref{SmallF},
$$\|A\|_{ \ell^{a,p}\rightarrow  \ell^{a,p}}\leq s \|X_F\|_{s;D(r-\sigma,s)\times\mathcal{O}}\leq Cs\sigma.$$
Thus for all $0\leq t \leq1$
$$\|Z(t)\|_{ \ell^{a,p}}\leq\me^{\|A\|_{ \ell^{a,p}\rightarrow  \ell^{a,p}}}\|Z^0\|_{ \ell^{a,p}}\leq \me^{Cs\sigma}\frac{s}{4}\leq \frac{s}{2}.$$

Consider the equation for $I.$ By \eqref{Feq2} and \eqref{FeqZsol}, we have
 \begin{equation}\label{FeqI}
\begin{cases}
\dot{I}(t)=-\frac{1}{2}Z^T\partial_\theta A(\theta)Z,\\
I(0)=I^0.\\
\end{cases}
\end{equation}
The integral form   of the above  equation  \eqref{FeqI} is
 \begin{equation}\label{FeqIsol}
I(t)=I^0-\frac{1}{2}\int^t_{0}Z^T(\tau)\partial_\theta A(\theta)Z(\tau)d\tau.
\end{equation}
Then for all $0\leq t \leq1$
$$|I(t)|\leq |I^0|+  \frac{1}{2\sigma}\|A\|_{\ell^{a,p}\rightarrow \ell^{a,p}}\|Z(t)\|^2_{\ell^{a,p}}\leq \frac{s}{4}+\frac{Cs^3}{8}\leq \frac{s}{2}.$$
Thus the flow $X^t_{F}$ exists for all $0\leq t\leq 1$ and it maps the domain $D(r-\sigma,\frac{s}{4})$ to
$D(r,\frac{s}{2})$. Denote the solution  $(\theta(t), I(t), z(t), \bar{z}(t))=X^t_{F}(\theta^0, I^0,  z^0, \bar{z}^0),$
then  for $0 \leq t\leq 1$ and  $(\theta^0, I^0,  z^0, \bar{z}^0)\in D(r-\sigma,\frac{s}{4}),$  the solution   $(\theta(t), I(t), z(t), \bar{z}(t))\in D(r,\frac{s}{2}).$

Now we prove the estimate \eqref{z11}.
Rewrite the solution   $Z(t)$ in  \eqref{FeqZsol}  as
 \begin{equation}\label{FeqZsol2}
Z(t)=(Id+B(t))Z^0,
\end{equation}
where $$B(t)=\me^{tA(\theta)}-Id=\sum^{\infty}_{k=1}\frac{t^k}{k!}A^k(\theta).$$
By Proposition  \ref{Malg} and Lemma  \ref{AF}, for all $0\leq t \leq1,$
\begin{equation}\label{Bt}
 \begin{split}
\langle\langle B\rangle\rangle_{\rho-\delta,r-\sigma, \mathcal{O}}
\leq&\sum^{\infty}_{k=1}\frac{(k-1)^{k-1}}{k!}(\frac{C}{\delta})^{k-1}\langle\langle A\rangle\rangle^k_{\rho,r-\sigma}\\
\leq&\sum^{\infty}_{k=1}\frac{\me^{k-1}}{k}(\frac{C}{\delta})^{k-1}\langle F\rangle^k_{\rho,D(r-\sigma,s)\times\mathcal{O}}\\
 \leq&C\langle F\rangle_{\rho,D(r-\sigma,s)\times\mathcal{O}}.
 \end{split}
\end{equation}
This completes the proof of the estimate \eqref{z11}.

\end{proof}

\begin{proposition}[\textbf{Canonical transformation}]\label{cabotra}
Let $0<3\delta<\rho,$ $0<\sigma<r$   and    $R,\,F\in \mathcal{T}^\rho_{ D(r,s)\times\mathcal{O}},$ where the Hamiltonian  $F$ is a   quadratic  function.   Assume that  the Hamiltonian $F$ satisfies \eqref{SmallF}.
Then the composition $R\circ X^1_F\in \mathcal{T}^{\rho-3\delta}_{ D(r-\sigma,s/4)\times\mathcal{O}}$ and
 there exists a constant  $C>0$ so that
\begin{equation}\label{Cano}
\langle R\circ X^1_F\rangle_{\rho-3\delta,D(r-\sigma,s/4)\times \mathcal{O}}\leq C\delta^{-2}\langle R\rangle_{\rho,D(r,s/2)\times \mathcal{O}}.
\end{equation}
\end{proposition}

\begin{proof}
By Proposition \ref{Flowest}, the time-1 mapping  $X^1_F$  maps $(\theta^0,I^0,Z^0)\in D(r-\sigma,\frac{s}{4})$ to
$(\theta, I, Z):=X^1_F(\theta^0,I^0,Z^0)\in D(r,\frac{s}{2}).$

Since  the mapping  $Z$ is linear in $Z^0,$ the Hessian
$\partial^2_{Z^0} Z =0.$
Then the  Hessian $\partial^2_{Z^0}(R\circ X^1_F)$ of $R\circ X^1_F$ with respect to $Z^0$ becomes
$$\partial^2_{Z^0}(R\circ X^1_F)=(\partial_{Z^0}Z)^T\partial^2_{Z}R(X^1_F)\partial_{Z^0}Z.$$
Note that $\langle\langle J^T(\partial_{Z^0} Z)^TJ\rangle\rangle_{\rho,r}=\langle\langle \partial_{Z^0} Z\rangle\rangle_{\rho,r},$
then by Lemma \ref{AF}  and Proposition \ref{Flowest}, we have
\begin{equation}\label{RFcom}
\begin{split}
& \langle R\circ X^1_F\rangle_{\rho-3\delta,D(r-\sigma,s/4)\times \mathcal{O}}\\
=&\langle\langle J\partial^2_{Z^0}(R\circ X^1_F)\rangle\rangle_{\rho-3\delta,D(r-\sigma,s/4)\times \mathcal{O}}\\
 \leq&C\delta^{-2}\langle\langle J^T(\partial_{Z^0} Z)^TJ\rangle\rangle_{\rho-\delta,r-\sigma}\langle\langle J\partial^2_{Z}R\rangle\rangle_{\rho,D(r,s/2)\times \mathcal{O}}
 \langle\langle \partial_{Z^0} Z\rangle\rangle_{\rho-\delta,r-\sigma}\\
 =&C\delta^{-2}\langle R\rangle_{\rho,D(r,s/2)\times \mathcal{O}}
 \langle\langle \partial_{Z^0} Z\rangle\rangle^2_{\rho-\delta,r-\sigma}\\
 \leq&C\delta^{-2}\langle R\rangle_{\rho,D(r,s/2)\times \mathcal{O}}.
 \end{split}
\end{equation}

\end{proof}

\section[]{A Reducibility  Theorem under  Brjuno Condition}\label{kamthem}

Consider the following  quadratic Hamiltonian with time quasi-periodic perturbation:
\begin{equation}\label{qpHam}
\begin{split}
  H(\omega t,z,\bar{z})=&\sum_{j\geq 1}\Omega_j z_j\bar{z}_j+P(\omega t,z,\bar{z})\\
  =&\sum_{j\geq 1}\Omega_j z_j\bar{z}_j+\sum\limits_{|\alpha|+|\beta|=2}P_{\alpha\beta}(\omega t)z^{\alpha}\bar{z}^{\beta},\\
\end{split}
\end{equation}
where  $(z,\bar{z})\in \ell^{a,p}_{0}\times \ell^{a,p}_{0},$  the space $\ell^{a,p}_{0}$
is the unilateral infinite sequences space  defined  in \eqref{uniseq}.
The forcing frequency vector  $\omega\in [0, 2\pi)^{n}$ and the normal   frequencies  $\Omega_j\in \mathbb{R}$ for all $j\geq1.$
Then the associated  linear Hamiltonian system reads
\begin{equation}\label{Heq0}
      \begin{cases}
            \dot{z}_j=\mi\Omega_j z_j+\mi \frac{\partial}{\partial \bar{z}_j} P(\omega t,z,\bar{z}),\quad           \\
            \dot{\bar{z}}_j=-\mi\Omega_j \bar{z}_j-\mi \frac{\partial }{\partial z_j} P(\omega t,z,\bar{z}),\quad    &  j\geq 1.\\
      \end{cases}
\end{equation}

Introducing the angle variables \,$\theta=\omega t\in \mathbb{T}^n,$\, and  the  auxiliary action variables \,$I\in \mathbb{R}^n,$
then we obtain an  autonomous Hamiltonian system
\begin{equation}\label{Heq1}
      \begin{cases}
           \dot{z}_j=\mi\Omega_j z_j+\mi \frac{\partial}{\partial \bar{z}_j} P(\omega t,z,\bar{z}),\quad           \\
            \dot{\bar{z}}_j=-\mi\Omega_j \bar{z}_j-\mi \frac{\partial }{\partial z_j} P(\omega t,z,\bar{z}),\quad    &  j\geq 1,\\
            \dot{\theta}_i=\omega_i,\quad                                                                              &   i=1\cdots n,\\
            \dot{I}_i=-\frac{\partial}{\partial \theta_i}P(\theta,z,\bar{z}),\quad                             &   i=1\cdots n.
      \end{cases}
\end{equation}
on the phase space \,$\mathcal{P}^{a,p}_0:=\mathbb{T}^n \times\mathbb{R}^n\times \ell^{a,p}_0\times \ell^{a,p}_0$\
with respect to the symplectic form
\begin{equation*}
      \sum_{i=1}^nd\theta_i\wedge dI_i+\mi\sum_{j\geq 1}dz_j\wedge d\bar{z}_j.
\end{equation*}
The new Hamiltonian is
\begin{equation}\label{Hf}
 \begin{split}
 H(\theta,I, z,\bar{z};\omega)=&N+P(\theta,z,\bar{z};\omega)\\
 =&\sum_{i=1}^n\omega_iI_i+\sum_{j\geq 1}\Omega_j z_j\bar{z}_j+\sum\limits_{|\alpha|+|\beta|=2}P_{\alpha\beta}(\theta; \omega)z^{\alpha}\bar{z}^{\beta}.
 \end{split}
\end{equation}
Given $s,r>0,$  in the following, we investigate Hamiltonian \eqref{Hf} on the domain $D(r,s) \subseteq \mathcal{P}^{a,p}_{0,\mathbb{C}}.$
The forcing frequency $\omega\in [0, 2\pi)^{n}$ will play the role of parameters.
Suppose $H(\theta,I, z,\bar{z};\omega)$ in \eqref{Hf} is real analytic on $(\theta,I, z,\bar{z};\omega)$ and $C^1_W-$smooth in compact subset
$\mathcal{O}\subseteq[0, 2\pi)^{n}$ with positive Lebesgue measure.
Furthermore,
suppose  Hamiltonian \eqref{Hf}  satisfies the following  assumptions.
 \begin{description}
   \item[(A1) Asymptotics of normal frequencies:]
   \begin{equation}\label{AsyNF}
   \Omega_j=j+\breve{\Omega}_j(\omega),\,\,j\geq 1,
   \end{equation}
where    $ \breve{\Omega}_j\in C^1_W(\mathcal{O})$ and there exist positive constants $A_0,$  such that
$\sup\limits_{j\geq 1,\omega\in \mathcal{O}}|\breve{\Omega}_j|\leq A_0.$
$\sup\limits_{j\geq 1}\sup\limits_{\omega\in \mathcal{O}}|\partial_\omega\breve{\Omega}_j|\leq \varepsilon_0.$

 \item[(A2) Non-resonance conditions: ]
There exist a constant $0<\gamma \leq 1$    and some fixed approximation function $\Delta$  such that  uniformly on $ \mathcal{O},$ for all $(k,l)\in \mathbb{Z}^n\times\mathbb{Z}^\infty\setminus\{0\},$
\begin{equation}\label{brjnonres}
  \begin{split}
  | k\cdot \omega|\geq&\frac{\gamma}{\Delta(|k|)},\, k\neq0,\\
 | k\cdot \omega+l\cdot \Omega(\omega)|\geq&\frac{\gamma}{\Delta(|k|)},\,\,|l|=2\\
   \end{split}
\end{equation}
where $|k|=|k_1|+\cdots+|k_n|,$ $|l|=\sum_{j}|l_j|.$

 \item[(A3) Regularity: ]
 The Hamiltonian vector field $X_P=(0, -P_\theta, \mi P_{\bar{z}}, -\mi P_z)^T$ of   perturbation $P$ defines a map $$X_P: D(r,s)\times \mathcal{O}\rightarrow   \mathcal{P}^{a,p}_{0,\mathbb{C}},$$
 $X_P(\cdot;\omega)$ is real analytic in $D(r,s)$ for each $\omega\in \mathcal{O},$ and
$P(\chi;\cdot)$ is $C^1_W-$smooth  in $\mathcal{O}$ for each $\chi\in D(r,s).$

\item[(A4) T\"{o}plitz-Lipschitz property: ]
 $\breve{\Omega}:=diag(\breve{\Omega}_j)_{j\geq 1}\in \mathfrak{M}^{\rho}_{r,\mathcal{O}}$ and
  $P\in \mathcal{T}^\rho_{ D(r,s)\times\mathcal{O}}$ for some $\rho>0.$

 \end{description}

Denote
\begin{equation}\label{NewPnorm}
  [ P ]^{\rho}_{s;D(r,s)\times\mathcal{O}} :=\|X_P\|_{s;D(r,s)\times\mathcal{O}} + \langle P\rangle_{\rho,D(r,s)\times \mathcal{O}}.
\end{equation}

\begin{theorem}\label{Redb}
Let  $\Delta$  be an  approximation function such that
\begin{equation}\label{appfun}
  \sum\limits_{k\in \mathbb{Z}^n}\frac{1}{\sqrt{\Delta(|k|)}}< +\infty.
\end{equation}
If the Hamiltonian  $H=N+P$ in \eqref{Hf} satisfies   the  above assumptions $(A1)-(A4)$  and there exists $0< \varepsilon_0 < \{\frac{\gamma}{4}(\sqrt{\Delta(1)}-1), (C_*\gamma2^5)^{\frac{3}{2}}\}$ so that
$$\langle\langle \breve{\Omega}\rangle\rangle_{\rho,r,\mathcal{O}}<\varepsilon_0\,\,
\hbox{and }\,\,
  [ P ]^{\rho}_{s;D(r,s)\times\mathcal{O}} <\varepsilon_0.$$
  Then there exist
\begin{description}
  \item[(i)]
 a Cantor subset $\mathcal{O}_\gamma\subset\mathcal{O}$ with Lebesgue measure $\mes(\mathcal{O}\setminus\mathcal{O}_\gamma)=O(\sqrt{\gamma})$ as $\gamma\rightarrow 0$;
  \item[(ii)] a $C^1_W-$smooth family of  real analytic,  symplectic  coordinate transformations
$\Phi: \mathcal{P}^{a,0}_{0}\times\mathcal{O}_\gamma\rightarrow  \mathcal{P}^{a,0}_{0}$ of the form
\begin{equation}\label{lintransform}
  \Phi_{\omega}\left(
    \begin{array}{c}
     \theta \\
     I\\
      Z \\
    \end{array}
  \right)
=
\left(
  \begin{array}{c}
   \theta  \\
   I+\frac{1}{2}Z^TM_{\omega}(\theta)Z\\
    L_{\omega}(\theta)Z \\
  \end{array}
\right)
\end{equation}
where   $Z=(Z_j)^T_{j\geq1}$ with $Z_j=(z_j, \bar{z}_j)^T.$    $M_{\omega}(\theta)$
and
$L_{\omega}(\theta)$ are linear bounded operators on $\ell^{a,p}_{0} \times \ell^{a,p}_{0}$ for all $p\geq 0$,   and $L_{\omega}(\theta)$ is also invertible;
  \item[(iii)]  a $C^1_W-$smooth family of  new normal forms
  \begin{equation}\label{Ninf}
   \begin{split}
 N^{\infty}=&\sum_{j=1}^n\omega_jI_j+\sum_{j\geq 1}\Omega^{\infty}_j z_j\bar{z}_j\\
 \end{split}
\end{equation}
such that on $\mathcal{P}^{a,0}_{0}\times\mathcal{O}_\gamma,$
$$H\circ \Phi=N^{\infty}.$$
Moreover the new normal frequencies are close to the original ones
$$|\Omega^{\infty}-\Omega|_{\mathcal{O}_\gamma}
\leq c \varepsilon,$$
and the the new frequencies satisfy a non-resonant condition: for all $\omega\in \mathcal{O}_\alpha,$
$$| k\cdot \omega|\geq\frac{\gamma}{2\Delta(|k|)},\,\forall k\neq0,$$
$$| k\cdot \omega+l\cdot \Omega^{\infty}(\omega)|\geq\frac{\gamma}{2\Delta(|k|)},\,\forall k\in \mathbb{Z}^n,\, |l|=2.$$
\end{description}
\end{theorem}

\section[]{Applications to some linear Hamiltonian PDEs}\label{App}

We give the proof of Theorem \ref{maintheorem} by  Theorem \ref{Redb}.

\subsection{The  Hamiltonian derivative wave equations}
We consider the wave equation \eqref{eq1}.
Let
\begin{equation}\label{chavar}
      \begin{cases}
           w= \frac{1}{\sqrt{2}}(\mathbf{D}_mu +\mi u_t) ,\\
            \bar{w}=\frac{1}{\sqrt{2}}(\mathbf{D}_mu -\mi u_t).
      \end{cases}
\end{equation}
Then the  equation \eqref{eq1} is  written as a  non-autonomous Hamiltonian equation
\begin{equation}\label{Heq1}
      \begin{cases}
            w_t= -\mi \frac{\partial}{\partial \bar{w}}H(t,w,\bar{w})=-\mi \mathbf{D}_mw- \frac{\mi\epsilon}{2} V(\omega t, x) (w+ \bar{w}),\\
           \bar{w}_t=\mi \frac{\partial}{\partial w}H(t,w,\bar{w})=\mi \mathbf{D}_mw+\frac{\mi\epsilon}{2} V(\omega t, x) (w+ \bar{w}).
      \end{cases}
\end{equation}
with  the Hamiltonian
\begin{equation*}
 \begin{split}
    H(t,w,\bar{w})
    =&\int^{\pi}_{0}\left[\bar{w}\mathbf{D}_mw+\frac{\epsilon}{2}V(\omega t, x) (w+ \bar{w})^2\right] dx.
 \end{split}
\end{equation*}

Recall the  function  space $\mathcal{H}^{a,p}_0$ in \eqref{fB1}.
Through  the inverse discrete Fourier transform $\mathcal{S}:\ell^{a,p}_0\rightarrow \mathcal{H}^{a,p}_0$, the space
$\mathcal{H}^{a,p}_0$  can  be identified with the space $\ell^{a,p}_0.$

We expand $w(t, x),\,\, \bar{w}(t, x)$ on the eigenfunctions
\,$$w(t, x)=\sum\limits_{j\geq1}q_j(t)\phi_j(x)\in \mathcal{H}^{a,p}_0,\,\,\bar{w}(t, x)=\sum\limits_{j\geq1}\bar{q}_j(t)\phi_j(x)\in \mathcal{H}^{a,p}_0$$
with $q=(q_j)_{j\geq1},\, \bar{q}=(\bar{q}_j)_{j\geq1}  \in  \ell^{a,p}_0.$
Then the  equation \eqref{Heq1} becomes
\begin{equation}
      \begin{cases}
            \dot{q}_j=-\mi\frac{\partial}{\partial \bar{q}_j}H(t,q,\bar{q})=-\mi \lambda_jq_j-\mi\frac{\partial}{\partial \bar{q}_j}G,\\
            \dot{\bar{q}}_j=\mi\frac{\partial}{\partial q_j}H(t,q,p)=\mi \lambda_jq_j+\mi\frac{\partial}{\partial \bar{q}_j}G,
      \end{cases}
\end{equation}
where
$$H(t,q, p)=\Lambda+G,$$
\begin{equation*}
 \begin{split}
    \Lambda=&\sum\limits_{j\geq1}\lambda_jq_j \bar{q}_j,\\
    G=&\frac{\epsilon}{2}\sum\limits_{i,j\geq1} (q_i+\bar{q}_i) (q_j+\bar{q}_j)\int^{\pi}_{0}  V(t\omega,x)\phi_i(x)\phi_j(x)dx.
 \end{split}
\end{equation*}
Now we
introduce the angle variables \,$\theta=\omega t\in \mathbb{T}^n,$\,  the auxiliary action variables \,$I\in \mathbb{R}^n$ and
the complex coordinates \,$z=(z_j)_{j\geq1},\,\,  \bar{z}=(\bar{z}_j)_{j\geq1}$\, via letting  $z_j = -q_j,\,\, \bar{z}_j = -\bar{q}_j.$
Then we obtain an  autonomous Hamiltonian system
\begin{equation}\label{atHeq1w}
      \begin{cases}
            \dot{z}_j=\mi\lambda_j z_j+\mi \frac{\partial}{\partial \bar{z}_j} P(\theta,z,\bar{z})\quad            &   j\geq 1,\\
            \dot{\bar{z}}_j=-\mi\lambda_j \bar{z}_j-\mi \frac{\partial }{\partial z_j} P(\theta,z,\bar{z})\quad    &  j\geq 1,\\
            \dot{\theta}_i=\omega_i\quad                                                                              &   i=1\cdots n,\\
            \dot{I}_i=-\frac{\partial}{\partial \theta_i}P(\theta,z,\bar{z})\quad                             &  i=1\cdots n.
      \end{cases}
\end{equation}
on the phase space \,$\mathcal{P}^{a,p}_0$\
with respect to the symplectic form
\begin{equation*}
      \sum_{i=1}^nd\theta_i\wedge dI_i+\mi\sum_{j\geq 1}dz_j\wedge d\bar{z}_j.
\end{equation*}
The Hamiltonian associated to the system \eqref{atHeq1w} is
\begin{equation}\label{HDWlam}
 \begin{split}
 H=&N+P\\
  \end{split}
\end{equation}
 where
 \begin{equation}\label{HDWlamNP}
 \begin{split}
 N=&\sum_{j=1}^n\omega_jI_j+\sum_{j\geq 1}\lambda_j z_j\bar{z}_j,\\
    P=&\frac{\epsilon}{2}\sum\limits_{i,j\geq1} (z_i+\bar{z}_i) (z_j+\bar{z}_j)\int^{\pi}_{0}  V(\theta,x)\phi_i(x)\phi_j(x)dx.
 \end{split}
\end{equation}


In the following, we check that the Hamiltonian \eqref{HDWlam} satisfies the  assumptions \textbf{(A1)-(A4)}.
 Let $r$ be that in Assumption 1.1 and $s>0$ be a suitable positive number.
Take  $\varepsilon_0=(2^{p+1}+2^{4}+ \frac{18 n}{r})C_V\epsilon>0.$
\begin{description}
  \item[(1)] \emph{Verifying  the assumption} \textbf{(A1)}.

  Since $\lambda_j =\sqrt{j^2+m}=j+\frac{m}{2j}-\frac{m^2}{8j^3}+\cdots,$ then we
take   $\Omega_j=j+\breve{\Omega}_j=j+O(\frac{1}{j}).$
Note that $\breve{\Omega}_j$ does not depend on $\omega\in [0, 2\pi)^{n},$
thus $\partial_\omega\breve{\Omega}_j=0$ and   $ \breve{\Omega}_j\in C^1_W([0, 2\pi)^{n}).$
Take  $A_0=1+m.$ Since  $\breve{\Omega}_j=O(\frac{1}{j})$ and $\partial_\omega\breve{\Omega}_j=0,$  then for all $j\geq 1$  and $\omega\in [0, 2\pi)^{n},$ we have
$|\breve{\Omega}_j|\leq A_0$ and
$|\partial_\omega\breve{\Omega}_j|\leq \varepsilon_0.$

  \item[(2)] \emph{Verifying the assumption}  \textbf{(A2)}.

Take the  vector $v=(\sgn(k_1),\cdots,\sgn(k_n))$ then $k\cdot v=|k|.$ Let $\omega=\omega_s=sv+w$ with $s\in \mathbb{R},$
$w\in v^{\perp}.$
Consider the function  $f(s)=k\cdot \omega_s + l\cdot \Omega=|k|s +k\cdot w + l\cdot \Omega.$
Thanks to $\partial_\omega\Omega=0,$ we have
$$ | f'(s)|=|k|.$$
By Lemma \ref{measest}   in Appendix, we have
$$\mes\{s:sv+w\in [0, 2\pi)^{n}, |f(s)|\leq \delta\}\leq \frac{4\delta}{|k|}.$$
It follows that  the  measure
\begin{equation}
\begin{split}
    &\mes\{\omega \in  [0, 2\pi)^{n} :|k\cdot \omega + l\cdot \Omega|\leq \frac{\gamma}{\Delta(|k|)},|l|=0,2\}\\
 \leq &\diam^{n-1}([0, 2\pi)^{n})\mes\{s:sv+w\in [0, 2\pi)^{n}, |f(s)|\leq \frac{\gamma}{\Delta(|k|)}\}\\
 \leq &  (2\pi)^{n(n-1)}\frac{4\gamma}{|k|\Delta(|k|)}.
 \end{split}
\end{equation}
Thus  there is a subset $\mathcal{O}\subset[0, 2\pi)^{n}$ of positive Lebesgue measure
with $\mes\mathcal{O}\geq (2\pi)^{n}(1-O(\gamma))$ such that  the assumption\textbf{(A2)} holds on $\mathcal{O}.$

  \item[(3)] \emph{Verifying the assumption}  (A3).

The perturbation $P$ in \eqref{HDWlamNP} reads
\begin{equation*}
 \begin{split}
    P(\theta, z, \bar{z})=&\frac{\epsilon}{2}\sum\limits_{ij\geq1}p^{20}_{ij}(\theta) z_iz_j + \epsilon\sum\limits_{ij\geq1}p^{11}_{ij}(\theta) z_i\bar{z}_j+\frac{\epsilon}{2}\sum\limits_{ij\geq1}p^{02}_{ij}(\theta) \bar{z}_i\bar{z}_j ,
 \end{split}
\end{equation*}
where
 \begin{equation}\label{dwpij1}
 \begin{split}
p^{20}_{ij}(\theta)=p^{11}_{ij}(\theta)=p^{02}_{ij}(\theta)
=&\int^{\pi}_{0}  V(\theta,x)\phi_i(x)\phi_j(x)dx \\
=&  \begin{cases}
           \frac{1}{2}(\widetilde{V}_{i-j}(\theta)- \widetilde{V}_{i+j}(\theta)), \quad     i>j, \\
          \widetilde{V}_{0}(\theta)  -\frac{1}{2} \widetilde{V}_{2j}(\theta),  \quad     i=j,  \\
          \frac{1}{2}(\widetilde{V}_{j-i}(\theta)- \widetilde{V}_{i+j}(\theta)), \quad     i<j.
         \end{cases}\\
\end{split}
\end{equation}

Now we investigate the regularity of the  perturbation  vector field $X_P=(0, -\frac{\partial P}{\partial \theta},  \mi\frac{\partial P}{\partial \bar{z}}, -\mi \frac{\partial P}{\partial z}).$
Note that the   vector field $X_P$ does not depends on $\omega.$  For the above $r, s>0,$  we estimate   the   vector field  norm
 \begin{equation*}
 \begin{split}
&\|X_P\|_{s;D(r,s)\times\mathcal {O}}\\
=&\frac{1}{s^2}\sum\limits^\textrm{n}_{h=1}\left\|\frac{\partial P}{\partial \theta_h}\right\|_{D(r,s)\times\mathcal {O}}
+\frac{1}{s}\sup\limits_{\|z\|_{a,p}<s,\atop \|\bar{z}\|_{a,p}<s}\sum\limits^{\infty}_{i=1}i^{p}\me^{ai }\left(\left\|\frac{\partial P}{\partial \bar{z}_i}\right\|_{D(r)\times\mathcal {O}}
+\left\|\frac{\partial P}{\partial z_i}\right\|_{D(r)\times\mathcal {O}}\right).
\end{split}
\end{equation*}

$\bullet$ We first estimate  the sum
\begin{equation*}
 \begin{split}
&\sum\limits^\textrm{n}_{h=1}\left\|\frac{\partial P}{\partial \theta_h}\right\|_{D(r,s)\times\mathcal {O}}\\
=&\epsilon\sup\limits_{\|z\|_{a,p}<s,\atop \|\bar{z}\|_{a,p}<s}\sum\limits^\textrm{n}_{h=1} \sum\limits_{i, j\geq1}\left(\frac{1}{2}\|\frac{\partial p^{20}_{ij}}{\partial \theta_h}\|_{D(r)}|z_i||z_j|
+\|\frac{\partial p^{11}_{ij}}{\partial \theta_h}\|_{D(r)}|z_i||\bar{z}_j|+\frac{1}{2}\|\frac{\partial p^{02}_{ij}}{\partial \theta_h}\|_{D(r)}|\bar{z}_i||\bar{z}_j|\right)\\
\leq&\frac{n \epsilon}{r}\sup\limits_{\|z\|_{a,p}<s,\atop \|\bar{z}\|_{a,p}<s}  \left(\sum\limits_{i, j\geq1}\|p^{20}_{ij}\|_{D(2r)}|z_i||z_j|
+\sum\limits_{i, j\geq1}\|p^{11}_{ij}\|_{D(2r)}|z_i||\bar{z}_j|+\sum\limits_{i, j\geq1}\|p^{02}_{ij}\|_{D(2r)}|\bar{z}_i||\bar{z}_j|\right).
\end{split}
\end{equation*}
For this purpose, it suffices to estimate each of three sums  on the last line:
\begin{equation*}
 \begin{split}
 &\sum\limits_{i, j\geq1}\|p^{11}_{ij}\|_{D(2r)}|z_i||\bar{z}_j|\\
=&\sum\limits_{j\geq1} \|p^{11}_{jj}\|_{D(2r)}|z_j||\bar{z}_j|
+\sum\limits_{j\geq1}  \sum\limits_{1\leq i\leq j-1} \|p^{11}_{ij}\|_{D(2r)}|z_i||\bar{z}_j|
+\sum\limits_{j\geq1}  \sum\limits_{i\geq j+1} \|p^{11}_{ij}\|_{D(2r)}|z_i||\bar{z}_j|
\\
\leq&\sum\limits_{j\geq1} (\| \widetilde{V}_{0}(\theta)\|_{D(2r)}+\| \widetilde{V}_{2j}(\theta)\|_{D(2r)})|z_j||\bar{z}_j|\\
&+\sum\limits_{j\geq1}  \sum\limits_{1\leq i\leq j-1} (\|\widetilde{V}_{j-i}\|_{D(2r)}+\|\widetilde{V}_{j+i}\|_{D(2r)})|z_i||\bar{z}_j|\\
&+\sum\limits_{j\geq1}  \sum\limits_{i\geq j+1} (\|\widetilde{V}_{i-j}\|_{D(2r)}+\|\widetilde{V}_{j+i}\|_{D(2r)})|z_i||\bar{z}_j|
\\
\leq&6C_V\|z\|_{a,p}\|\bar{z}\|_{a,p}.
\end{split}
\end{equation*}
Similarly, we have  $$\sum\limits_{i, j\geq1}\|p^{20}_{ij}\|_{D(2r)}|z_i||z_j|\leq  6C_V\|z\|_{a,p}\|z\|_{a,p}$$ and
 $$\sum\limits_{i, j\geq1}\|p^{02}_{ij}\|_{D(2r)}|\bar{z}_i||\bar{z}_j|\leq  6C_V\|\bar{z}\|_{a,p}\|\bar{z}\|_{a,p}.$$
This  shows  that
\begin{equation}\label{Ptheta}
  \frac{1}{s^2}\sum\limits^\textrm{n}_{h=1}\|\frac{\partial P}{\partial \theta_h}\|_{D(r,s)\times\mathcal {O}}\leq  \frac{18 n}{r}C_V \epsilon.
\end{equation}

$\bullet$ We turn to the estimate for  $$\frac{1}{s}\sup\limits_{\|z\|_{a,p}<s,\atop \|\bar{z}\|_{a,p}<s}\sum\limits^{\infty}_{i=1}i^{p}\me^{ai }\left(\left\|\frac{\partial P}{\partial \bar{z}_i}\right\|_{D(r)\times\mathcal {O}}
+\left\|\frac{\partial P}{\partial z_i}\right\|_{D(r)\times\mathcal {O}}\right).$$
It suffices to consider
\begin{equation*}
 \begin{split}
&\sum\limits^{\infty}_{i=1}i^{p}\me^{ai }\left\|\frac{\partial P}{\partial z_i}\right\|_{D(r)\times\mathcal {O}}
=\epsilon\sum\limits^{\infty}_{i=1}\sum\limits^{\infty}_{j=1}i^{p}\me^{ai }\|p^{11}_{ij}\|_{D(r)}|z_j|+\epsilon\sum\limits^{\infty}_{i=1}\sum\limits^{\infty}_{j=1}i^{p}\me^{ai }\|p^{02}_{ij}\|_{D(r)}|\bar{z}_j|.\\
\end{split}
\end{equation*}
By \eqref{dwpij1},
\begin{equation*}
 \begin{split}
&\sum\limits^{\infty}_{i=1}\sum\limits^{\infty}_{j=1}i^{p}\me^{ai }\|p^{11}_{ij}\|_{D(r)}|z_j|\\
=&\sum\limits_{j\geq1}j^{p}\me^{aj } \| \widetilde{V}_{0}(\theta)  -\frac{1}{2} \widetilde{V}_{2j}(\theta)\|_{D(2r)}|z_j|\,\,\,\cdots(\ast1)\\
+&\sum\limits_{j\geq1}  \sum\limits_{1\leq i\leq j-1} i^{p}\me^{ai }\|\frac{1}{2}(\widetilde{V}_{j-i}(\theta)- \widetilde{V}_{i+j}(\theta))\|_{D(2r)}|z_j|\,\,\,\cdots(\ast2)\\
+&\sum\limits_{j\geq1}  \sum\limits_{i\geq j+1}i^{p}\me^{ai } \|\frac{1}{2}(\widetilde{V}_{i-j}(\theta)- \widetilde{V}_{i+j}(\theta))\|_{D(2r)}|z_j|\,\,\,\cdots(\ast3),\\
\end{split}
\end{equation*}
where
\begin{equation*}
 \begin{split}
(\ast1)=\sum\limits_{j\geq1}j^{p}\me^{aj } \| \widetilde{V}_{0}(\theta)  -\frac{1}{2} \widetilde{V}_{2j}(\theta)\|_{D(2r)}|z_j|
\leq 2\|V\|_{D(2r),b,p}\|z\|_{a,p}.
\end{split}
\end{equation*}
\begin{equation*}
 \begin{split}
(\ast2)=&\sum\limits_{j\geq1}  \sum\limits_{1\leq i\leq j-1} i^{p}\me^{ai }\|\frac{1}{2}(\widetilde{V}_{j-i}(\theta)- \widetilde{V}_{i+j}(\theta))\|_{D(2r)}|z_j|\\
\leq &   \sum\limits_{j\geq1}  \sum\limits_{1\leq i\leq j-1} i^{p}\me^{ai }\frac{1}{2}  \|V\|_{D(2r),b,p} (j-i)^{-p} \me^{-b(j-i)}       |z_j|  \\
     & +  \sum\limits_{j\geq1}  \sum\limits_{1\leq i\leq j-1} i^{p}\me^{ai }\frac{1}{2}  \|V\|_{D(2r),b,p} (j+i)^{-p} \me^{-b(j+i)}       |z_j|\\
 \leq & \frac{1}{2}  \|V\|_{D(2r),b,p}  \sum\limits_{j\geq1}  \sum\limits_{1\leq i\leq j-1} (\frac{i}{j-i})^p \me^{ai }  \me^{-b(j-i)}       |z_j|  \\
     & +\frac{1}{2}  \|V\|_{D(2r),b,p}  \sum\limits_{j\geq1}  \sum\limits_{1\leq i\leq j-1} (\frac{i}{j+i})^p \me^{ai }  \me^{-b(j+i)}       |z_j|\\
  \leq & \frac{1}{2}  \|V\|_{D(2r),b,p}  \sum\limits_{j\geq1}   |z_j| j^p 2\me^{aj}
     +\frac{1}{2}  \|V\|_{D(2r),b,p}  \sum\limits_{j\geq1} 2 |z_j|      \\
  \leq  & 2C_V\|z\|_{a,p}.
\end{split}
\end{equation*}

\begin{equation*}
 \begin{split}
(\ast3)=&\sum\limits_{j\geq1}  \sum\limits_{i\geq j+1}i^{p}\me^{ai } \|\frac{1}{2}(\widetilde{V}_{i-j}(\theta)- \widetilde{V}_{i+j}(\theta))\|_{D(2r)}|z_j|\\
\leq &   \sum\limits_{j\geq1}  \sum\limits_{i\geq j+1} i^{p}\me^{ai }\frac{1}{2}  \|V\|_{D(2r),b,p} (i-j)^{-p} \me^{-b(i-j)}       |z_j|  \\
     & +  \sum\limits_{j\geq1}  \sum\limits_{i\geq j+1} i^{p}\me^{ai }\frac{1}{2}  \|V\|_{D(2r),b,p} (j+i)^{-p} \me^{-b(j+i)}       |z_j|\\
 \end{split}
\end{equation*}
\begin{equation*}
 \begin{split}
 \leq & \frac{1}{2}  \|V\|_{D(2r),b,p}  \sum\limits_{j\geq1}  \sum\limits_{i\geq j+1} (\frac{i}{i-j})^p \me^{ai }  \me^{-b(i-j)}       |z_j|  \\
     & +\frac{1}{2}  \|V\|_{D(2r),b,p}  \sum\limits_{j\geq1}  \sum\limits_{i\geq j+1} (\frac{i}{j+i})^p \me^{ai }  \me^{-b(j+i)}       |z_j|\\
  \leq  & (2^p +2)C_V\|z\|_{a,p}.
\end{split}
\end{equation*}
Then
$$\sum\limits^{\infty}_{i=1}\sum\limits^{\infty}_{j=1}i^{p}\me^{ai }\|p^{11}_{ij}\|_{D(r)}|z_j|\leq (\ast1)+(\ast2)+(\ast3) \leq  (2^p +6)C_V\|z\|_{a,p}.$$
By the similar argument, we get
$$\sum\limits^{\infty}_{i=1}\sum\limits^{\infty}_{j=1}i^{p}\me^{ai }\|p^{02}_{ij}\|_{D(r)}|\bar{z}_j|\leq (2^p +6)C_V\|\bar{z}\|_{a,p}. $$
It follows that
\begin{equation}\label{Pz}
  \frac{1}{s}\sup\limits_{\|z\|_{a,p}<s,\atop \|\bar{z}\|_{a,p}<s}\sum\limits^{\infty}_{i=1}i^{p}\me^{ai }\left(\|\frac{\partial P}{\partial \bar{z}_i}\|_{D(r)\times\mathcal {O}}
+\|\frac{\partial P}{\partial z_i}\|_{D(r)\times\mathcal {O}}\right)\leq 2(2^p +6)C_V\epsilon .
\end{equation}

We conclude   from \eqref{Ptheta} and \eqref{Pz} that
\begin{equation*}
 \begin{split}
&\|X_P\|_{s;D(r,s)\times\mathcal {O}}\\
=&\frac{1}{s^2}\sum\limits^\textrm{n}_{h=1}\|\frac{\partial P}{\partial \theta_h}\|_{D(r,s)\times\mathcal {O}}
+\frac{1}{s}\sup\limits_{\|z\|_{a,p}<s,\atop \|\bar{z}\|_{a,p}<s}\sum\limits^{\infty}_{i=1}i^{p}\me^{ai }\left(\|\frac{\partial P}{\partial \bar{z}_i}\|_{D(r)\times\mathcal {O}}
+\|\frac{\partial P}{\partial z_i}\|_{D(r)\times\mathcal {O}}\right)\\
\leq&(2^{p+1}+12+ \frac{18 n}{r})C_V\epsilon\leq \varepsilon_0.
\end{split}
\end{equation*}
Thus we  complete the verification of the regularity for  $X_P.$

  \item[(4)] \emph{Verifying the assumption}  (A4).

  $\bullet$We verify  $\breve{\Omega}:=diag(\breve{\Omega}_j)_{j\geq1}$  satisfies  T\"{o}plitz-Lipschitz property.
During the verification of the   assumption  (A1), we have obtained   $|\breve{\Omega}_j|\leq \frac{C_0}{j},$  where $C_0$ is a constant depending on $m.$
It is evident that $\lim\limits_{t\rightarrow\infty}\breve{\Omega}_{j+t}=0$
and $$\left\|\lim_{t\rightarrow\infty}\breve{\Omega}_{j+t}\right\|_{\mathcal{O}}\leq C_0.$$
 $$\left\|\breve{\Omega}_{j+t}- \lim_{t\rightarrow\infty}\breve{\Omega}_{j+t}\right\|_{\mathcal{O}}=\left\|\breve{\Omega}_{j+t}- \lim_{t\rightarrow\infty}\breve{\Omega}_{j+t}\right\|_{\mathcal{O}}\leq \frac{C_0}{|j+t|}\leq\frac{C_0}{|t|}.$$

$\bullet$Taking $\rho=2a,$
 we verify the perturbation $P\in \mathcal{T}^\rho_{ D(r,s)\times\mathcal{O}}.$

We first  consider $\frac{\partial^2 P}{\partial z_i \partial \bar{z}_j}.$
By \eqref{dwpij1}, we have for $t\geq1,$
\begin{equation*}
   \frac{\partial^2 P}{\partial z_{i+t} \partial \bar{z}_{j+t}}=\epsilon p^{11}_{i+t,j+t}(\theta)
 = \begin{cases}
           \frac{\epsilon }{2}(\widetilde{V}_{i-j}(\theta)- \widetilde{V}_{i+j+2t}(\theta)), \quad     i>j ,\\
          \epsilon \widetilde{V}_{0}(\theta)  -\frac{\epsilon }{2} \widetilde{V}_{2j+2t}(\theta),  \quad     i=j , \\
          \frac{\epsilon }{2}(\widetilde{V}_{j-i}(\theta)- \widetilde{V}_{i+j+2t}(\theta)), \quad     i<j.
         \end{cases}\\
\end{equation*}

Due to $  \|\widetilde{V}_{j}\|_{D(2r)}\leq   C_V  \me^{-2aj}, j\geq 1, $ the limit  $\lim\limits_{t\rightarrow\infty}\frac{\partial^2 P}{\partial z_{i+t} \partial \bar{z}_{j+t}}$ exists and
\begin{equation*}
   \lim\limits_{t\rightarrow\infty}\frac{\partial^2 P}{\partial z_{i+t} \partial \bar{z}_{j+t}}
 = \begin{cases}
           \frac{\epsilon }{2}\widetilde{V}_{i-j}(\theta), \quad     i>j ,\\
          \epsilon \widetilde{V}_{0}(\theta)  ,  \quad     i=j , \\
          \frac{\epsilon }{2}\widetilde{V}_{j-i}(\theta), \quad     i<j.
         \end{cases}\\
\end{equation*}
Moreover,
\begin{equation*}
  \left\|\lim\limits_{t\rightarrow\infty} \frac{\partial^2 P}{\partial z_{i+t} \partial \bar{z}_{j+t}}\right\|_{D{(r,s)}\times\mathcal{O}}
  \leq\epsilon \left\|\widetilde{V}_{|i-j|}(\theta)\right\|_{D{(r,s)}\times\mathcal{O}}
  \leq\varepsilon_0 \me^{-\rho|i-j|}.\\
\end{equation*}

Thanks  to the exponentially decay of $\widetilde{V}_{j},$ we also have
\begin{equation*}
 \begin{split}
&\left\|\frac{\partial^2 P}{\partial z_{i+t} \partial \bar{z}_{j+t}}-\lim\limits_{t\rightarrow\infty}\frac{\partial^2 P}{\partial z_{i+t} \partial \bar{z}_{j+t}}\right\|_{D{(r,s)}\times\mathcal{O}}\\
= &  \begin{cases}
           \frac{\epsilon }{2}\left\| \widetilde{V}_{i+j+2t}(\theta))\right\|_{D{(r,s)}\times\mathcal{O}}, \quad     i>j ,\\
           \frac{\epsilon }{2} \left\| \widetilde{V}_{2j+2t}(\theta)\right\|_{D{(r,s)}\times\mathcal{O}},  \quad     i=j , \\
          \frac{\epsilon }{2}\left\| \widetilde{V}_{i+j+2t}(\theta))\right\|_{D{(r,s)}\times\mathcal{O}}, \quad     i<j,
         \end{cases}\\
\leq &  \frac{\epsilon }{2}   \|V\|_{D(2r),b,p}   \me^{-2a(i+j+2t)} \leq   \frac{\varepsilon_0}{t}\me^{-\rho|i-j|}.
\end{split}
\end{equation*}
where we use the inequality
$\me^{-2a(i+j+2t)}=\me^{-2a(i+j)}\me^{-2 t}\leq \frac{1}{t}\me^{-2a|i-j|}.$

As to the second derivative   $\frac{\partial^2 P}{\partial z_i \partial z_j}=\frac{\epsilon }{2} p^{20}_{ij}(\theta),$
we consider the lift  $\widetilde{P}(\theta, \tilde{z}, \bar{\tilde{z}})= P(\theta, z, \bar{z}),$ where $(\tilde{z}, \bar{\tilde{z}})\in \ell^{a,p}\times\ell^{a,p}$ and $\tilde{z}=z_j, \bar{\tilde{z}}=\bar{z}_j$ when $ j\geq1.$
(recall the Definition \ref{liftF}).  Then
\begin{equation}\label{liftpsecd}
\frac{\partial^2 \widetilde{P}}{\partial \tilde{z}_i \partial \tilde{z}_j}=
\begin{cases}
\frac{\partial^2 P}{\partial z_i \partial z_j},  i\geq1,j\geq1,\\
0,\,\, \hbox{otherwise}.\\
\end{cases}
\end{equation}

When $|t|$ is sufficiently large, we have either $i+t<0$  or $j-t<0,$ then  $ \frac{\partial^2 \widetilde{P}}{\partial \tilde{z}_{i+t} \partial \tilde{z}_{j-t}}=0$
  and thus the limit  $\lim\limits_{t\rightarrow\infty}\frac{\partial^2 \widetilde{P}}{\partial \tilde{z}_{i+t} \partial \tilde{z}_{j-t}}=0.$ It is obvious that
\begin{equation*}
  \left\|\lim\limits_{t\rightarrow\infty} \frac{\partial^2 \widetilde{P}}{\partial \tilde{z}_{i+t} \partial \tilde{z}_{j-t}}\right\|_{D{(r,s)}\times\mathcal{O}}
  \leq\varepsilon_0 \me^{-\rho|i+j|}\\
\end{equation*}
and
\begin{equation*}
 \begin{split}
\left\|\frac{\partial^2 \widetilde{P}}{\partial \tilde{z}_{i+t} \partial \tilde{z}_{j-t}}-\lim\limits_{t\rightarrow\infty}\frac{\partial^2 \widetilde{P}}{\partial \tilde{z}_{i+t} \partial \tilde{z}_{j-t}}\right\|_{D{(r,s)}\times\mathcal{O}}
\leq&    \frac{\varepsilon_0}{t}\me^{-\rho|i+j|}.
\end{split}
\end{equation*}
Similar argument also applies to  the second derivative   $\frac{\partial^2 P}{\partial \bar{z}_i \partial \bar{z}_j}.$

It follows that  $P\in \mathcal{T}^\rho_{ D(r,s)\times\mathcal{O}}$ and  $\langle P\rangle_{\rho,D(r,s)\times \mathcal{O}}\leq\varepsilon_0  .$

\end{description}

%

\subsection{The  half-wave equations}
Denote  the inner product $\langle u,v\rangle=Re\int^\pi_0u(x)\overline{v(x)}dx.$
The  half-wave equation \eqref{eq2} can be written as
\begin{equation}\label{Hhweq2}
            u_t= \mi \nabla H(t,u)=\mi \mathbf{D}_0u+\mi\varepsilon V(\omega t, x)u.\\
\end{equation}
where  the Hamiltonian
\begin{equation*}
    H(t,u)=\frac{1}{2}\langle\mathbf{D}_0u, u\rangle+\frac{\varepsilon }{2}\int^{\pi}V(\omega t, x)|u|^2 dx.
\end{equation*}

We expand $u(t, x)$ on the eigenfunctions
\,$$u(t, x)=\sum\limits_{j\geq1}q_j(t)\phi_j(x)\in \mathcal{H}^{a,p}_0,$$
(see \eqref{fB1} on the space $\mathcal{H}^{a,p}_0$) where $q=(q_j)_{j\geq1}  \in  \ell^{a,p}_0.$
Then the equation \eqref{Hhweq2} becomes
\begin{equation}
            \dot{q}_j=2\mi\frac{\partial}{\partial \bar{q}_j}H(t,q,\bar{q})=\lambda_jq_j+2\frac{\partial}{\partial \bar{q}_j}G,
\end{equation}
where
$$H(t,q, \bar{q})=\Lambda+G,$$
\begin{equation*}
 \begin{split}
    \Lambda=&\sum\limits_{j\geq1}\frac{\lambda_j}{2}q_j\bar{q}_j,\\
    G=&\frac{\varepsilon}{2}\sum\limits_{j,k\geq1} q_j\bar{q}_k \int^{\pi}_{0}  V(t\omega,x)\phi_j(x)\phi_k(x)dx.
 \end{split}
\end{equation*}
To rewrite the above  equation   as an autonomous Hamiltonian system, we
introduce the angle variables \,$\theta=\omega t\in \mathbb{T}^n,$\, the  action variables \,$I\in \mathbb{R}^n$ and
the complex coordinates \,$z=(z_j)_{j\geq1},\,\,  \bar{z}=(\bar{z}_j)_{j\geq1}$\,through
\begin{equation*}
    z_j = \frac{1}{\sqrt{2}}q_j,\,\, \bar{z}_j = \frac{1}{\sqrt{2}}\bar{q}_j.
\end{equation*}
Then we obtain an  autonomous Hamiltonian system
\begin{equation}\label{atHeq1}
      \begin{cases}
            \dot{z}_j=\mi\lambda_j z_j+\mi \frac{\partial}{\partial \bar{z}_j} P(\theta,z,\bar{z})\quad            &   j\geq 1,\\
            \dot{\bar{z}}_j=-\mi\lambda_j \bar{z}_j-\mi \frac{\partial }{\partial z_j} P(\theta,z,\bar{z})\quad    &  j\geq 1,\\
            \dot{\theta}_i=\omega_i\quad                                                                              &   i=1\cdots n,\\
            \dot{I}_i=-\frac{\partial}{\partial \theta_i}P(\theta,z,\bar{z})\quad                             & i=1\cdots n.
      \end{cases}
\end{equation}
on the phase space \,$\mathcal{P}^{a,p}_0$\
with respect to the symplectic form
\begin{equation*}
      \sum_{i=1}^nd\theta_i\wedge dI_i+\mi\sum_{j\geq 1}dz_j\wedge d\bar{z}_j.
\end{equation*}
The new Hamiltonian  associated to the system \eqref{atHeq1} is
\begin{equation}\label{HHWlam}
 \begin{split}
 H=&N+P\\
  \end{split}
\end{equation}
 where
 \begin{equation*}
 \begin{split}
 N=&\sum_{i=1}^n\omega_iI_i+\sum_{j\geq 1}\lambda_j z_j\bar{z}_j,\\
    P=&\varepsilon\sum\limits_{l,k\geq1} z_l\bar{z}_k \int^{\pi}_{0}  V(t\omega,x)\phi_l(x)\phi_k(x)dx.
 \end{split}
\end{equation*}

The next is  the verification of the  assumptions \textbf{(A1)-(A4)} for the Hamiltonian \eqref{HHWlam}.
Let $r$ be that in Assumption 1.1 and $s>0$ be a suitable positive number.
Take  $\varepsilon_0=(2^{p+1}+2^{4}+ \frac{n}{2r})C_V\epsilon>0.$

\begin{description}
  \item[(1)] \emph{Verifying  the assumption} \textbf{(A1)}.

Since $\lambda_j =j,$ then we
take   $\Omega_j=j+\breve{\Omega}_j$ with
 $\breve{\Omega}_j=0,$
thus   $ \breve{\Omega}_j\in C^1_W([0, 2\pi)^{n}).$
Let $A_0=1.$ It is obvious that  for all $j\geq 1$  and $\omega\in [0, 2\pi)^{n},$
$|\breve{\Omega}_j|\leq A_0$ and
$|\partial_\omega\breve{\Omega}_j|\leq \varepsilon_0.$

 \item[(2)] \emph{Verifying  the assumption} \textbf{(A2)}.

Following the  verification of the  assumption \textbf{(A2)} in Section \ref{DLW}, we can also prove that
 there is a subset $\mathcal{O}\subset[0, 2\pi)^{n}$ of positive Lebesgue measure
with $\mes\mathcal{O}\geq (2\pi)^{n}(1-O(\gamma))$ such that  the  assumption \textbf{(A2)} holds for \eqref{HHWlam} on $\mathcal{O}.$

\item[(3)] \emph{Verifying  the assumption} \textbf{(A3)}.

The perturbation   $P$ in \eqref{HHWlam} reads
\begin{equation}\label{hwpij2}
 \begin{split}
    P=&\varepsilon\sum\limits_{ij\geq1}p_{ij}(\theta) z_i\bar{z}_j,
 \end{split}
\end{equation}
where
 \begin{equation*}
 \begin{split}
p_{ij}(\theta):=&\int^{\pi}_{0}  V(\theta,x)\phi_i(x)\phi_j(x)dx \\
=&  \begin{cases}
           \frac{1}{2}(\widetilde{V}_{i-j}(\theta)- \widetilde{V}_{i+j}(\theta)), \quad     i>j, \\
          \widetilde{V}_{0}(\theta)  -\frac{1}{2} \widetilde{V}_{2j}(\theta),  \quad     i=j , \\
          \frac{1}{2}(\widetilde{V}_{j-i}(\theta)- \widetilde{V}_{i+j}(\theta)), \quad     i<j.
         \end{cases}\\
\end{split}
\end{equation*}

Following the arguments in the  verification of the  assumption \textbf{(A3)} for the wave equation \eqref{eq1}, one can prove  that
\begin{equation*}
 \begin{split}
\|X_P\|_{s;D(r,s)\times\mathcal {O}}
\leq(2^{p+1}+12+\frac{n}{2r})\|V\|_{D(2r),b,p}\epsilon\leq\varepsilon_0.
\end{split}
\end{equation*}
This shows
the regularity of  Hamiltonian vector field $X_P.$

\item[(4)] \emph{Verifying  the assumption} \textbf{(A4)}.

Let $\rho=2a.$
Thanks to $\breve{\Omega}_j\equiv0,$ it is obvious that
 $\breve{\Omega}:=diag(\breve{\Omega}_j)_{j\geq1}\in \mathfrak{M}^{\rho}_{r,\mathcal{O}}.$

Now we verify $P\in \mathcal{T}^\rho_{ D(r,s)\times\mathcal{O}}.$
By \eqref{hwpij2}, we have
 $$\frac{\partial^2 P}{\partial z_i \partial \bar{z}_j}=\epsilon p_{ij}(\theta)\,\,
   \hbox{and}\,\,  \frac{\partial^2 P}{\partial z_i \partial z_j}=0=\frac{\partial^2 P}{\partial \bar{z}_i \partial \bar{z}_j}.$$

 Following the arguments in verifying  the  assumption \textbf{(A4)} in Section \ref{DLW}, we have
  the limit  $\lim\limits_{t\rightarrow\infty}\frac{\partial^2 P}{\partial z_{i+t} \partial \bar{z}_{j+t}}$ exists.
  Moreover,
\begin{equation*}
  \left\|\lim\limits_{t\rightarrow\infty} \frac{\partial^2 P}{\partial z_{i+t} \partial \bar{z}_{j+t}}\right\|_{D{(r,s)}\times\mathcal{O}}
  \leq\varepsilon_0 \me^{-\rho|i-j|}
\end{equation*}
and
\begin{equation*}
\left\|\frac{\partial^2 P}{\partial z_{i+t} \partial \bar{z}_{j+t}}-\lim\limits_{t\rightarrow\infty}\frac{\partial^2 P}{\partial z_{i+t} \partial \bar{z}_{j+t}}\right\|_{D{(r,s)}\times\mathcal{O}}
\leq   \frac{\epsilon }{2}   \|V\|_{D(2r),b,p}   \me^{-b(i+j+2t)} \leq   \frac{\varepsilon_0}{t}\me^{-\rho|i-j|}.
\end{equation*}
This together with  $\frac{\partial^2 P}{\partial z_i \partial z_j}=0=\frac{\partial^2 P}{\partial \bar{z}_i \partial \bar{z}_j}$   shows that the perturbation
$P\in \mathcal{T}^\rho_{ D(r,s)\times\mathcal{O}}$
and
$\langle P\rangle_{\rho,D(r,s)\times \mathcal{O}}\leq\varepsilon_0.$

\end{description}

\section[]{Proof of the reducibility Theorem \ref{Redb}}

\subsection{Basic strategy}

The reducibility  Theorem \ref{Redb} is proved  by KAM method.
We construct a sequence of Hamiltonian $H=N+P$ of  the form \eqref{Hf}.
Suppose the perturbation $P=O(\varepsilon),$ then we construct a symplectic   coordinate transformation $\Phi$
such that it transforms $H=N+P$ into a  new Hamiltonian $H_+=H\circ\Phi=N_++P_+$
with new normal form $N_+$ and a smaller perturbation  $P_+=O(\varepsilon^\kappa),\,\, 1<\kappa<2,$ than the old perturbation $P.$

The above transformation $\Phi$ is constructed  via the flow  $X^t_F$ generated by
 a quadratic Hamiltonian $F.$ Taking  $\Phi=X^1_{F}$  and denoting $R=T_KP$, then
\begin{equation*}
  \begin{split}
     H\circ\Phi = H\circ X^1_{F}=&N\circ X^1_{F}+R\circ X^1_{F}+(P-R)\circ X^1_{F}\\
               =&N+\{N,F\}+\int^1_0(1-t)\{\{N,F\},F\}\circ X^t_{F}dt\\
               &+R+\int^1_0\{R,F\}\circ X^t_{F}dt+(P-R)\circ X^1_{F}.
   \end{split}
\end{equation*}
The new normal form is defined as $N_+=N+\hat{N}.$ This leads to the following  homological equation
\begin{equation*}
    \{N, F\}+R=\hat{N}.
\end{equation*}
where  the  unknowns  are $F$ and $\hat{N}.$
We solve this  homological equation in the next section.


\subsection{Solving the Homological equation}\label{Homeqapprsol}

Consider the  homological equation
\begin{equation}\label{vhomeq}
    \{N, F\}+R=\hat{N}
\end{equation}
on $D(r,s)\times \mathcal{O},$
where
$$ N=\sum_{i=1}^n\omega_iI_i+\sum_{j\geq 1}\Omega_j(\xi) z_j\bar{z}_j$$
with the fixed tangential frequencies  $\omega(\xi)\in \mathbb{R}^n.$ The normal  frequencies $\Omega_j(\xi)\in \mathbb{R},\,j\geq1$ satisfy \eqref{AsyNF}.
 The Hamiltonian $R$ is a quadratic  on $(z,\bar{z})$ of the  form
\begin{equation}
  \begin{split}
        R(\theta, z, \bar{z}; \xi)=&\langle R^{20}(\theta)z, z\rangle  +\langle R^{11}(\theta)z, \bar{z}\rangle+\langle R^{02}(\theta)z, \bar{z}\rangle\\
        =&\sum\limits_{|k|\leq K}\sum\limits_{i,j\geq1}[R^{20}_{kij}(\xi)z_iz_j + R^{11}_{kij}(\xi)z_i\bar{z}_j +  R^{02}_{kij}(\xi)\bar{z}_i\bar{z}_j]\me^{\mi k\cdot\theta}.\\
   \end{split}
\end{equation}
It does not depend on the action  variables $I$ and satisfies $R=\mathcal{T}_KR.$
We define its mean value $[ R]$ with respect to $\theta$  by
\begin{equation*}
  \begin{split}
        [ R]=\sum\limits_{j\geq1}R^{11}_{0jj}(\xi)z_j\bar{z}_j.
   \end{split}
\end{equation*}

In the following,  we use the notations  $$\Gamma_{ab}(\sigma)=\sup\limits_{t\geq0}(1+t)^a\Delta^b(t)\me^{-t\sigma},\,\,a,\, b\in \mathbb{N}.$$

\begin{proposition}\label{slhomeq}
  Let   $\gamma>0 $ and $0<5\sigma<r.$  Suppose  $N$ and $R$  satisfy the above conditions (A1)--(A2),
then the homological equation  \eqref{vhomeq} has the  unique  solutions $F$ and $\widehat{N}$ satisfying $[F]=0$
 and the estimates
\begin{equation}\label{F}
    \|X_F\|_{s;D(r-\sigma,s)\times \mathcal{O}}\leq (1+\mathcal{C}_0)\gamma^{-2}\Gamma_{12}(\sigma) \|X_R\|_{s;D(r,s)\times \mathcal{O}},
\end{equation}
\begin{equation}\label{N}
    \|X_{\widehat{N}}\|_{s;D(r,s)\times \mathcal{O}}\leq  \|X_R\|_{s;D(r,s)\times \mathcal{O}},
\end{equation}
where the constant $\mathcal{C}_0$ depends only on $E$ and $L.$
\end{proposition}
\begin{proof}
We look for a Hamiltonian $F$ of the form
\begin{equation}
  \begin{split}
        F(\theta, z, \bar{z}; \xi)=&\langle F^{20}(\theta)z, z\rangle  +\langle F^{11}(\theta)z, \bar{z}\rangle+\langle F^{02}(\theta)z, \bar{z}\rangle\\
        =&\sum\limits_{|k|\leq K}\sum\limits_{i,j\geq1}[F^{20}_{kij}(\xi)z_iz_j + F^{11}_{kij}(\xi)z_i\bar{z}_j +  F^{02}_{kij}(\xi)\bar{z}_i\bar{z}_j]\me^{\mi k\cdot\theta}.\\
   \end{split}
\end{equation}

Denote $\omega\cdot\nabla f(\theta):=\sum^n_{b=1}\omega_b \frac{\partial f}{\partial \theta_b}.$
We take  $\hat{N}=[R].$
By the comparison of coefficients,
the homological equation \eqref{vhomeq} is  equivalent to the following scalar  form: For all $i,j\geq1,$
\begin{equation}\label{sheq1}
   \omega\cdot\nabla  F^{20}_{ij}+\mi (\Omega_i+\Omega_j) F^{20}_{ij}= R^{20}_{ij},
\end{equation}
\begin{equation}\label{sheq2}
   \omega\cdot\nabla  F^{11}_{ij}+\mi (\Omega_i-\Omega_j) F^{11}_{ij}= R^{11}_{ij}-\delta_{ij}[R^{11}_{ij}],
\end{equation}
and
\begin{equation}\label{sheq3}
   \omega\cdot\nabla  F^{02}_{ij}-\mi (\Omega_i+\Omega_j) F^{02}_{ij}= R^{02}_{ij},
\end{equation}
here $\delta_{ij}=1, \hbox{if}\ i=j,\hbox{and}\ 0, \hbox{otherwise}.$

Consider the  equation \eqref{sheq2}.
For $i=j,$ the equation \eqref{sheq2} becomes
\begin{equation}\label{sheq21}
   \partial_{\omega} F^{11}_{jj}= R^{11}_{jj}-[R^{11}_{jj}],
\end{equation}
then by Fourier expansion,
 \begin{equation*}
     F^{11}_{kjj}=
                       \begin{cases}
                        0,\quad                 &  k=0,\\
                       \frac{R^{11}_{kjj}}{\mi k\cdot \omega},\quad                   & 0<|k|\leq K.
                       \end{cases}
\end{equation*}
and we obtain the form solution
 \begin{equation*}
     F^{11}_{jj}=\sum\limits_{0<|k|\leq K}\frac{R^{11}_{kjj}}{\mi k\cdot \omega}\me^{\mi k\cdot \theta}.
\end{equation*}

For $i\neq j,$ by Fourier expansion, the equation \eqref{sheq2} becomes
 \begin{equation*}
     F^{11}_{kij}=\frac{R^{11}_{kij}}{\mi (k\cdot \omega + \Omega_i-\Omega_j)}
\end{equation*}
and we obtain the form solution
 \begin{equation}\label{Fsol}
     F^{11}_{ij}(\theta)=\sum\limits_{0\leq|k|\leq K}\frac{R^{11}_{kij}}{\mi (k\cdot \omega + \Omega_i-\Omega_j)}\me^{\mi k\cdot \theta}.
\end{equation}

Now we give the estimate for $ F^{11}_{ij}.$
Denote $S_{ij}=k\cdot \omega + \Omega_{i}- \Omega_{j}.$
For all $1\leq a\leq n,$
\begin{equation*}
\begin{split}
    |\partial_{\xi_a} S_{i,j}|=|k\cdot \partial_{\xi_a}\omega + \partial_{\xi_a}\breve{\Omega}_{i}- \partial_{\xi_a}\breve{\Omega}_{j}|
   \leq \mathcal{C}_0(1+|k|),
\end{split}
\end{equation*}
where the constant $\mathcal{C}_0=\mathcal{C}_0(E,L)$ depends only on $E$ and $L.$

Then
\begin{equation}\label{F11ij}
\begin{split}
     &\|F^{11}_{ij}\|_{D(r-\sigma)\times \mathcal{O}}\\
    \leq&\sum\limits_{|k|\leq K}\left(\frac{|R^{11}_{k,ij}|_{\mathcal{O}}}{| S_{ij}|}+\frac{| \partial_{\xi}S_{ij}||R^{11}_{k,ij}|_{\mathcal{O}}}{| S^2_{ij}|} \right)\me^{|k|(r-\sigma)}\\
    \leq&\sum\limits_{|k|\leq K}(1+\mathcal{C}_0)(1+|k|)\gamma^{-2}\Delta^2(|k|)\me^{-|k|\sigma}   |R^{11}_{k,ij}|_{\mathcal{O}}\me^{|k|r}\\
  \leq& (1+\mathcal{C}_0)\gamma^{-2}\Gamma_{12}(\sigma)\|R^{11}_{ij}\|_{D(r)\times \mathcal{O}} .
\end{split}
\end{equation}
Similarly, we have
\begin{equation*}
     \|F^{20}_{ij}\|_{D(r-\sigma)\times \mathcal{O}}
    \leq(1+\mathcal{C}_0)\gamma^{-2}\Gamma_{12}(\sigma) \|R^{20}_{ij}\|_{D(r)\times \mathcal{O}},
\end{equation*}
\begin{equation*}
     \|F^{02}_{ij}\|_{D(r-\sigma)\times \mathcal{O}}
    \leq(1+\mathcal{C}_0)\gamma^{-2}\Gamma_{12}(\sigma) \|R^{02}_{ij}\|_{D(r)\times \mathcal{O}}.
\end{equation*}

Note that the derivative
\begin{equation}
        \frac{\partial F}{\partial z_i}
        =\sum\limits_{j\geq1}F^{20}_{ji}z_j +  F^{20}_{ij}z_j+ F^{11}_{ij}\bar{z}_j,
\end{equation}
then
\begin{equation}\label{Fzder}
\lfloor \frac{\partial F}{\partial z_i}\rfloor_{D(r-\sigma)\times\mathcal {O}}
        \leq     (1+\mathcal{C}_0)\gamma^{-2}\Gamma_{12} (\sigma) \lfloor \frac{\partial R}{\partial z_i}\rfloor_{D(r-\sigma)\times\mathcal {O}}.
\end{equation}
Similarly,
\begin{equation}\label{Fzbarder}
\lfloor \frac{\partial F}{\partial \bar{z}_i}\rfloor_{D(r-\sigma)\times\mathcal {O}}
        \leq    (1+\mathcal{C}_0)\gamma^{-2}\Gamma_{12} (\sigma) \lfloor \frac{\partial R}{\partial \bar{z}_i}\rfloor_{D(r-\sigma)\times\mathcal {O}}.
\end{equation}

 For each  $1\leq b\leq n,$  by \eqref{Fsol},   the norm of the derivative $\frac{\partial F^{11}_{ij}}{\partial \theta_b}$ is
\begin{equation*}
\begin{split}
     \|\frac{\partial F^{11}_{ij}}{\partial \theta_b}\|_{D(r-\sigma)\times \mathcal{O}}
    \leq (1+\mathcal{C}_0)\gamma^{-2}\Gamma_{12} (\sigma) \|\frac{\partial R^{11}_{ij}}{\partial \theta_b}\|_{D(r)\times \mathcal{O}}.
\end{split}
\end{equation*}
Similarly, we have
$$\|\frac{\partial F^{20}_{ij}}{\partial \theta_b}\|_{D(r-\sigma)\times \mathcal{O}}\leq (1+\mathcal{C}_0)\gamma^{-2}\Gamma_{12} (\sigma) \|\frac{\partial R^{20}_{ij}}{\partial \theta_b}\|_{D(r)\times \mathcal{O}},$$
and
$$\|\frac{\partial F^{02}_{ij}}{\partial \theta_b}\|_{D(r-\sigma)\times \mathcal{O}}\leq (1+\mathcal{C}_0)\gamma^{-2}\Gamma_{12} (\sigma) \|\frac{\partial R^{02}_{ij}}{\partial \theta_b}\|_{D(r)\times \mathcal{O}}.$$
It follows  that
\begin{equation}\label{Fthetader}
  \|\frac{\partial F}{\partial \theta}\|_{D(r-\sigma,s)\times \mathcal{O}}\leq (1+\mathcal{C}_0)\gamma^{-2}\Gamma_{12} (\sigma) \|\frac{\partial R}{\partial \theta}\|_{D(r,s)\times \mathcal{O}}
\end{equation}

From  \eqref{Fzder}, \eqref{Fzbarder} and \eqref{Fthetader},
   we obtain the estimate for the Hamiltonian vector field $X_F:$
\begin{equation*}
  \|X_F\|_{s;D(r-\sigma,s)\times \mathcal{O}}\leq (1+\mathcal{C}_0)\gamma^{-2}\Gamma_{12} (\sigma) \|X_R\|_{s;D(r,s)\times \mathcal{O}}.
\end{equation*}

The estimates of $X_{\hat{N}}$ follow from the observation that $\hat{N}_{z\bar{z}}$ is the diagonal of the mean value of $R_{z\bar{z}}.$

\end{proof}

The above lemma implies the  estimate for the Jacobian $DX_F:$
\begin{equation}\label{DfJacob}
\begin{split}
  \|DX_F\|_{s;D(r-2\sigma,s)\times \mathcal{O}}  \leq  C\sigma^{-1}(1+\mathcal{C}_0)\gamma^{-2}\Gamma_{12} (\sigma) \|X_R\|_{s;D(r,s)\times \mathcal{O}}.
\end{split}
\end{equation}

Now we verify the T\"{o}plitz-Lipschitz property of the solutions of homological equation   \eqref{vhomeq}.
\begin{proposition}\label{TL}
Suppose  $N$ and $R$  satisfy the above conditions (A1)--(A2) and   $R\in \mathcal{T}^\rho_{D(r,s)\times\mathcal{O}},$
   then there exists a constant  $C:= 5+4\mathcal{C}_0$ such that for any $0<\sigma<r,$ the solutions $F$ and $\hat{N}$   of homological equation   \eqref{vhomeq}  are  T\"{o}plitz-Lipschitz on $D(r,s)\times\mathcal{O},$ i.e.,  $F\in \mathcal{T}^\rho_{D(r-\sigma,s)\times\mathcal{O}},\,\,\hat{N}\in \mathcal{T}^\rho_{D(r,s)\times\mathcal{O}},$  and
\begin{equation}\label{FTL}
  \langle F\rangle_{\rho,D(r-\sigma,s)\times\mathcal{O}} \leq  C\gamma^{-3}\Gamma_{13} (\sigma) \langle R\rangle_{\rho, D(r,s)\times\mathcal{O}},
\end{equation}
\begin{equation}\label{NTL}
  \langle \hat{N}\rangle_{\rho,D(r,s)\times\mathcal{O}} \leq    \langle R\rangle_{\rho, D(r,s)\times\mathcal{O}} .
\end{equation}
\end{proposition}
 \begin{proof}
 The estimation  of $\hat{N}$ follows from the observation that $\hat{N}_{z\bar{z}}$ is the diagonal of the mean value of $R_{z\bar{z}}.$
 In the following, we prove  the estimation  \eqref{FTL}.

From   \eqref{Fsol} in the proof of Lemma \ref{slhomeq}, the second derivative of $F$ w.r.t. $ z_i,\,\,\bar{z}_j$ is
 \begin{equation*}
     \frac{\partial^2 F}{\partial z_i \partial \bar{z}_j}=F^{11}_{ij}(\theta)=\sum\limits_{0\leq|k|\leq K}\frac{R^{11}_{kij}}{\mi (k\cdot \omega + \Omega_i-\Omega_j)}\me^{\mi k\cdot \theta}.
\end{equation*}

$\bullet$ We first verify the  exponentially off-diagonal decay of $\frac{\partial^2 F}{\partial z_i \partial \bar{z}_j}.$

Since $R\in \mathcal{T}^\rho_{D(r,s)\times\mathcal{O}},$ we have
\begin{equation*}
\left\|\frac{\partial^2 R}{\partial z_{i}\partial \bar{z}_{j}}\right\|_{D(r,s)\times\mathcal{O}}
\leq \langle R\rangle_{\rho,D(r,s)\times \mathcal{O}}\me^{-\rho|i-j|}.
\end{equation*}
Then
\begin{equation*}
\begin{split}
     \left\| \frac{\partial^2 F}{\partial z_i \partial \bar{z}_j}\right\|_{D(r-\sigma)\times \mathcal{O}}
     \leq&\sum\limits_{|k|\leq K}\gamma^{-1}\Delta(|k|)\me^{-|k|\sigma}   |R^{11}_{k,ij}|_{\mathcal{O}}\me^{|k|r}\\
    \leq&\gamma^{-1}\Gamma_{11}(\sigma) \left\|\frac{\partial^2 R}{\partial z_{i}\partial \bar{z}_{j}}\right\|_{D(r,s)\times\mathcal{O}}\\
  \leq& \gamma^{-1}\Gamma_{11}(\sigma) \langle R\rangle_{\rho, D(r,s)\times\mathcal{O}}\me^{-\rho|i-j|} .
\end{split}
\end{equation*}

$\bullet$We then check  the  asymptotically T\"{o}plitz property of $\frac{\partial^2 F}{\partial z_i \partial \bar{z}_j}.$

Since
$\Omega_j=j+\breve{\Omega}_j,\,\,j\geq1,$ and  $\langle\langle \breve{\Omega}\rangle\rangle_{\rho,r,\mathcal{O}}<\varepsilon_0,$
the limits  $\lim_{t\rightarrow\infty}\breve{\Omega}_{j+t}$ exist and satisfy
\begin{equation}\label{tl0}
\left\|\lim_{t\rightarrow\infty}\breve{\Omega}_{j+t}\right\|_{\mathcal{O}}
\leq \varepsilon_0,
\end{equation}
\begin{equation}\label{tl00}
\left\|\breve{\Omega}_{j+t}-\lim_{t\rightarrow\infty}\breve{\Omega}_{j+t}\right\|_{\mathcal{O}}\leq \frac{\varepsilon_0}{|t|}.
\end{equation}
Note that
 $$\Omega_{i+t}-\Omega_{j+t}=i-j+\breve{\Omega}_{i+t}-\breve{\Omega}_{j+t},$$
then for all $i, j$ the limits  $\Omega_{i,j,\infty}:=\lim_{t\rightarrow\infty}(\Omega_{i+t}-\Omega_{j+t})$ exist and satisfy
the non-resonance conditions
\begin{equation}\label{tnonre}
| k\cdot\omega+ \Omega_{i,j,\infty}|\geq\frac{\gamma}{\Delta(|k|)}.
\end{equation}

Denote  $S_{ij,\infty}:=k\cdot \omega + \Omega_{i,j,\infty}.$
For  $1\leq a\leq n, $
\begin{equation*}
\begin{split}
    |\partial_{\xi_a} S_{ij,\infty}|=&|k\cdot \partial_{\xi_a}\omega + \partial_{\xi_a}\breve{\Omega}_{i,\infty}- \partial_{\xi_a}\breve{\Omega}_{j,\infty}|\\
  \leq& |k||\omega|_{\mathcal{O}}+2|\breve{\Omega}|_{\mathcal{O}}
   \leq \mathcal{C}_0(1+|k|),
\end{split}
\end{equation*}
where the constant $\mathcal{C}_0=\mathcal{C}_0(E,L).$

Since $R\in \mathcal{T}^\rho_{D(r,s)\times\mathcal{O}},$
the limit  $R^{11}_{ij,\infty}:=\lim_{t\rightarrow\infty}R^{11}_{i+t,j+t}.$ exists.
Consider a similar equation to   the  equation \eqref{sheq2}:
\begin{equation*}
   \partial_{\omega} u+\mi \Omega_{i,j,\infty}u= R^{11}_{ij,\infty}.
\end{equation*}
By the non-resonance conditions \eqref{tnonre}, the solution $F^{11}_{ij,\infty}$ of the above  equation exists:
\begin{equation}\label{Finfty}
  F^{11}_{ij,\infty}=\sum\limits_{0\leq|k|\leq K}\frac{R^{11}_{k,ij,\infty}}{\mi S_{ij,\infty}}\me^{\mi k\cdot \theta}.
\end{equation}
Moreover, similar to the  estimate for $\|F^{11}_{ij}\|_{D(r-\sigma)\times \mathcal{O}}$ in \eqref{F11ij}, we obtain
\begin{equation*}
\begin{split}
     \|F^{11}_{ij,\infty}\|_{D(r-\sigma)\times \mathcal{O}}
  \leq (1+\mathcal{C}_0)\gamma^{-2}\Gamma_{12}\langle R\rangle_{\rho, D(r,s)\times\mathcal{O}}\me^{-\rho|i-j|}  ,
\end{split}
\end{equation*}
thus
$$\|\lim_{t\rightarrow\infty}\frac{\partial^2 F}{\partial z_{i+t} \partial \bar{z}_{j+t}}\|_{D(r-\sigma)\times \mathcal{O}}\leq (1+\mathcal{C}_0)\gamma^{-2}\Gamma_{12}(\sigma)\langle R\rangle_{\rho, D(r,s)\times\mathcal{O}}\me^{-\rho|i-j|} .$$

$\bullet$ Finally,  we check that   $\frac{\partial^2 F}{\partial z_i \partial \bar{z}_j}$ is Lipschitz at infinity.

By \eqref{Fsol} and \eqref{Finfty}, we write  the difference $ F^{11}_{i+t,j+t}-F^{11}_{ij,\infty}$  as
\begin{equation*}
\begin{split}
    F^{11}_{i+t,j+t}-F^{11}_{ij,\infty}=&\sum\limits_{|k|\leq K}\mathcal{F}_{k,ij}(\xi) \me^{k\cdot\theta}.
\end{split}
\end{equation*}
where
\begin{equation*}
\begin{split}
    \mi \mathcal{F}_{k,ij}(\xi)=\frac{R^{11}_{k,i+t,j+t}}{S_{i+t,j+t}}
    - \frac{ R^{11}_{k,ij,\infty}}{S_{ij,\infty}} .
\end{split}
\end{equation*}
For  $ a=1,\cdots, n,$
the Whitney derivatives of $\mathcal{F}_{k,ij}(\xi)$ with respect to $\xi_a$ are
\begin{equation*}
\begin{split}
  \mi \partial_{\xi_a} \mathcal{F}_{k,ij}(\xi)
 =&\frac{\partial_{\xi_a}( R^{11}_{k,i+t,j+t}-R^{11}_{k,ij,\infty})}{S_{i+t,j+t}}- \frac{\partial_{\xi_a}S_{i+t,j+t}}{S^2_{i+t,j+t}}( R^{11}_{k,i+t,j+t}-R^{11}_{k,ij,\infty})\\
   & +\left(\frac{S_{ij,\infty}-S_{i+t,j+t}}{S_{i+t,j+t}S_{ij,\infty}} \right)   \partial_{\xi_a}  R^{11}_{k,ij,\infty} -\left(\frac{\partial_{\xi_a} S_{i+t,j+t}-\partial_{\xi_a} S_{ij,\infty}}{S^2_{i+t,j+t}} \right)   R^{11}_{k,ij,\infty}\\
    &-\partial_{\xi_a} S_{ij,\infty}  (S_{ij,\infty}-S_{i+t,j+t})
    \left(\frac{1}{S^2_{i+t,j+t} S_{ij,\infty}} + \frac{1}{S_{i+t,j+t} S^2_{ij,\infty}}\right)   R^{11}_{k,ij,\infty}.
\end{split}
\end{equation*}

In view of  $\langle\langle \breve{\Omega}\rangle\rangle_{\rho,r,\mathcal{O}}<\varepsilon_0,$ $|\omega|_{\mathcal{O}}\leq E$ and $|\widetilde{\Omega}|_{\mathcal{O}}\leq L,$ we have
\begin{equation*}
\begin{split}
    |S_{ij,\infty}-S_{i+t,j+t}|
   \leq 2|t|^{-1}\varepsilon_0.
\end{split}
\end{equation*}
and for  $ a=1,\cdots, n,$
\begin{equation*}
\begin{split}
    |\partial_{\xi_a} S_{ij,\infty}-\partial_{\xi_a} S_{i+t,j+t}|
   \leq 2|t|^{-1}\varepsilon_0.
\end{split}
\end{equation*}
It follows that
\begin{equation*}
\begin{split}
    | \mathcal{F}_{k,ij}(\xi)|
 \leq&  \gamma^{-1}\Delta(|k|)|R^{11}_{k,i+t,j+t}-R^{11}_{k,ij,\infty}|
 + \gamma^{-2}\Delta^2(|k|)2|t|^{-1}\varepsilon_0 |R^{11}_{k,ij,\infty}|\\
\end{split}
\end{equation*}
and
\begin{equation*}
\begin{split}
   &| \partial_{\xi_a}  \mathcal{F}_{k,ij}(\xi) |\\
\leq&   \Delta(|k|) \gamma^{-1} |\partial_{\xi_a}( R^{11}_{k,i+t,j+t}-R^{11}_{k,ij,\infty})|+\mathcal{C}_0(1+|k|)\Delta^2(|k|) \gamma^{-2}  |R^{11}_{k,i+t,j+t}-R^{11}_{k,ij,\infty}|\\
&+2|t|^{-1}\varepsilon_0\Delta^2(|k|)  \gamma^{-2} ( | \partial_{\xi_a}  R^{11}_{k,ij,\infty}|+  | R^{11}_{k,ij,\infty}|)+4\mathcal{C}_0(1+|k|) |t|^{-1}\varepsilon_0\Delta^3(|k|)\gamma^{-3} | R^{11}_{k,ij,\infty}|.
\end{split}
\end{equation*}
Therefore,
\begin{equation*}
\begin{split}
     \| F^{11}_{i+t,j+t}-F^{11}_{ij,\infty}\|_{D(r-\sigma)\times \mathcal{O}}
    =&\sum\limits_{|k|\leq K}| \mathcal{F}_{k,ij}(\xi) |_{\mathcal{O}}\me^{|k|(r-\sigma)}\\
 \leq& (2\gamma^{-1}\Gamma_{01}+\mathcal{C}_0\gamma^{-2}\Gamma_{12})\|R^{11}_{i+t,j+t}-R^{11}_{ij,\infty}\|_{D(r)\times \mathcal{O}}\\
 &+( 5\gamma^{-2}\Gamma_{02}+ 4\mathcal{C}_0\gamma^{-3}\Gamma_{13})|t|^{-1}\|R^{11}_{ij,\infty}\|_{D(r)\times \mathcal{O}}.
\end{split}
\end{equation*}
This together with  $R\in \mathcal{T}^\rho_{D(r,s)\times\mathcal{O}}$ shows  that
\begin{equation*}
\begin{split}
     &\left\| \frac{\partial^2 F}{\partial z_{i+t} \partial \bar{z}_{j+t}}-\lim_{t\rightarrow\infty}\frac{\partial^2 F}{\partial z_{i+t} \partial \bar{z}_{j+t}}\right\|_{D(r-\sigma)\times \mathcal{O}}\\
 \leq& ( 5+4\mathcal{C}_0)\gamma^{-3}\Gamma_{13}  |t|^{-1}\langle R\rangle_{\rho, D(r,s)\times\mathcal{O}}\me^{-\rho|i-j|}.
\end{split}
\end{equation*}

Similarly, we have
\begin{equation*}
\begin{split}
     &\| \lim_{t\rightarrow\infty}\frac{\partial^2 F}{\partial z_{i+t} \partial z_{j-t}}\|_{D(r-\sigma)\times \mathcal{O}},\,\| \lim_{t\rightarrow\infty}\frac{\partial^2 F}{\partial \bar{z}_{i+t} \partial \bar{z}_{j-t}}\|_{D(r-\sigma)\times \mathcal{O}}\\
     &\leq(1+\mathcal{C}_0)\gamma^{-2}\Gamma_{12}(\sigma) \langle R\rangle_{\rho, D(r,s)\times\mathcal{O}}\me^{-\rho|i+j|}
\end{split}
\end{equation*}
and
\begin{equation*}
\begin{split}
     &\left\| \frac{\partial^2 F}{\partial \bar{z}_{i+t} \partial \bar{z}_{j-t}}-\lim_{t\rightarrow\infty}\frac{\partial^2 F}{\partial \bar{z}_{i+t} \partial \bar{z}_{j-t}}\right\|_{D(r-\sigma)\times \mathcal{O}},\,\,
     \left\| \frac{\partial^2 F}{\partial z_{i+t} \partial z_{j-t}}-\lim_{t\rightarrow\infty}\frac{\partial^2 F}{\partial z_{i+t} \partial z_{j-t}}\right\|_{D(r-\sigma)\times \mathcal{O}}\\
     \leq&( 5+4\mathcal{C}_0)\gamma^{-3}\Gamma_{13}  |t|^{-1}\langle R\rangle_{\rho, D(r,s)\times\mathcal{O}}\me^{-\rho|i+j|}.
\end{split}
\end{equation*}
Thus we complete the proof of the estimation \eqref{FTL}.

\end{proof}

\subsection{KAM iteration  and convergence}

Let $C_*$ be a constant that is twice the maximum of all implicit constants used during the KAM step,
and it depends only on $n, E, L$ and $\rho_0.$

We take the Hamiltonian $H=N+P$ in \eqref{Hf} as the initial Hamiltonian  $H_0=N_0+P_0$ .
Similarly, we set other initial  quantities  as those in Section \ref{kamthem}. Namely,
we set $r_0=r,\, s_0=s,\, \gamma_0=\gamma,\,   \rho_0=\rho,\, K_0=K,\,\, \mathcal{O}_{0}=\mathcal{O}.$

For $\nu\geq0$,
$$\gamma_\nu=\frac{\gamma_0}{2}(1+2^{-\nu}),$$
$$\delta_\nu=2^{-(\nu+4)}\rho_{0},\,\,\rho_{\nu+1}=\rho_{\nu}-4\delta_\nu.$$

Denote $$\Gamma(\sigma)=\Gamma_{23}(\sigma)=\sup\limits_{t\geq0}(1+t)^2 \Delta^3(t) \me^{-\sigma t}.$$

 $$\kappa_\nu=\kappa^{-(\nu+1)},\,\,\kappa=\frac{4}{3}.$$

Given $\sigma>0$ with $6\sigma<r_0.$ There exists a non-increasing positive sequence
 $$\sigma_0 \geq \sigma_1 \geq \sigma_2 \geq \cdots \geq\sigma_\nu \geq \sigma_{\nu+1}\geq\cdots>0 $$
such that
\begin{equation}\label{sigma}
  \sum\limits^\infty_{\nu=0}\sigma_\nu=\sigma
\end{equation}
and
\begin{equation}\label{Xi}
  \Xi(\sigma)=\inf\limits_{\tilde{\sigma}_0 \geq \tilde{\sigma}_1 \geq \cdots >0,\atop\,\tilde{\sigma}_0 +\tilde{\sigma}_1+\cdots\leq \sigma}\prod\limits^{\infty}_{\mu=0}\Gamma^{\kappa_\mu}(\tilde{\sigma}_\mu)=\prod\limits^{\infty}_{\mu=0}\Gamma^{\kappa_\mu}(\sigma_\mu)<\infty,
\end{equation}
see Appendix for the proof.
For such a fixed sequence $\{\sigma_\nu\}$,  we define
\begin{equation}\label{Gamman}
  \Gamma_\nu=2C_* \Gamma(\sigma_\nu),
\end{equation}
and
\begin{equation}\label{varepsilon1}
 \varepsilon_{\nu+1}=\Gamma_\nu \varepsilon^{\kappa}_{\nu},
\end{equation}
then
\begin{equation}\label{varepsilon2}
 \varepsilon_{\nu}=\left(\prod\limits^{\nu-1}_{\mu=0}\Gamma^{\kappa_\mu}_\mu \varepsilon_0\right)^{\kappa^\nu},\,\, \nu\geq1.
\end{equation}

 The order  $K_\nu$  of  Fourier truncation is defined implicitly by
\begin{equation}\label{Knu}
C_*  \me^{-K_\nu\sigma_\nu}=\Gamma_\nu \varepsilon^{1/2}_{\nu}.
\end{equation}

Finally, we set $$r_\nu=r_0-3\sum\limits^{\nu-1}_{\mu=0}\sigma_\mu,\,\,s_{\nu+1}=\frac{1}{4}s_{\nu}$$
and denote the domain  $D_\nu=D(r_\nu, s_\nu).$ It is obvious that
$$D_{0} \supseteq D_{1}  \supseteq D_\nu \supseteq \cdots,\,\,\hbox{and}\,\,\mathcal{O}_{0} \supseteq \mathcal{O}_{1}  \supseteq \mathcal{O}_\nu \supseteq \cdots.$$

\begin{lemma}[Iterative Lemma]\label{IteLem}
Let $0<\varepsilon_0<\min \{(C_*\gamma_02^5)^{\frac{3}{2}}, \delta^{12}_{0}, (\gamma_{0}\delta_{0})^{9/2}\}.$
Given a sequence of parameter domains $$\mathcal{O}_{0} \supseteq \mathcal{O}_{1} \supseteq \cdots \supseteq \mathcal{O}_\nu .$$
Suppose for $\ell=0,1,\cdots,\nu,$  the  Hamiltonian  $H_\ell=N_\ell+P_\ell$ are regular  on $D_\ell\times \mathcal{O}_\ell,$ where
the normal forms
\begin{equation}\label{Nell}
  N_\ell=\sum_{j=1}^n\omega_jI_j+\sum_{j\in \mathbb{Z}}\Omega_{\ell,j}(\omega) z_j\bar{z}_j
\end{equation}
 with  $\Omega_{\ell,j}(\omega)=j+\widetilde{\Omega}_{\ell,j}(\omega)$ satisfies
\begin{equation}\label{widetOmega}
  |\widetilde{\Omega}_{\ell}|_{\mathcal{O}}\leq A_0+\sum\limits^{\ell-1}_{b=1}\varepsilon_b\,\,\hbox{and}\,\, \langle\langle \widetilde{\Omega}_{\ell}\rangle\rangle_{\rho_\ell,r_\ell,\mathcal{O}_\ell}\leq \varepsilon_0+\sum\limits^{\ell-1}_{b=1}\varepsilon_b,
\end{equation}
\begin{equation}\label{iteranonresoncond}
  \begin{split}
   &| k\cdot \omega|\geq\frac{\gamma_\ell}{\Delta(|k|)},\,\forall 0<|k|\leq K_\ell,\\
   &|k\cdot \omega+\Omega_{\ell,i}(\omega)+\Omega_{\ell,j}(\omega)|\geq\frac{\gamma_\ell}{\Delta(|k|)},\,\forall |k|\leq K_\ell,\, i,j\geq 1, \\
   &|k\cdot \omega+\Omega_{\ell,i}(\omega)-\Omega_{\ell,j}(\omega)|\geq\frac{\gamma_\ell}{\Delta(|k|)},\,\forall |k|\leq K_\ell,\, i\neq j, \\
   \end{split}
\end{equation}
on $\mathcal{O}_\ell,$  and the perturbation $P_\ell$  satisfies
\begin{equation}\label{Pell}
  P_\ell\in \mathcal{T}^{\rho_\ell}_{D_\ell\times\mathcal{O}_\ell}\,\, \hbox{and}\,\,  [ P_\ell]^{\rho_\ell}_{s_\ell;D_\ell\times\mathcal{O}_\ell} <\varepsilon_{\ell}.
\end{equation}
Then there exists a Whitney smooth  family of real analytic symplectic transformations
$\Phi_{\nu+1}: D_{\nu+1} \times \mathcal{O}_{\nu} \rightarrow D_\nu$ satisfying
\begin{equation}\label{nearidenty}
    \|\Phi_{\nu+1}-id\|_{s_\nu;D_{\nu+1}\times \mathcal{O}_\nu},\,\,\|D\Phi_{\nu+1}-I\|_{s_\nu;D_{\nu+1}\times \mathcal{O}_\nu} \leq \varepsilon^{5/6}_\nu,
\end{equation}
and a closed subset of $\mathcal{O}_\nu:$
\begin{equation}\label{newnonresonset}
    \mathcal{O}_{\nu+1}=\mathcal{O}_\nu \setminus \bigcup_{|k|>K_\nu}\left(\mathcal{R}^{\nu+1}_{k}(\gamma_{\nu+1}) \cup \bigcup_{i,j}\mathcal{R}^{+,\nu+1}_{kij}(\gamma_{\nu+1})\cup \bigcup_{i\neq j}\mathcal{R}^{-,\nu+1}_{kij}(\gamma_{\nu+1})\right),
\end{equation}
where
\begin{equation*}
    \mathcal{R}^{\nu+1}_{k}(\gamma_{\nu+1})=\left\{\omega\in \mathcal{O}_\nu:|k\cdot \omega|< \frac{\gamma_{\nu+1}}{\Delta(|k|)}\right\},
\end{equation*}
\begin{equation*}
    \mathcal{R}^{+,\nu+1}_{kij}(\gamma_{\nu+1})=\left\{\omega\in \mathcal{O}_\nu:|k\cdot \omega+\Omega_{\nu+1,i}(\omega)+\Omega_{\nu+1,j}(\omega)|<\frac{\gamma_{\nu+1}}{\Delta(|k|)}\right\},
\end{equation*}
\begin{equation*}
    \mathcal{R}^{-,\nu+1}_{kij}(\gamma_{\nu+1})=\left\{\omega\in \mathcal{O}_\nu:|k\cdot \omega+\Omega_{\nu+1,i}(\omega)-\Omega_{\nu+1,j}(\omega)|<\frac{\gamma_{\nu+1}}{\Delta(|k|)}\right\},
\end{equation*}
such that $\Phi_{\nu+1}$ transforms  $H_\nu$ into
$$H_{\nu+1}=H_\nu\circ \Phi_{\nu+1}=N_{\nu+1}+P_{\nu+1},$$
and on the domain $D_{\nu+1} \times \mathcal{O}_{\nu+1},$  $N_{\nu+1}$ and $P_{\nu+1}$
satisfy the  conditions  $\eqref{Nell}_{\nu+1},$  $\eqref{widetOmega}_{\nu+1},$   $\eqref{iteranonresoncond}_{\nu+1}$ and $\eqref{Pell}_{\nu+1}.$

\end{lemma}
\begin{proof}

$\blacklozenge$ \emph{The construction of  symplectic transformation} $\Phi_{\nu+1}.$

Let $R_{\nu}=T_{K_\nu}P_{\nu}$ be the Fourier truncation of order $K_{\nu}$ of $P_{\nu}.$
Using the inequalities
$$\|X_{R_{\nu}}\|_{s_{\nu};D_{\nu}\times \mathcal{O}_{\nu}}\leq \|X_{P_{\nu}}\|_{s_{\nu};D_{\nu}\times \mathcal{O}_{\nu}}\leq \varepsilon_\nu,$$
$$\langle R_{\nu}\rangle_{\rho_{\nu}, D_{\nu}\times\mathcal{O}_{\nu}}\leq \langle P_{\nu}\rangle_{\rho_{\nu}, D_{\nu}\times\mathcal{O}_{\nu}}  \leq \varepsilon_\nu,$$
and by Propositions \ref{slhomeq} and \ref{TL}, under the non-resonance conditions $\eqref{iteranonresoncond}_{\nu},$
the  homological equation
\begin{equation}\label{itehomeq}
    \{N_{\nu}, F\}+R_{\nu}=\hat{N}
\end{equation}
has a set of unique solutions $F=F_{\nu}$ and $\hat{N}=\hat{N}_{\nu}$ satisfying the estimates
\begin{equation}\label{Fite}
    \|X_{F_{\nu}}\|_{s_{\nu};D(r_{\nu}-\sigma_{\nu},s_{\nu})\times \mathcal{O}_{\nu}}\leq C\gamma^{-2}_{\nu}\Gamma_{12} (\sigma_{\nu})\|X_{R_{\nu}}\|_{s_{\nu};D_{\nu}\times \mathcal{O}_{\nu}}  \leq C\gamma^{-2}_{\nu}\Gamma_{12} (\sigma_{\nu})\varepsilon_\nu,
\end{equation}
\begin{equation}\label{Nite}
    \|X_{\widehat{N}_{\nu}}\|_{s_{\nu};D_{\nu}\times \mathcal{O}_{\nu}}\leq  \|X_{R_{\nu}}\|_{s_{\nu};D_{\nu}\times \mathcal{O}_{\nu}}\leq \varepsilon_{\nu},
\end{equation}
\begin{equation}\label{FTLite}
  \langle F_{\nu}\rangle_{\rho_{\nu},D(r_{\nu}-\sigma_{\nu},s_{\nu})\times\mathcal{O}_{\nu}} \leq  C\gamma^{-3}_{\nu}\Gamma_{13}(\sigma_{\nu})  \langle R_{\nu}\rangle_{\rho_{\nu}, D(r_{\nu},s_{\nu})\times\mathcal{O}_{\nu}}\leq  C\gamma^{-3}_{\nu}\Gamma_{13}(\sigma_{\nu})  \varepsilon_\nu,
\end{equation}
and
\begin{equation}\label{NTLite}
  \langle \hat{N}_{\nu}\rangle_{\rho_{\nu},D_{\nu}\times\mathcal{O}_{\nu}} \leq    \langle R_{\nu}\rangle_{\rho_{\nu}, D_{\nu}\times\mathcal{O}_{\nu}} \leq \varepsilon_\nu.
\end{equation}
Since $\Gamma_{12}\leq\Gamma_{13}$ by the definition and $\Gamma_{13}\leq  \sigma\Gamma_{23}$ by Lemma ? in Appendix,
 we have
\begin{equation}\label{Fplusite}
 \begin{split}
 [ F_\nu]^{\rho_\nu}_{s_\nu;D(r_{\nu}-\sigma_{\nu},s_{\nu})\times\mathcal{O}_\nu}=& \|X_{F_{\nu}}\|_{s_{\nu};D(r_{\nu}-\sigma_{\nu},s_{\nu})\times \mathcal{O}_{\nu}}+\langle F_{\nu}\rangle_{\rho_{\nu},D(r_{\nu}-\sigma_{\nu},s_{\nu})\times\mathcal{O}_{\nu}}\\
   \leq&C\gamma^{-2}_{\nu}\Gamma_{12} (\sigma_{\nu})\varepsilon_\nu+C\gamma^{-3}_{\nu}\Gamma_{13}(\sigma_{\nu})  \varepsilon_\nu\\
   \leq& C\sigma_{\nu}.\\
  \end{split}
\end{equation}
Then  by Lemma \ref{Flowest},  the flow $X^t_{F_{\nu}}$ generated by the  Hamiltonian vector field $X_{F_{\nu}}$ exists  on $D(r_{\nu}-\sigma_{\nu},\frac{s_{\nu}}{4})$ for all $0\leq t\leq1.$
Taking  $\Phi_{\nu+1}=X^1_{F_{\nu}},$
it maps $D(r_{\nu}-\sigma_{\nu},\frac{s_{\nu}}{4})$ into $ D(r_{\nu},\frac{s_{\nu}}{2}).$

Now we prove the estimate \eqref{nearidenty}.
\begin{equation}
 \begin{split}
\|\Phi_{\nu+1}-id\|_{s_\nu;D_{\nu+1}\times \mathcal{O}_\nu}\leq& 2 \|X_{F_{\nu}}\|_{s_{\nu};D(r_{\nu}-\sigma_{\nu},s_{\nu})\times \mathcal{O}_{\nu}}\\
\leq&2C\gamma^{-2}_{\nu}\Gamma_{12} (\sigma_{\nu})\varepsilon_\nu\\
   \leq&\varepsilon^{\kappa-1/2}_\nu=\varepsilon^{5/6}_\nu.
  \end{split}
\end{equation}
\begin{equation}
 \begin{split}
\|D\Phi_{\nu+1}-I\|_{s_\nu;D_{\nu+1}\times \mathcal{O}_\nu}\leq& 2 \|DX_{F_{\nu}}\|_{s_{\nu};D(r_{\nu}-\sigma_{\nu},s_{\nu})\times \mathcal{O}_{\nu}}\\
\leq&\sigma^{-1}_{\nu}2C\gamma^{-2}_{\nu}\Gamma_{12} (\sigma_{\nu})\varepsilon_\nu\\
   \leq& C\frac{\varepsilon_{\nu+1}}{\varepsilon_\nu}\\
   \leq&\varepsilon^{\kappa-1/2}_\nu=\varepsilon^{5/6}_\nu.
  \end{split}
\end{equation}

$\blacklozenge$ \emph{The new Hamiltonian} $H_{\nu+1}.$

Using the Taylor formula together with the   homological equation \eqref{itehomeq}, we define the new Hamiltonian
\begin{equation}
 \begin{split}
H_{\nu+1}=H_\nu\circ \Phi_{\nu+1}=&N_\nu\circ \Phi_{\nu+1}+R_\nu\circ \Phi_{\nu+1}+(P_\nu-R_\nu)\circ \Phi_{\nu+1}\\
               =&N_\nu+\{N_\nu,F_\nu\}+\int^1_0(1-t)\{\{N_\nu, F_\nu\},F_\nu\}\circ X^t_{F_\nu}dt\\
               &+R_\nu+\int^1_0\{R_\nu, F_\nu\}\circ X^t_{F_\nu}dt+(P_\nu-R_\nu)\circ X^1_{F_\nu}\\
=&N_{\nu+1}+P_{\nu+1},
  \end{split}
\end{equation}
where
$$N_{\nu+1}=N_{\nu}+\hat{N}_{\nu},$$
$$P_{\nu+1}=\int^1_0\{\widehat{R}_{\nu}(t),F_{\nu}\}\circ X^t_{F_{\nu}}dt+(P_{\nu}-R_{\nu})\circ X^1_{F_{\nu}}$$
with $\widehat{R}_{\nu}(t)=(1-t)\hat{N}_{\nu}+tR_\nu.$

$\bullet$ \emph{The estimation for} $P_{\nu+1}.$

We first consider the  estimation for $\|X_{P_{\nu+1}}\|_{s_{\nu+1};D_{\nu+1}\times \mathcal{O}_{\nu+1}}.$
Note that $$X_{P_{\nu+1}}=\int^1_0 (X^t_{F_{\nu}})^*[X_{\widehat{R}_{\nu}(t)}, X_{F_{\nu}}] dt+(X^1_{F_{\nu}})^*(X_{P_{\nu}}-X_{R_{\nu}}) .$$
Then
\begin{equation}
 \begin{split}
&\|X_{P_{\nu+1}}\|_{s_{\nu+1};D_{\nu+1}\times \mathcal{O}_{\nu+1}}\\
  \leq&\int^1_0 \|(X^t_{F_{\nu}})^*[X_{\widehat{R}_{\nu}(t)}, X_{F_{\nu}}] \|_{s_{\nu+1};D_{\nu+1}\times \mathcal{O}_{\nu+1}}dt\\
     &+\|(X^1_{F_{\nu}})^*(X_{P_{\nu}}-X_{R_{\nu}})\|_{s_{\nu+1};D_{\nu+1}\times \mathcal{O}_{\nu+1}}\\
   \leq&2\|[X_{R_{\nu}}, X_{F_{\nu}}] \|_{s_{\nu+1};D(r_{\nu}-2\sigma_{\nu}, s_{\nu+1})\times \mathcal{O}_{\nu+1}}
                +2\|X_{P_{\nu}}-X_{R_{\nu}}\|_{s_{\nu+1};D(r_{\nu}-\sigma_{\nu}, s_{\nu+1})\times \mathcal{O}_{\nu+1}}\\
       \leq&2C\sigma^{-1}_{\nu}\|X_{R_{\nu}} \|_{s_{\nu};D_\nu\times \mathcal{O}_{\nu}} \|X_{F_{\nu}}\|_{s_{\nu};D_\nu\times \mathcal{O}_{\nu}}
                +2\me^{-K_{\nu}\sigma_{\nu}}\|X_{P_{\nu}}\|_{s_{\nu};D_\nu\times \mathcal{O}_{\nu}}\\
            \leq&C\sigma^{-1}_{\nu} \gamma^{-2}_{\nu}   \Gamma_{12}(\sigma_{\nu}) \varepsilon^{2}_{\nu}   +     2\me^{-K_{\nu}\sigma_{\nu}}    \varepsilon_{\nu}.\\
  \end{split}
\end{equation}

We then consider  the  estimation for  $\langle P_{\nu+1}\rangle^{\rho_{\nu+1}}_{D_{\nu+1}\times \mathcal{O}_{\nu+1}}.$
\begin{equation}
 \begin{split}
&\langle P_{\nu+1}\rangle^{\rho_{\nu+1}}_{D_{\nu+1}\times \mathcal{O}_{\nu+1}}\\
  \leq&\int^1_0 \langle \{\widehat{R}_{\nu}(t),F_{\nu}\}\circ X^t_{F_{\nu}}\rangle^{\rho_{\nu+1}}_{D_{\nu+1}\times \mathcal{O}_{\nu+1}}dt\\
     &+\langle (P_{\nu}-R_{\nu})\circ X^1_{F_{\nu}}\rangle^{\rho_{\nu+1}}_{D_{\nu+1}\times \mathcal{O}_{\nu+1}}\\
   \leq&C\delta^{-2}_{\nu}\int^1_0 \langle \{\widehat{R}_{\nu}(t),F_{\nu}\}\rangle^{\rho_{\nu}-\delta_{\nu}}_{D(r_{\nu}-2\sigma_{\nu}, s_{\nu})\times \mathcal{O}_{\nu+1}}dt\\
     &+C\delta^{-2}_{\nu}\langle P_{\nu}-R_{\nu}\rangle^{\rho_{\nu}-\delta_{\nu}}_{D(r_{\nu}-2\sigma_{\nu}, s_{\nu})\times \mathcal{O}_{\nu+1}}\\
   \leq&C \delta^{-3}_{\nu}\gamma^{-3}_{\nu}   \Gamma_{13}(\sigma_{\nu}) \varepsilon^{2}_{\nu}   +   C\delta^{-2}_{\nu}  \me^{-K_{\nu}\sigma_{\nu}}    \varepsilon_{\nu}.\\
  \end{split}
\end{equation}
It follows that
\begin{equation}
 \begin{split}
 &[P_{\nu+1}]^{\rho_{\nu+1}}_{s_{\nu+1};D_{\nu+1}\times\mathcal{O}_{\nu+1}} \\
   \leq&C\delta^{-3}_{\nu}\sigma^{-1}_{\nu} \gamma^{-3}_{\nu}   \Gamma_{13}(\sigma_{\nu}) \varepsilon^{2}_{\nu}   +    C\delta^{-2}_{\nu} \me^{-K_{\nu}\sigma_{\nu}}    \varepsilon_{\nu}\\
      \leq&C  \Gamma_{23}(\sigma_{\nu}) \varepsilon^{\kappa}_{\nu}   +     C\delta^{-2}_{\nu}  \me^{-K_{\nu}\sigma_{\nu}}    \varepsilon_{\nu}\\
      \leq& \Gamma_{\nu} \varepsilon^{\kappa}_{\nu}=\varepsilon_{\nu+1}.
  \end{split}
\end{equation}

$\bullet$ \emph{The new frequency and non-resonance condition.}

In the new normal form $N_{\nu+1},$  the  frequencies $\Omega_{\nu+1,j}=j+\breve{\Omega}_{\nu+1,j}=\Omega_{\nu,j}+\widehat{\Omega}_{\nu,j},$ where
$\widehat{\Omega}_{\nu,j}=\frac{\partial^2 \widehat{N}_{\nu}}{\partial z_j \partial \bar{z}_j}.$
Thus $$|\widehat{\Omega}_{\nu,j}|_{\mathcal{O}}\leq \|X_{\widehat{N}_{\nu}}\|_{s_{\nu};D_{\nu}\times \mathcal{O}_{\nu}}\leq  \|X_{R_{\nu}}\|_{s_{\nu};D_{\nu}\times \mathcal{O}_{\nu}}\leq \varepsilon_{\nu}.$$

\begin{equation}
 \begin{split}
 |\lim_{t\rightarrow\infty}\widehat{\Omega}_{\nu,j+t}|_{\mathcal{O}_{\nu}}\leq \|\lim_{t\rightarrow\infty}\frac{\partial^2 \widehat{N}_{\nu}}{\partial z_{j+t} \partial \bar{z}_{j+t}}\|_{s_{\nu};D_{\nu}\times \mathcal{O}_{\nu}}
 \leq\langle \hat{N}_{\nu}\rangle_{\rho_{\nu},D_{\nu}\times\mathcal{O}_{\nu}} \leq \varepsilon_\nu.
  \end{split}
\end{equation}
\begin{equation}
 \begin{split}
 |\widehat{\Omega}_{\nu,j+t}-\lim_{t\rightarrow\infty}\widehat{\Omega}_{\nu,j+t}|_{\mathcal{O}_{\nu}}\leq &\|\frac{\partial^2 \widehat{N}_{\nu}}{\partial z_{j+t} \partial \bar{z}_{j+t}}-\lim_{t\rightarrow\infty}\frac{\partial^2 \widehat{N}_{\nu}}{\partial z_{j+t} \partial \bar{z}_{j+t}}\|_{s_{\nu};D_{\nu}\times \mathcal{O}_{\nu}}\\
 \leq&|t|^{-1}\langle \hat{N}_{\nu}\rangle_{\rho_{\nu},D_{\nu}\times\mathcal{O}_{\nu}} \leq |t|^{-1}\varepsilon_\nu.
  \end{split}
\end{equation}
These imply
$$|\widehat{\Omega}_{\nu}|_{\mathcal{O}}\leq  \varepsilon_{\nu},\,\,  \langle\langle \widehat{\Omega}_{\nu}\rangle\rangle_{\rho_{\nu},r_{\nu},\mathcal{O}_{\nu}}\leq \varepsilon_{\nu}.$$
Therefore,

\begin{equation*}
  |\breve{\Omega}_{\nu+1}|_{\mathcal{O}}\leq A_0+\sum\limits^{\nu}_{b=1}\varepsilon_b\,\,\hbox{and}\,\, \langle\langle \breve{\Omega}_{\nu+1}\rangle\rangle_{\rho_{\nu},r_{\nu},\mathcal{O}_{\nu}}\leq \varepsilon_0+\sum\limits^{\nu}_{b=1}\varepsilon_b,
\end{equation*}

$\blacklozenge$Finally, we consider the construction of $\mathcal{O}_{\nu+1}$.
It suffices to verify
\begin{equation*}
  \begin{split}
   |k\cdot \omega+\Omega_{\nu+1,i}(\omega)-\Omega_{\nu+1,j}(\omega)|\geq\frac{\gamma_\ell}{\Delta(|k|)},\,\forall |k|\leq K_{\nu},\, i\neq j.
   \end{split}
\end{equation*}

By the definition of $\gamma_\nu,$ $\Gamma_\nu$ and $K_\nu,$  we have
\begin{equation*}
  \begin{split}
   \frac{\gamma_0}{2^{\nu+3} \varepsilon_\nu \Delta(K_\nu)}=& \frac{\gamma_0\me^{-K_\nu\sigma_\nu}}{2^{\nu+3} \varepsilon_\nu \Delta(K_\nu)\me^{-K_\nu\sigma_\nu}}\\
   \geq&  \frac{\gamma_0\me^{-K_\nu\sigma_\nu}}{2^{\nu+3} \varepsilon_\nu \Gamma(\sigma_\nu)}\\
   =&\frac{2\gamma_0  \varepsilon^{\kappa-5/6}_{\nu}}{2^{\nu+3} \varepsilon_\nu }\\
      =&\frac{\gamma_0  }{2^{\nu+2} \varepsilon^{1/2}_{\nu} }\geq1.\\
   \end{split}
\end{equation*}
This implies $\gamma_{\nu}-\gamma_{\nu+1}\geq 2 \varepsilon_\nu \Delta(|k|) $ for all $0<|k|\leq K_\nu,$ thus
\begin{equation*}
  \begin{split}
   |k\cdot \omega+\Omega_{\nu+1,i}(\omega)-\Omega_{\nu+1,j}(\omega)|\geq& |k\cdot \omega+\Omega_{\nu,i}(\omega)-\Omega_{\nu,j}(\omega)|- |\widehat{\Omega}_{\nu,i}(\omega)|-|\widehat{\Omega}_{\nu,j}(\omega)|\\
   \geq& \frac{\gamma_\nu}{\Delta(|k|)}-2\varepsilon_{\nu}
   \geq\frac{\gamma_{\nu+1}}{\Delta(|k|)}.
   \end{split}
\end{equation*}
Then after removing the resonance zones for $K_\nu<|k|\leq K_{\nu+1},$  we get a closed set $\mathcal{O}_{\nu+1}\subseteq\mathcal{O}_{\nu}$ with the desired properties.

\end{proof}

\emph{The Convergence Proof.}

By the iterative Lemma \ref{IteLem}, we obtain  a sequence of decreasing domains $D_\nu \times \mathcal{O}_\nu$ and symplectic  transformations
$\Phi^{\nu}=\Phi_{1}\circ \Phi_{2} \circ \cdots \circ \Phi_{\nu}:D_\nu \times \mathcal{O}_{\nu-1} \rightarrow D_{\nu-1},\,\nu\geq 1.$
Then by \eqref{nearidenty} and following the arguments in \cite{Poschel96a}, the sequence $\Phi^{\nu}$ of symplectic  transformations  converge uniformly on
$D(r/2) \times \mathcal{O}_\gamma $ to a real analytic torus embedding $\Phi:\mathbb{T}^n\rightarrow \mathcal{P}^{a,p},$    for which
we also need to   verify
\begin{description}
  \item[(a)] the symplectic coordinate transformation  $\Phi$ is of the form given in \eqref{lintransform};
  \item[(b)] the new Hamiltonian eventually  reduces to the new normal form, i.e., $P^\infty=0;$
  \item[(c)]  the symplectic coordinate transformation $\Phi$, which is defined by Theorem \ref{Redb} on
each $\mathcal{P}^{a,p}$, extends to $\mathcal{P}^{a,0}.$
\end{description}
In fact, by \eqref{FeqZsol} and \eqref{FeqIsol} in Section \ref{FeqIsol}, the  the symplectic coordinate transformation  $\Phi_{\nu}$  at the $\nu^{th}-$step
has the form
the form
\begin{equation}\label{lintransform2}
  \Phi_{\nu}\left(
    \begin{array}{c}
     \theta \\
     I\\
      Z \\
    \end{array}
  \right)
=
\left(
  \begin{array}{c}
   \theta  \\
    \Phi^{(I)}_{\nu}\\
    \Phi^{(Z)}_{\nu} \\
  \end{array}
\right)=
\left(
  \begin{array}{c}
   \theta  \\
   I+\frac{1}{2}Z^TM_{\nu}(\theta)Z\\
    L_{\nu}(\theta)Z \\
  \end{array}
\right).
\end{equation}
In particular, the linear operator  $L_{\nu}(\theta)=\me^{JA_{\nu}(\theta)}$ is invertible.
Then property (\textbf{a}) is satisfied at each step, and thus  we can iterate
the process.  It follows that  the limiting
 transformation $\Phi=\Phi_{1}\circ \Phi_{2} \circ \cdots$
also satisfies the property (\textbf{a}).
 Similar to the initial Hamiltonian, the transformed  Hamiltonian
is linear in $I$ and quadratic in $Z,$  which implies that the new Hamiltonian eventually  reduces to the new normal form, i.e., $P^\infty=0.$

Since $ \Phi^{(Z)}$ is a linear symplectomorphism, then following Prop.1.3 (\cite{Kuksin00}) by duality, it extends on $\ell^{a,p}\times\ell^{a,p}$ for all $p\in [-2, 2]$ and thus the conclusion (\textbf{c})  holds  if we take $p=0.$

The sequence of closed subset $\mathcal{O}_\nu$  converges to a closed set
 $$\mathcal{O}_\gamma=\bigcap_{\nu\geq0}\mathcal{O}_\nu.$$
By the construction of $\gamma_\nu$ and $|\Omega_{\nu+1} - \Omega_{\nu}|_{\mathcal{O}}=|\widehat{\Omega}_{\nu}|_{\mathcal{O}}\leq \varepsilon_{\nu},$
we have $|\Omega^{\infty} - \Omega|_{\mathcal{O}}\leq \varepsilon^{1/2}_{0},$ and thus
for all $\omega\in \mathcal{O}_\gamma,$
$$|\langle k, \omega\rangle|\geq\frac{\gamma}{2\Delta(|k|)},\,\forall k\neq0,$$
$$|\langle k, \omega\rangle+l\cdot \Omega^{\infty}(\omega)|\geq\frac{\gamma}{2\Delta(|k|)},\,\forall k\in \mathbb{Z}^n,\, |l|=2.$$

The measure estimate of $\mathcal{O} \setminus \mathcal{O}_\gamma$ of bad frequencies is given in the next section.

\subsection{Measure estimate}

In this subsection, we complete the  Lebesgue measure estimate   of the parameter set $\mathcal{O} \setminus \mathcal{O}_\gamma.$
 In the process of constructing iterative sequences,   we obtain a decreasing sequence of closed sets $\mathcal{O}_0\supset \mathcal{O}_1\supset \cdots$ such that $\mathcal{O}_\gamma=\bigcap\limits_{\nu \geq 0}\mathcal{O}_\nu$ and
 \begin{equation}\label{newresonset}
    \mathcal{O} \setminus \mathcal{O}_\gamma=\bigcup_{\nu\geq0} \bigcup_{K_{\nu-1}<|k|\leq K_{\nu},i,j }\left(\mathcal{R}^{\nu}_{k}(\gamma_{\nu}) \cup \bigcup_{i,j}\mathcal{R}^{+,\nu}_{kij}(\gamma_{\nu})\cup \bigcup_{i\neq j}\mathcal{R}^{-,\nu}_{kij}(\gamma_{\nu})\right),
\end{equation}
where  $\mathcal{R}^{\nu}_{k},  \mathcal{R}^{+,\nu}_{kij},  \mathcal{R}^{-,\nu}_{kij}$
are defined in \eqref{newnonresonset}.

Below we only consider the most difficult resonance set $\mathcal{R}^{-,\nu}_{kij}(\gamma_{\nu}).$
The proof for other resonance sets $\mathcal{R}^{\nu}_{k},\,\,\mathcal{R}^{+,\nu}_{kij}$
are more simple, and thus omitted.

Since   $\Omega_{\nu,j}=j+\breve{\Omega}_{\nu,j},$ then  by \eqref{widetOmega}, there is a constant $A_1>0$ such that
$|\Omega_{\nu,i}-\Omega_{\nu,j}|\geq A_1|i-j|.$ Denote $A_2=(1+2A_1+2A_0)/ A_1.$
Note that when $|i-j|>A_2|k|,$
 $$|k\cdot \omega +\Omega_{\nu,i}-\Omega_{\nu,j}|\geq (1+A_0+A_1)|k|,$$
thus in this case there is no small divisor, and in the following it remains to consider the case of $|i-j|\leq A_2|k|.$

Denote $$S^{\nu}_{k,i,j}=k\cdot \omega + \Omega_{\nu,i}- \Omega_{\nu,j},$$
 $$S^{\nu}_{k,i,j,\infty}=k\cdot \omega + \lim_{t\rightarrow\infty}(\Omega_{\nu,i+t}- \Omega_{\nu,j+t})$$
and introduce the following resonant sets
\begin{equation*}
    \mathcal{R}^{-,\nu}_{k,i+t,j+t}(\gamma_{\nu})=\left\{\omega\in \mathcal{O}_{\nu-1}:|S^{\nu}_{k,i+t,j+t}|<\frac{\gamma_{\nu}}{\Delta(|k|)}\right\},
\end{equation*}

\begin{lemma}\label{mea2}
For $i,j\geq 1$ with $|i-j|\leq  A_2|k|$,
  there exist  $i',\,\, j'\geq 1$ satisfying  $i'\leq 2A_2|k|,\,\, j'\leq 2A_2|k|$ and $t\geq 1$
such that $i=i'+ t,\,\, j=j'+ t.$  Consequently,
\begin{equation}\label{Rijt}
  \bigcup\limits_{i,j,\atop\,|i-j|\leq  A_2|k|} \mathcal{R}^{\nu}_{kij}\subset \bigcup\limits_{i', j'\leq 2A_2|k|}\bigcup\limits_{t\geq1} \mathcal{R}^{\nu}_{k,i'+t,j'+t}.
\end{equation}
\end{lemma}
\begin{proof}
Without loss of generalization, we assume $j>i.$
  For given $i,\,\, j,$ choosing  a $t_0\geq 1$ such that $0\leq i-t_0\leq A_2|k|.$
Let $i'=i-t_0$ and  $j'=i'+j-i=j-t_0$, then $$j'\leq i'+|j-i|\leq 2A_2|k|.$$
It follows that  \eqref{Rijt} holds.
\end{proof}

\begin{lemma}\label{mea3}
For fixed $k,  i', j',$
$$ \mes\left(\bigcup\limits_{ t\geq1} \mathcal{R}^{\nu}_{k,i'+t,j'+t}\right)\leq (20+8B_0)(2\pi)^{n(n-1)}\frac{\sqrt{\gamma}}{|k|^2\sqrt{\Delta(|k|)}}.$$
\end{lemma}
\begin{proof}
For
$\omega \in\bigcup\limits_{t>\sqrt{\frac{\Delta(|k|)}{\gamma}}}\mathcal{R}^{-,\nu}_{k,i'+t,j'+t}(\gamma_{\nu}),$ suppose
$\omega \in \mathcal{R}^{-,\nu}_{k,i'+t_0,j'+t_0}(\gamma_{\nu})$ for some  $t_0>\sqrt{\frac{\Delta(|k|)}{\gamma}}.$

From the T\"{o}plitz-Lipschitz property of $P_\nu$ and $\breve{\Omega}_\nu$, we conclude  that
$$|S^{\nu}_{k,i'+t,j'+t}-S^{\nu}_{k,i,j,\infty}|<\frac{2(1+B_0)}{|t|}.$$ Then
 \begin{equation*}
\begin{split}
 |S^{\nu}_{k,i',j',\infty}|
  \leq&|S^{\nu}_{k,i'+t_0,j'+t_0}|+|S^{\nu}_{k,i'+t_0,j'+t_0}-S^{\nu}_{k,i',j',\infty}|\\
  \leq& \frac{\gamma_{\nu}}{\Delta(|k|)}+\frac{2(1+B_0)}{|t_0|}\leq (3+2B_0)\frac{\sqrt{\gamma}}{\sqrt{\Delta(|k|)}}.
\end{split}
\end{equation*}
Thus
 \begin{equation*}
\begin{split}
\bigcup\limits_{t>\sqrt{\frac{\Delta(|k|)}{\gamma}}}\mathcal{R}^{-,\nu}_{k,i+t,j+t}(\gamma_{\nu})
\subseteq&\left\{\omega\in \mathcal{O}_{\nu-1}:|S^{\nu}_{k,i',j',\infty}|<(3+2B_0)\frac{\sqrt{\gamma}}{\sqrt{\Delta(|k|)}}\right\}\\
=:&\mathcal{Q}^{\nu}_{k,i',j',\infty}.
\end{split}
\end{equation*}

We give the estimate of $\mathcal{Q}^{\nu}_{k,i',j',\infty}.$
Taking  the  vector $v=|k|(\sgn(k_1),\cdots,\sgn(k_n)),$ then $k\cdot v=|k|^2.$ Let $\omega=\omega_s=sv+w$ with $s\in \mathbb{R},$
$w\in v^{\perp}.$
Let  $$f(s)=S^{\nu}_{k,i,j,\infty}=k\cdot \omega_s + \lim_{t\rightarrow\infty}(\Omega_{\nu,i+t}(\omega_s)- \Omega_{\nu,j+t}(\omega_s)).$$
Due to  $\sup_{\omega\in \mathcal{O}}|\lim_{t\rightarrow\infty}\partial_\omega\widetilde{\Omega}_{\nu,i+t}|\leq 3\varepsilon_0$
and $\varepsilon_0\leq \frac{1}{12n},$ the derivative
\begin{equation}\label{derfs}
\begin{split}
 | f'(s)|=&||k|^2+ \lim_{t\rightarrow\infty}v\cdot(\partial_\omega\widetilde{\Omega}_{\nu,i+t}(\omega_s)- \partial_\omega\widetilde{\Omega}_{\nu,j+t}(\omega_s))|\\
 \geq&|k|^2-6n|k|\varepsilon_0\\
 \geq&\frac{1}{2}|k|^2.
 \end{split}
\end{equation}
Then by Lemma \ref{measest}, one has
$$\mes\{s:sv+w\in \mathcal{O}_{\nu-1}, |f(s)|\leq \delta\}\leq \frac{4\delta}{|k|^2}.$$
It follows that, by Fubini's theorem,
\begin{equation}\label{Qmeas}
\begin{split}
 &\mes\left(\mathcal{Q}^{\nu}_{k,i',j',\infty}\right)\\
 \leq &\diam^{n-1}(\mathcal{O}_{\nu-1})\mes\{s:sv+w\in \mathcal{O}_{\nu-1}, |f(s)|\leq (3+2B_0)\frac{\sqrt{\gamma}}{\sqrt{\Delta(|k|)}}\}\\
 \leq &  4(2\pi)^{n(n-1)}(3+2B_0)\frac{\sqrt{\gamma}}{|k|^2\sqrt{\Delta(|k|)}}.
 \end{split}
\end{equation}

Similarly, for  the resonant set
$\mathcal{R}^{-,\nu}_{k,i'+t,j'+t},$
following the argument of estimating  $\mes\left(\mathcal{Q}^{\nu}_{k,i',j',\infty}\right),$ we have
\begin{equation}\label{Rfinmeas}
\begin{split}
 \mes\left(\mathcal{R}^{-,\nu}_{k,i'+t,j'+t}\right)
 \leq &  (2\pi)^{n(n-1)}\frac{4\gamma_{\nu}}{|k|^2\Delta(|k|)}.
 \end{split}
\end{equation}
Then
\begin{equation}\label{Rtfinmeas}
\begin{split}
 \mes\left(\bigcup\limits_{t\leq\sqrt{\frac{\Delta(|k|)}{\gamma}}}\mathcal{R}^{-,\nu}_{k,i'+t,j'+t}\right)
 \leq& 2\sqrt{\frac{\Delta(|k|)}{\gamma}} (2\pi)^{n(n-1)}\frac{4\gamma}{|k|^2\Delta(|k|)}\\
 \leq & 8(2\pi)^{n(n-1)}\frac{\sqrt{\gamma}}{|k|^2\sqrt{\Delta(|k|)}}.
 \end{split}
\end{equation}

Using  \eqref{Qmeas}  and  \eqref{Rtfinmeas}, we complete the proof.

\end{proof}

Finally, we give the estimate of  $\mes\left(\mathcal{O} \setminus \mathcal{O}_\gamma\right).$
\begin{lemma}\label{mea4}
 Let $\Delta$ be an approximation function satisfying \eqref{appfun},i.e., $ \sum\limits_{k\in \mathbb{Z}^n}\frac{1}{\sqrt{\Delta(|k|)}}< \infty.$
 Then the total measure of resonant set should be excluded during the KAM iteration is
$$ \mes\left(\mathcal{O} \setminus \mathcal{O}_\gamma\right)=O(\sqrt{\gamma}),$$
where the implicit constants  in $``O"$ depend only on $n, A_2, B_0, \Delta$ and  are made explicit in the proof.
\end{lemma}
\begin{proof}

By Lemma \ref{mea3},
\begin{equation*}
\begin{split}
 &\mes\left(\bigcup\limits_{1\leq i', j'\leq 2A_2|k|}\bigcup\limits_{ t\geq1} \mathcal{R}^{\nu}_{k,i'+t,j'+t}\right)\\
 \leq &(2A_2|k|)^2 (20+8B_0)(2\pi)^{n(n-1)}\frac{\sqrt{\gamma}}{|k|^2\sqrt{\Delta(|k|)}}\\
 \leq & A^2_2 (80+32B_0)(2\pi)^{n(n-1)}\frac{\sqrt{\gamma}}{\sqrt{\Delta(|k|)}}.
\end{split}
\end{equation*}
Then
\begin{equation*}
\begin{split}
   \mes\left(\bigcup_{\nu\geq0} \bigcup_{K_{\nu-1}<|k|\leq K_{\nu} }\bigcup\limits_{i,j} \mathcal{R}^{\nu}_{kij}\right)
  \leq& \sum\limits_{\nu\geq0}\sum_{K_{\nu-1}<|k|\leq K_{\nu} }\mes\left(\bigcup\limits_{|i'|, |j'|\leq 2A_2|k|}\bigcup\limits_{|t|\geq1} \mathcal{R}^{\nu}_{k,i'+t,j'+t}\right)\\
  \leq&  A^2_2 (80+32B_0)(2\pi)^{n(n-1)} \sum\limits_{\nu\geq0}\sum_{K_{\nu-1}<|k|\leq K_{\nu} }\frac{\sqrt{\gamma}}{\sqrt{\Delta(|k|)}}\\
   \leq&  A^2_2 (80+32B_0)(2\pi)^{n(n-1)} \sqrt{\gamma} \sum\limits_{k}\frac{1}{\sqrt{\Delta(|k|)}}.
\end{split}
\end{equation*}

Consequently,  the measure of the set $\mathcal{O} \setminus \mathcal{O}_\gamma $ is
\begin{equation*}
   \mathcal{O} \setminus \mathcal{O}_\gamma=O(\sqrt{\gamma}).
\end{equation*}

\end{proof}

\begin{remark}
Below we list three  typical approximation functions:
$\Delta_1=\exp(t^\alpha/\alpha),\,\,0<\alpha<1,$ $\Delta_2=\exp\left(\frac{t}{1+\log^\alpha(1+t)}\right),\,\,\alpha>1$ and  $ \Delta_3=\exp\left(\frac{t}{\log^\alpha t}\right),\,\,\alpha>1.$

\end{remark}

\section{Acknowledgments}
The author wishes  to thank  Prof. Jiansheng Geng  for valuable comments and  suggestions.
The research was  supported by  the National Natural Science Foundation of China (NSFC) (Grant No. 11901291)  and the Natural Science Foundation of Jiangsu Province, China (Grant No. BK20190395).

\section{Appendix}

\subsection{Some properties of  approximation functions}

%
%


\begin{lemma}\label{proAP}
For all integers $k\geq1,$ $l\geq 0,$
 $$\Gamma_{k}(\sigma)\leq  \sigma^l \Gamma_{k+l}(\sigma),\,\,$$
where $\Gamma_{k}(\sigma)=\Gamma_{k3}(\sigma)=\sup\limits_{t\geq0}(1+t)^k\Delta^3(t)\me^{-t\sigma}.$
\end{lemma}
\begin{proof}
Let $$f(t)=k\log(1+t)+\log\Delta^3(t)-t\sigma.$$
Its derivative  $$f^{\prime}(t)=\frac{k}{1+t}+\frac{d}{dt}\log\Delta^3(t)-\sigma.$$

 If $\sigma(1+t)\leq 1,$ then
$$f^{\prime}(t)\geq  \frac{k-1}{t+1}+\frac{d}{dt}\log\Delta^3(t)\geq 0.$$
It follows that $(1+t)^k\Delta^3(t)\me^{-t\sigma}$ arrive at its supremum at some point $t_*$
with $\sigma(1+t_*)\geq 1.$ Therefore, for all $l\geq 0,$
\begin{equation*}
 \begin{split}
 \Gamma_{k}(\sigma)= & (1+t_*)^k  \Delta^3(t_*)\me^{-t_*\sigma}\\
 \leq&\sigma^l(1+t_*)^{k+l}  \Delta^3(t_*)\me^{-t_*\sigma}\leq  \sigma^l \Gamma_{k+l}(\sigma).
  \end{split}
\end{equation*}

\end{proof}

Recall $\Xi$ defined in \eqref{Xi}:
\begin{equation*}
  \Xi(\sigma)=\inf\limits_{\tilde{\sigma}_0 \geq \tilde{\sigma}_1 \geq \cdots >0,\atop\,\tilde{\sigma}_0 +\tilde{\sigma}_1+\cdots\leq \sigma}\prod\limits^{\infty}_{\mu=0}\Gamma^{\kappa_\mu}(\tilde{\sigma}_\mu)=\prod\limits^{\infty}_{\mu=0}\Gamma^{\kappa_\mu}(\sigma_\mu)<\infty,
\end{equation*}
where $\Gamma(\sigma)=\Gamma_{2}(\sigma)=\Gamma_{23}(\sigma)=\sup\limits_{t\geq0}(1+t)^2\Delta^3(t)\me^{-t\sigma}.$

\begin{lemma}
  The $\Xi$ defined in \eqref{Xi} is finite for all $\sigma>0.$ In particular, let $T>0,$ if
  \begin{equation*}
  \frac{1}{\log \kappa}\int\limits^{\infty}_{T}\frac{\log\Delta(t)}{t^2}dt<\sigma,
\end{equation*}
then $$ \Xi(\sigma)\leq \me^{\sigma T}.$$
\end{lemma}
\begin{proof}
Let  $\delta(t)=\log(1+t)^2\Delta^3(t)$   and
$$t_\nu=\kappa^{\nu+1}T,\,\, \sigma_\nu=\frac{\delta(t_\nu)}{t_\nu}$$
for $\nu\geq0.$
By condition \eqref{decay1} and the hypotheses, $\sigma_0 \geq \sigma_1 \geq   \cdots>0$ and
$$\sum\limits^\infty_{\nu=0}\sigma_\nu \leq  \int\limits^{\infty}_{-1}\frac{\delta(t_\nu)}{t_\nu}d\nu\leq \frac{1}{\log \kappa}\int\limits^{\infty}_{T}\frac{\delta(t)}{t^2}dt\leq\sigma. $$
Since $\delta(t)-\sigma_\nu t\leq0$ for $t\geq t_\nu,$  then by condition \eqref{decay1} again the supremum of $\delta(t)-\sigma_\nu t$ is obtained
on the interval $[0, t_\nu]$ and thus  smaller than $\delta(t_\nu).$
It follows that
$$\Gamma(\sigma_\nu)=\sup\limits_{t\geq0}\me^{\delta(t)-\sigma_\nu t}\leq \me^{\delta(t_\nu)}=\me^{\sigma_\nu t_\nu}$$
in view of the definition of $\sigma_\nu$ and hence by $\kappa_\nu t_\nu =  T,$
\begin{equation*}
  \Xi(\sigma)\leq\prod\limits^{\infty}_{\mu=0}\me^{\kappa_\nu \sigma_\nu t_\nu}\leq \me^{\sigma T}.
\end{equation*}
\end{proof}

\begin{lemma}
$$ \sum\limits_{k\in \mathbb{Z}^n}\frac{1}{\sqrt{\Delta(|k|)}}\leq 2^n \int\limits^{\infty}_{0}\left(
                                                                                                \begin{array}{c}
                                                                                                  n+t \\
                                                                                                  n \\
                                                                                                \end{array}
                                                                                              \right)
\frac{d\log\sqrt{\Delta(t)}}{\sqrt{\Delta(t)}}dt$$
 provided that $t^n/\Delta(t)$ as $t\rightarrow0.$
\end{lemma}
\begin{proof}

Note that $$\sum\limits_{k\in \mathbb{Z}^n}\frac{1}{\sqrt{\Delta(|k|)}}\leq 2^n   \sum\limits_{k\in \mathbb{N}^n}\frac{1}{\sqrt{\Delta(|k|)}}.$$
Let $V_n(t)=\card\{k\in \mathbb{N}^n: |k|\leq t\}.$ Then
By the monotonicity of approximation  functions  the sum above may be written as a Stieltjes integral
\begin{equation*}
\begin{split}
  \sum\limits_{k\in \mathbb{N}^n}\frac{1}{\sqrt{\Delta(|k|)}}\leq& \inf\limits_{0=t_0<t_1<t_2<\cdots}\{1+\sum\limits^{\infty}_{\nu=0}\frac{V_n(t_{\nu+1})-V_n(t_{\nu})}{\sqrt{\Delta(t_{\nu})}}\}\\
  \leq&1+\int\limits^{\infty}_{0}\frac{dV_n(t)}{\sqrt{\Delta(t)}}=\int\limits^{\infty}_{0}V_n(t)\frac{d\log\sqrt{\Delta(t)}}{\sqrt{\Delta(t)}}\\
\end{split}
\end{equation*}
by partial integration.
From the proof of Lemma 8.3 in  \cite{Poschel89},
$$V_n(t)\leq \left(
               \begin{array}{c}
                 n+t \\
                 n \\
               \end{array}
             \right),
$$
this prove the lemma.
\end{proof}

\begin{lemma}
There are approximation  functions $\Delta$ such that
$$ \sum\limits_{k\in \mathbb{Z}^n}\frac{1}{\sqrt{\Delta(|k|)}}\leq K^{n\log\log n}$$
for all sufficiently large $n$ with some constant $K.$
\end{lemma}
\begin{proof}
For $t\leq n, $
$$ \left(
               \begin{array}{c}
                 n+t \\
                 n \\
               \end{array}
             \right)
 \leq
   \left(
               \begin{array}{c}
                 2n \\
                 n \\
               \end{array}
             \right)
 \leq\frac{(2n)!}{(n!)^2} \leq 4^n
$$
for all $n\geq1.$
Hence
$$\int\limits^{n}_{0}\left(
                                                                                                \begin{array}{c}
                                                                                                  n+t \\
                                                                                                  n \\
                                                                                                \end{array}
                                                                                              \right)
\frac{d\log\sqrt{\Delta(t)}}{\sqrt{\Delta(t)}}dt
\leq
4^n \int\limits^{\infty}_{0}  \frac{d\log\sqrt{\Delta(t)}}{\sqrt{\Delta(t)}}dt=4^n
$$
for every approximation  function $\Delta.$

For $t\geq n, $
$$ \left(
               \begin{array}{c}
                 n+t \\
                 n \\
               \end{array}
             \right)
 =
   \frac{1}{n!}(t+1)\cdots (t+n)
 \leq\frac{2^n}{n!} t^n.
$$
Let $\varphi$ be given by $\varphi(s)=\log^2s,$ and define $\Delta$ by
stipulating that $t\mapsto s=\log \sqrt{\Delta(t)}$ is the inverse function of $s\mapsto t=s\varphi(s),$
at least for large $t$ and $s$ respectively.
Let $s_n=n/\varphi(n)$
Since $$s\varphi(s)|_{s_n}=\frac{n}{\varphi(n)} \varphi (\frac{n}{\varphi(n)})\leq n$$
by the monotonicity of $\varphi,$ then
\begin{equation}\label{cvf}
 \begin{split}
 \int\limits^{n}_{0}\left(
                                                                                                \begin{array}{c}
                                                                                                  n+t \\
                                                                                                  n \\
                                                                                                \end{array}
                                                                                              \right)
\frac{d\log\sqrt{\Delta(t)}}{\sqrt{\Delta(t)}}dt
\leq&
\frac{2^n}{n!}\int\limits^{\infty}_{0}  t^n  \frac{d\log\sqrt{\Delta(t)}}{\sqrt{\Delta(t)}}dt\\
\leq&  \frac{2^n}{n!}\int\limits^{\infty}_{s_n} s^n\varphi^n(s)\me^{-s}ds.
\end{split}
\end{equation}
For all large $n$ and $s\geq s_n,$
$$\varphi(s)=\log^2s\leq s^{h_n},\,\,h_n=\frac{\log \varphi(s_n)}{\log s_n}\leq \frac{4\log\log n}{\log n}.$$
Thus, for all large $n,$
\begin{equation*}
 \begin{split}
 \int\limits^{n}_{0}\left(
                                                                                                \begin{array}{c}
                                                                                                  n+t \\
                                                                                                  n \\
                                                                                                \end{array}
                                                                                              \right)
\frac{d\log\sqrt{\Delta(t)}}{\sqrt{\Delta(t)}}dt
\leq&
\frac{2^n}{n!}\int\limits^{\infty}_{0} s^{n+nh_n}\me^{-s}ds\\
\leq&  \frac{2^n}{n^n}(n+nh_n)^{n+nh_n+1}=2^nA^n_nn^{nh_n+1}
\end{split}
\end{equation*}
here $A_n=(1+h_n)^{1+h_n+1/n}.$
The final estimate follows, since $A_n\rightarrow 1$ as $n\rightarrow\infty$
and $nh_n \log n= 4n \log\log n.$

\end{proof}


\subsection{ Proof of Proposition \ref{Malg}}\label{pfMalg}
\begin{proof}
We only give the proof for the  estimate of  $(AB)^{(11)}_{ij}$ and $(AB)^{(12)}_{ij}$, the  proofs for the  estimates of $(AB)^{(21)}_{ij}$ and $(AB)^{(22)}_{ij}$ are similar.

By the matrix multiplication, we have
 $$(AB)^{(11)}_{ij}=\sum_{k\in \mathbb{Z}}(A^{11}_{ik}B^{11}_{kj}+A^{12}_{ik}B^{21}_{kj})$$
and
$$(AB)^{(12)}_{ij}=\sum_{k\in \mathbb{Z}}(A^{11}_{ik}B^{12}_{kj}+A^{12}_{ik}B^{22}_{kj}).$$

$\bullet$Verifying the property (\textbf{T1}$^{\prime}$).  In view of   $A,B\in \mathfrak{M}^{\rho}_r$ and  the inequality in Lemma \ref{expineq},  we have
\begin{equation*}
\begin{split}
\|(AB)^{(11)}_{i,j}\|_{D(r)\times\mathcal{O}}
\leq&\sum\limits_{k\in \mathbb{Z}}\|A^{11}_{i,k}\|_{D(r)\times\mathcal{O}}\|B^{11}_{k,j}\|_{D(r)\times\mathcal{O}} +\sum\limits_{k}\|A^{12}_{i,k}\|_{D(r)\times\mathcal{O}}\|B^{21}_{k,j}\|_{D(r)\times\mathcal{O}}\\
\leq&  \langle\langle A\rangle\rangle_{\rho,r}\langle\langle B\rangle\rangle_{\rho,r}(\sum\limits_{k} \me^{-\rho(|i-k|+|k-j|)}
+\sum\limits_{k\in \mathbb{Z}} \me^{-\rho(|i+k|+|k+j|)})\\
\leq&  C\delta^{-1}\langle\langle A\rangle\rangle_{\rho,r}\langle\langle B\rangle\rangle_{\rho,r}\me^{-(\rho-\delta)(|i-j|)}
 \end{split}
\end{equation*}
and
\begin{equation*}
\begin{split}
\|(AB)^{(12)}_{i,j}\|_{D(r)\times\mathcal{O}}
\leq&\sum\limits_{k\in \mathbb{Z}}\|A^{11}_{i,k}\|_{D(r)\times\mathcal{O}}\|B^{12}_{k,j}\|_{D(r)\times\mathcal{O}} +\sum\limits_{k}\|A^{12}_{i,k}\|_{D(r)\times\mathcal{O}}\|B^{22}_{k,j}\|_{D(r)\times\mathcal{O}}\\
\leq&  \langle\langle A\rangle\rangle_{\rho,r}\langle\langle B\rangle\rangle_{\rho,r}(\sum\limits_{k} \me^{-\rho(|i-k|+|k+j|)}
+\sum\limits_{k\in \mathbb{Z}} \me^{-\rho(|i+k|+|k-j|)})\\
\leq&  C\delta^{-1}\langle\langle A\rangle\rangle_{\rho,r}\langle\langle B\rangle\rangle_{\rho,r}\me^{-(\rho-\delta)(|i+j|)}.
 \end{split}
\end{equation*}

$\bullet$Verifying  the property (\textbf{T2}$^{\prime}$).
In view of   $A, B\in \mathfrak{M}^{\rho}_r,$  then  following the verification of  Property (\textbf{T1}$^{\prime}$), we have
\begin{equation*}
\begin{split}
\|\lim_{t\rightarrow\infty}(AB)^{(11)}_{i+t,j+t}\|_{D(r)\times\mathcal{O}}
\leq&\sum\limits_{k}\|\lim_{t\rightarrow\infty}A^{11}_{i+t,k+t}\|_{D(r)\times\mathcal{O}}\|\lim_{t\rightarrow\infty}B^{11}_{k+t,j+t}\|_{D(r)\times\mathcal{O}}\\
& +\sum\limits_{k}\|\lim_{t\rightarrow\infty}A^{12}_{i+t,k-t}\|_{D(r)\times\mathcal{O}}\|\lim_{t\rightarrow\infty}B^{21}_{k-t,j+t}\|_{D(r)\times\mathcal{O}}\\
\leq&  \langle\langle A\rangle\rangle_{\rho,r}\langle\langle B\rangle\rangle_{\rho,r}(\sum\limits_{k} \me^{-\rho(|i-k|+|k-j|)}
+\sum\limits_{k} \me^{-\rho(|i+k|+|k+j|)})\\
\leq&  C\delta^{-1}\langle\langle A\rangle\rangle_{\rho,r}\langle\langle B\rangle\rangle_{\rho,r}\me^{-(\rho-\delta)(|i-j|)}
 \end{split}
\end{equation*}
and
\begin{equation*}
\begin{split}
\|\lim_{t\rightarrow\infty}(AB)^{(12)}_{i+t,j-t}\|_{D(r)\times\mathcal{O}}
\leq&\sum\limits_{k}\|\lim_{t\rightarrow\infty}A^{11}_{i+t,k+t}\|_{D(r)\times\mathcal{O}}\|\lim_{t\rightarrow\infty}B^{12}_{k+t,j-t}\|_{D(r)\times\mathcal{O}}\\
& +\sum\limits_{k}\|\lim_{t\rightarrow\infty}A^{12}_{i+t,k-t}\|_{D(r)\times\mathcal{O}}\|\lim_{t\rightarrow\infty}B^{22}_{k-t,j-t}\|_{D(r)\times\mathcal{O}}\\
\leq&  \langle\langle A\rangle\rangle_{\rho,r}\langle\langle B\rangle\rangle_{\rho,r}(\sum\limits_{k} \me^{-\rho(|i-k|+|k+j|)}
+\sum\limits_{k} \me^{-\rho(|i+k|+|k-j|)})\\
\leq&  C\delta^{-1}\langle\langle A\rangle\rangle_{\rho,r}\langle\langle B\rangle\rangle_{\rho,r}\me^{-(\rho-\delta)(|i-j|)}.
 \end{split}
\end{equation*}
These imply the property (T2) holds.

$\bullet$Verifying   the property (\textbf{T3}$^{\prime}$).
Denote
 $ A^{11}_{i,j,\infty}:=\lim_{t\rightarrow\infty}A^{11}_{i+t,j+t}$ and
 $ A^{12}_{i,j,\infty}:=\lim\limits_{t\rightarrow\infty} A^{12}_{i+t,j-t}.$
Similarly for other terms.

Then by the difference equality \eqref{difference} and the inequality in Lemma \ref{expineq}, we have
\begin{equation*}
\begin{split}
&\|(AB)^{(11)}_{i+t,j+t}-\lim_{t\rightarrow\infty}(AB)^{(11)}_{i+t,j+t}\|_{D(r)\times\mathcal{O}}\\
\leq&\sum\limits_{k}\|A^{11}_{i+t,k+t}-A^{11}_{i,k,\infty}\|_{D(r)\times\mathcal{O}}\|B^{11}_{k,j,\infty}\|_{D(r)\times\mathcal{O}}\\
& +\sum\limits_{k}\|A^{11}_{i,k,\infty}\|_{D(r)\times\mathcal{O}}\|B^{11}_{k+t,j+t}-B^{11}_{k,j,\infty}\|_{D(r)\times\mathcal{O}}\\
& +\sum\limits_{k}\|A^{11}_{i+t,k+t}-A^{11}_{i,k,\infty}\|_{D(r)\times\mathcal{O}} \|B^{11}_{k+t,j+t}-B^{11}_{k,j,\infty}\|_{D(r)\times\mathcal{O}}\\
 \end{split}
\end{equation*}
\begin{equation*}
\begin{split}
&+\sum\limits_{k}\|A^{12}_{i+t,k-t}-A^{12}_{i,k,\infty}\|_{D(r)\times\mathcal{O}}\|B^{21}_{k,j,\infty}\|_{D(r)\times\mathcal{O}}\\
& +\sum\limits_{k}\|A^{12}_{i,k,\infty}\|_{D(r)\times\mathcal{O}}\|B^{21}_{k-t,j+t}-B^{21}_{k,j,\infty}\|_{D(r)\times\mathcal{O}}\\
& +\sum\limits_{k}\|A^{12}_{i+t,k-t}-A^{12}_{i,k,\infty}\|_{D(r)\times\mathcal{O}} \|B^{21}_{k-t,j+t}-B^{21}_{k,j,\infty}\|_{D(r)\times\mathcal{O}}\\
 \end{split}
\end{equation*}
\begin{equation*}
\begin{split}
\leq& |t|^{-1} \langle\langle A\rangle\rangle_{\rho,r}\langle\langle B\rangle\rangle_{\rho,r}(\sum\limits_{k} \me^{-\rho(|i-k|+|k-j|)}
+\sum\limits_{k} \me^{-\rho(|i+k|+|k+j|)})\\
\leq&  |t|^{-1}C\delta^{-1}\langle\langle A\rangle\rangle_{\rho,r}\langle\langle B\rangle\rangle_{\rho,r}\me^{-(\rho-\delta)(|i-j|)}
 \end{split}
\end{equation*}
and
\begin{equation*}
\begin{split}
&\|(AB)^{(12)}_{i+t,j-t}-\lim_{t\rightarrow\infty}(AB)^{(12)}_{i+t,j-t}\|_{D(r)\times\mathcal{O}}\\
\leq&\sum\limits_{k}\|A^{11}_{i+t,k+t}-A^{11}_{i,k,\infty}\|_{D(r)\times\mathcal{O}}\|B^{12}_{k,j,\infty}\|_{D(r)\times\mathcal{O}}\\
&+\sum\limits_{k}\|A^{11}_{i,k,\infty}\|_{D(r)\times\mathcal{O}}\|B^{12}_{k+t,j-t}-B^{12}_{k,j,\infty}\|_{D(r)\times\mathcal{O}}\\
&+\sum\limits_{k}\|A^{11}_{i+t,k+t}-A^{11}_{i,k,\infty}\|_{D(r)\times\mathcal{O}}\|B^{12}_{k+t,j-t}-B^{12}_{k,j,\infty}\|_{D(r)\times\mathcal{O}}\\
& +\sum\limits_{k}\|A^{12}_{i+t,k-t}-A^{12}_{i,k,\infty}\|_{D(r)\times\mathcal{O}}\|B^{22}_{k,j,\infty}\|_{D(r)\times\mathcal{O}}\\
& +\sum\limits_{k}\|A^{12}_{i,k,\infty}\|_{D(r)\times\mathcal{O}}\|B^{22}_{k-t,j-t}-B^{22}_{k,j,\infty}\|_{D(r)\times\mathcal{O}}\\
& +\sum\limits_{k}\|A^{12}_{i+t,k-t}-A^{12}_{i,k,\infty}\|_{D(r)\times\mathcal{O}}\|B^{22}_{k-t,j-t}-B^{22}_{k,j,\infty}\|_{D(r)\times\mathcal{O}}\\
\leq&|t|^{-1}  \langle\langle A\rangle\rangle_{\rho,r}\langle\langle B\rangle\rangle_{\rho,r}(\sum\limits_{k} \me^{-\rho(|i-k|+|k+j|)}
+\sum\limits_{k} \me^{-\rho(|i+k|+|k-j|)})\\
\leq& |t|^{-1}C\delta^{-1}\langle\langle A\rangle\rangle_{\rho,r}\langle\langle B\rangle\rangle_{\rho,r}\me^{-(\rho-\delta)(|i+j|)}.
 \end{split}
\end{equation*}

\end{proof}

\subsection{ Some Technical Lemmas}

 \begin{lemma}\label{expineq}
\begin{equation}
   \sum\limits_{k\in \mathbb{Z}}\me^{-\delta(|i-k|+|k-j|)}\leq C\delta^{-1} ,
\end{equation}
where $C$ is a positive constant that depends on $\beta$ and does not depend on $i,j.$
\end{lemma}

\begin{lemma}\label{measest}\cite{Russmann}
 Let $f:[a,b]\rightarrow \mathbb{R}$ be a $q-$times continuously differentiable function satisfying
 $$|f^{(q)}(t)|\geq \beta,\,\forall\,\, t\in [a,b]$$
 for some $q\in \mathbb{N}$ and  $\beta>0.$ Then we have the estimate
 $$\mes\{ t\in [a,b]: |f(t)|\leq \varepsilon  \}\leq   4\left(  \frac{q!}{2\beta}\varepsilon\right)^{\frac{1}{q}},\,\,\forall\,\,\varepsilon>0.$$
\end{lemma}

%

\end{document}